\documentclass[12pt,leqno]{article}
\usepackage{amsmath,amsfonts}
\usepackage[dvips]{graphics,color}
\usepackage{latexsym}
\usepackage{amssymb}
\usepackage{graphicx}
\def\a{\mathcal{A}}
\def\b{\mathcal{B}}
\def\c{\mathcal{C}}
\def\d{\mathcal{D}}
\def\e{\mathcal{E}}
\def\f{\mathcal{F}}
\def\g{\mathcal{G}}
\def\h{\mathcal{H}}

\def\j{\mathcal{J}} 
\def\k{\mathcal{K}} 
\def\l{\mathcal{L}}
\def\m{\mathcal{M}}
\def\n{\mathcal{N}}
\def\p{\mathcal{P}}
\def\q{\mathcal{Q}}
\def\r{\mathcal{R}}
\def\s{\mathcal{S}}
\def\t{\mathcal{T}}
 \def\u{\mathcal{U}}
\def\v{\mathcal{V}}
\def\w{\mathcal{W}}

\def\y{\mathcal{Y}}
\def\z{\mathcal{Z}}

\def\sp{(s,\psi)}
\def\ss{(\rho,\psi)}
\def\bp{\bar{\psi}}
\def\pc{\chi\psi}
\def\tp{\frac{1}{2\pi i}}
\def\al{\alpha}
\def\ep{\varepsilon}

\def\la{\lambda}

\def\tg{\tilde{g}}

\def\tz{\tilde{Z}}
\def\vk{\varkappa}
\def\bt{\bar{\theta}}

\def\be{\beta}
\def\tk{\tilde{\kappa}}

\def\aa{\mathfrak{a}}

\def\ff{\mathfrak{f}}
\def\gh{\mathfrak{g}}

\def\yy{\mathfrak{y}}
\date{}
\begin{document}
\title{Discrete mean estimates and the Landau-Siegel zero}
\maketitle
\thanks{}
\author{\centerline {Yitang Zhang}}
\maketitle

\section *\small{{\centerline{Abstract} }
\medskip\noindent

 \small{Let $\chi$ be a real primitive character to the modulus $D$. It is proved that}
$$
L(1,\chi)\gg (\log D)^{-2022}
$$
where the implied constant is  absolute and effectively computable.

In the proof, the lower bound for $L(1,\chi)$ is first related to the distribution of zeros of a family of Dirichlet $L$-functions in a certain region, and some results on the gaps between consecutive zeros are derived. Then, by  evaluating certain  discrete means of the large sieve type, a contradiction can be obtained if $L(1,\chi)$ is too small.}

\section *{\centerline{  Table of Contents} }

\medskip\noindent
1. Introduction                                                                                       

\medskip\noindent
 2. Notation and outline of the proof
 
\medskip\noindent
  3. The set $\Psi_1$

\medskip\noindent
4. Zeros of $L\sp L(s,\pc)$ in $\Omega$

 \medskip\noindent
 5. Some analytic lemmas
 
  \medskip\noindent
 6. Approximate formula for $L\sp$
 
  \medskip\noindent
 7. Mean value formula I
 
 \medskip\noindent
 8. Evaluation of $\Xi_{11}$

  \medskip\noindent
 9. Evaluation of $\Xi_{12}$

  \medskip\noindent
 10. Proof of Proposition 2.4

 \medskip\noindent
 11. Proof of Proposition 2.6
 
  \medskip\noindent
 12. Evaluation of $\Xi_{15}$
 
  \medskip\noindent
 13. Approximation to $\Xi_{14}$

 \medskip\noindent
 14. Mean value formula II
 
  \medskip\noindent
 15. Evaluation of $\Phi_1$

 \medskip\noindent
 16. Evaluation of $\Phi_2$ 
 
 \medskip\noindent
 17. Evaluation of $\Phi_3$ 

 \medskip\noindent
 18. Proof of Proposition 2.5

  \medskip\noindent
 Appendix A. Some Euler products

   \medskip\noindent
 Appendix B. Some arithmetic sums

\medskip\noindent
 References

\section *{\centerline{1. Introduction}}
\medskip\noindent

Let $\chi$ be a real primitive character to the modulus $D$. It is known that the Dirichlet $L$-function 
$L(s,\chi)$ has at most one real and simple zero $\tilde{\rho}$ satisfying
$$
1-\tilde{\rho}<c_0(\log D)^{-1}
$$
where $c_0>0$ is an absolute constant.  Such a zero is called the Landau-Siegel zero. 
The typical methods to determine zero-free regions for Dirichlet $L$-functions are unable to 
eliminate the Landau-Siegel zero for an intrinsic reason. 

The non-vanishing of $L(s,\chi)$ near $s=1$  is closely related to the lower bound 
for the value of $L(s,\chi)$ at $s=1$. The well-known
Siegel theorem [19] asserts that, for any $\ep>0$, there exists a positive number $C_1(\ep)$ such that
$$
L(1,\chi)>C_1(\ep)D^{-\ep}.
$$
This implies, for any $\ep>0$, that there exists a positive number $C_2(\ep)$ such
that
$$
L(\sigma,\chi)\ne 0\qquad\text{if}\qquad \sigma>1-C_2(\ep)D^{-\ep}.
$$
However, the process of Siegel's proof essentially involves certain assumptions, such as 
the existence of a real primitive character $\hat{\chi}$ (mod $\hat{D}$) with $L(s,\hat{\chi})$ 
vanishing on the interval $[1-\ep,1)$ (see  [8]),
and  the constant $C_1(\ep)$  relies on the unknown modulus $\hat{D}$ necessarily .
This makes the result of Siegel's theorem entirely ineffective.

It is known (see  [9]) that the non-existence of the Landau-Siegel zero implies
 $$L(1,\chi)\gg(\log D)^{-1}.\eqno(1.1)$$   
 In the case $\chi(-1)=-1$, 
 Goldfeld [10] and Gross and Zagier [12] proved that
$$
L(1,\chi)\gg D^{-1/2}(\log D)^{1-\ep}
$$ 
for any $\ep>0$, where the implied constant is effectively computable. This result surpasses the trivial bound
$L(1,\chi)\gg D^{-1/2}$ and solves Gauss' class number problem for imaginary quadratic fields.
Granville and Stark [11] proved that the uniform {\it abc} conjecture for number fields implies (1.1) when $\chi(-1)=-1$.

The main result of this paper is 
\medskip\noindent

 {\bf Theorem 1}\,\, {\it If $\chi$ is a real primitive character  to the modulus $D$, then
$$
L(1,\chi)>c_1(\log D)^{-2022}
$$
where $c_1>0$ is an absolute, effectively computable constant.}
\medskip\noindent

As a direct consequence of Theorem 1 we have 
\medskip\noindent

{\bf Theorem 2}\,\, {\it If $\chi$ is a real primitive character  to the modulus $D$, then
$$
L(\sigma,\chi)\ne 0 
$$
for
$$
\sigma>1-c_2(\log D)^{-2024}
$$
where $c_2>0$ is an absolute, effectively computable constant.} 
\medskip\noindent

It is   possible to replace the exponent $-2022$ in Theorem 1 by a larger (negative) value if the  current arguments are modified,  but we  will not  discuss it in this paper. On the other hand, it seems that the lower bound (1.1) can not be achieved by the present methods.

\medskip\noindent

{\it Acknowledgements.}\,\, The basic ideas of the present proof were initially formed during my visit to  the Institute for Advanced Study in the spring of 2014. I thank the Institute for Advanced Study for providing me with excellent conditions. I also thank Professor Peter  Sarnak for his encouragement.

\section *{\centerline{2. Notation and outline of the proof}}

\medskip\noindent
{\it   Notation and conventions.}
\medskip\noindent

Throughout, $t$, $u$, $v$,  $z$ and $\sigma$ denote real variables;   $x$ and $y$ denote positive real variables; $s$ and $w$ denote  the complex variables
$s=\sigma+it$ and  $w=u+iv$; $d$, $h$, $j$,  $k$, $l$, $m$, $n$ and $r$
denote natural numbers, while $p$ and $q$ denote  primes;
 the arithmetic functions  $\mu(n)$ and  $\varphi(n)$ are defined as usual, while
$\tau_k(n)$ denotes the $k$-fold divisor function. Let $*$ denote the Dirichlet convolution of arithmetic functions.  Let  $e(s)=\exp\{2\pi is\}$. We  write
$$
\sum_n\qquad\text{and}\qquad\sum_{(n,k)=1}
$$
for
$$
\sum_{n=1}^{\infty}\qquad\text{and}\qquad\sum_{\substack{n=1\\(n,k)=1}}^{\infty}
$$
respectively. Let $\int_{(z)}$ denote  an integration from $z-i\infty$ to $z+i\infty$.

We let $\chi$ denote a real primitive character to the modulus $D$ with $D$ greater than a sufficiently large and 
effectively computable number. Write
$$
\l=\log D.\eqno(2.1)
$$

 Let $c$  denote a positive, effectively computable and 
absolute constant, not necessarily the same at each occurrence.  All implied constants, unless specified, are absolute and effectively computable. 
In this paper, an error term is usually given in the form that is small enough for our purpose, although shaper estimates are possible.

 For any  Dirichlet character  $\theta$ to the modulus $k$,  let $\tau(\theta)$ denote the Gauss sum
$$
\tau(\theta)=\sum_{a(\bmod\, k)}\theta(a)e(a/k)
$$  
(thus we shall  not write $\tau$ for $\tau_2$, the ordinary divisor function). 
The expressions 
$$
\sum_{\theta(\bmod\, k)},\qquad \sideset{}{'}\sum_{\theta(\bmod \,k)}\qquad\text{and}\qquad \sideset{}{^*}\sum_{\theta(\bmod\, k)}
$$
denote, respectively, a sum over all ${\theta(\bmod \,k)}$, a sum over all non-principal ${\theta(\bmod\, k)}$,
and a sum over all primitive ${\theta(\bmod \,k)}$.  Let $\psi_k^0$ denote the principle character $(\bmod\, k)$.

 In case   $\theta(\bmod\, k)$ is  a primitive character, the functional equation for $L(s,\theta)$ is
$$
L(s,\theta)=Z(s,\theta)L(1-s,\bar{\theta}), \eqno (2.2)
$$
where
$$
Z(s,\theta)=\tau(\theta)\pi^{s-1/2}k^{-s}\,
\Gamma\bigg(\frac{1-s}{2}\bigg)\Gamma\bigg(\frac{s}{2}\bigg)^{-1}$$
if $\theta(-1)=1$, and
$$
Z(s,\theta)=-i\tau(\theta)\pi^{s-1/2}k^{-s}\,
\Gamma\bigg(\frac{2-s}{2}\bigg)\Gamma\bigg(\frac{1+s}{2}\bigg)^{-1}$$
if $\theta(-1)=-1$. When $t$ is large,
it is convenient to use an asymptotic expression for $Z(s,\theta)$ as follows.  Let
$$
\vartheta(s)=2(2\pi)^{s-1}\Gamma(1-s)\sin(\pi s/2)\eqno(2.3)
$$
(this function is usually written as $\chi(s)$ in the literature). It is known that 
$$
\pi^{s-1/2}\Gamma\bigg(\frac{1-s}{2}\bigg)\Gamma\bigg(\frac{s}{2}\bigg)^{-1}=\vartheta(s)
$$
and 
$$
\pi^{s-1/2}\Gamma\bigg(\frac{2-s}{2}\bigg)\Gamma\bigg(\frac{1+s}{2}\bigg)^{-1}=\vartheta(s)\cot(\pi s/2).
$$
Assume $t>1$. Then
$$
i\cot(\pi s/2)=1+O(e^{-\pi t}).
$$
Thus we have uniformly
$$
Z(s,\theta)=\theta(-1)\tau(\theta)k^{-s}\vartheta(s)(1+O(e^{-\pi t})),\eqno(2.4)
$$
and
 $$
Z(s,\theta)^{-1}=\tau(\bar{\theta})k^{s-1}\vartheta(1-s)(1+O(e^{-\pi t})),\eqno(2.5)
$$ 
the $O(e^{-\pi t})$ term being identically zero in case $\theta(-1)=1$.

\medskip\noindent
{\it  Outline of the proof: Initial steps} 
\medskip\noindent

Our aim is to derive a contradiction from the following 
\medskip\noindent

{\bf Assumption (A)}
$$
L(1,\chi)<\l^{-2022}. 
$$
\medskip\noindent

The underlying ideas are inspired by the work of Goldfeld [10], Iwaniec and Sarnak [15] and
Conrey, Ghosh and Gonek [4]. 
The  paper [10] exhibits  a relationship between
the lower bound for $L(1,\chi)$ and the order of zeros of the function $L_E(s)L_E(s,\chi)$ at
 the central point, where $E$ is an elliptic curve;  the  paper [15] exhibits a relationship between the lower 
 bound for $L(1,\chi)$ and the non-vanishing of central values of a family of automorphic $L$-functions; 
 the paper [4] has definite influence on the idea of reducing the problem to evaluating 
 certain discrete means that is employed in the present work (the reader is also referred to [3]).

 The initial part of this paper consists in relating the lower bound for $L(1,\chi)$ to the distribution 
 of zeros of the function $L\sp L(s,\psi\chi)$, with $\psi$ belonging to a large 
 set $\Psi$ of primitive characters, in a domain $\Omega$.  Thus we introduce
$$
P=\exp\{\l^{9}\}.\eqno(2.6)
$$
For notational simplicity we write
$$
p\sim P\quad\text{for}\quad P<p<P(1+\l^{-68}).
$$
Let $\Psi$ be the set of all primitive characters $\psi(\bmod\,p)$ with $p\sim P$, and let
$$
\Omega=\big\{s:\quad |\Re(s-s_0)|<1/2,\quad |\Im(s-s_0)|<\l_1+2\big\}\eqno(2.7)
$$
where
$$
\l_1=\l^{405},\quad s_0=\frac12+2\pi it_0\quad\text{with}\quad t_0=\l^{519}.\eqno(2.8)
$$
Throughout, we will assume that $p\sim P$ and $\psi$ is a primitive characters $(\bmod\,p)$, i.e., $\psi\in\Psi$. 
Note that
$$
\p:=\sum_{p\sim P}p=(1+o(1))P^2\l^{-77}.\eqno(2.9)
$$
  In Section 3 we shall introduce a subset $\Psi_1$ of $\Psi$  and prove 
\medskip\noindent

{\bf Proposition 2.1.}\,\, {\it Let $\Psi_2$ be the complement of $\Psi_1$ in $\Psi$. If (A) holds, then}
$$
\sum_{\psi\in\Psi_2}1\ll \p\l^{-739}.
$$
\medskip\noindent

For $\psi(\bmod\,p)\in \Psi$, the average gap between consecutive zeros of $L(s,\psi)L(s,\psi\chi)$ in $\Omega$
is
$$
\simeq\frac{2\pi}{\log(Dp^2t_0^2)}=\al(1+O(\al\l)).
$$
where
$$
\al=\frac{\pi}{\log P}.\eqno(2.10)
$$
In Section 4 we  shall prove
\medskip\noindent

{\bf Proposition 2.2.}\,\, {\it Suppose that $\psi\in\Psi_1$. Then
 the following hold.} 

\medskip\noindent

(i). {\it \,\,All the zeros of  $L(s,\psi)L(s,\psi\chi)$ in $\Omega$ lie on the critical line.}

\medskip\noindent

(ii). {\it \,All the zeros of $L(s,\psi)L(s,\psi\chi)$ in $\Omega$ are simple.}

\medskip\noindent

(iii). {\it The gap between any consecutive zeros of $L(s,\psi)L(s,\psi\chi)$ in $\Omega$ is of the form}
$$\al+O(\al^2\l).$$
\medskip\noindent

It should be stressed that the set $\Psi_1$ is formally defined. In fact, Proposition 2.2 and some other results for $\psi\in\Psi_1$ are unconditionally derived from the definition of $\Psi_1$. However,   one is even unable
to determine whether $\Psi_1$ is empty or not without assuming (A).

The results of Proposition 2.1 and 2.2 are not surprising and they go back to Heath-Brown [15]. One may believe that the assertions (i) and
 (ii) of Proposition 2.2  hold for all characters $\psi$ in $\Psi$. On the other hand, however,
   in comparison with some conjectures on the vertical distribution of zeros of the 
  Riemann zeta function (see [2] and [17]), one may claim that 
  the gap assertion (iii) of Proposition 2.2 fails to hold 
  for most of $\psi$ in $\Psi$. Our proof thus consists in deriving 
  a contradiction  from the gap assertion.

Assume $\psi\in\Psi$. Recall that the functional equation for $L\sp$ is given by (2.2) with $\theta=\psi$.
On the upper-half plane, since $Z\sp$ is analytic and non-vanishing, there is an analytic function $Y\sp$ such that
 $$
 Y\sp^2=Z\sp^{-1}.
 $$
 Let 
$$M\sp=Y\sp L\sp.
$$
Note that $M\sp$ is defined for $t>0$ only, and, for each $\psi$, there are two choices of $M\sp$ up to $\pm$,  but the expressions
$$
M(s_1,\psi)M(s_2,\psi)\qquad\text{and}\qquad \frac{M(s_1,\psi)}{M(s_2,\psi)}
$$
are independent of these choices. The functional equation (2.2)  gives
$$
M\sp=Y\sp^{-1}L(1-s,\bp).
$$
 Since $|Y(1/2+it,\psi)|=1$, it follows  that 
 $$
 |M(1/2+it,\psi)|=|L(1/2+it,\psi)|,
 $$
$$
M(1/2+it,\psi)\in\mathbf{R},\eqno(2.11)
$$
and, consequently,
$$
iM'(1/2+it,\psi)\in\mathbf{R}.\eqno(2.12)
$$

The gap assertion (iii) of Proposition 2.2  can be restated as follows. If  $\psi\in\Psi_1$ and $(1/2+i\gamma, 1/2+i\gamma')$ is a pair of consecutive zeros of $L\sp L(s,\pc)$ in $\Omega$ with  $\gamma'>\gamma$, then for some (large) constant $c'>0$,
$$
|\gamma'-\gamma-\al|<c'\al^2\l.
$$
Write
$$
\be_1=i\al(1-5c'\al\l),\qquad\be_2=2 i\al(1+c'\al\l),\qquad \be_3=3 i\al(1-c'\al\l).\eqno(2.13)
$$
In what follows we  assume $\psi\in\Psi_1$. Let $\z(\psi)$ denote the set of zeros of $L\sp$ in the region $$|\sigma-1/2|<1/2,\qquad |t-2\pi t_0|<\l_1\eqno(2.14)$$ 
(this is slightly smaller than $\Omega$).
Assume  $\rho\in\z(\psi)$. Note that  $$M'\ss=Y\ss L'\ss\ne 0.$$ Write
 $$
 \c^*\ss=\frac{-iM(\rho+\be_1,\psi)M(\rho+\be_2,\psi)M(\rho+\be_3,\psi)}{M'\ss}.
 $$
By (2.11) and (2.12) we have
$$
\c^*\ss\in\mathbf{R}.
$$
\medskip\noindent

It will be proved that
\medskip\noindent

{\bf Lemma 2.3.}\,\,{\it Suppose  $\psi\in\Psi_1$ and $\rho\in\z(\psi)$. Then}
$$
\c^*\ss\ge 0.
$$

We introduce the smooth weight
 $$
\omega(s)=\frac{\sqrt{\pi}}{\l_2}\exp\bigg\{\frac{(s-s_0)^2}{4\l_2^2}\bigg\}\quad\text{with}\quad\l_2=\l^{400},\eqno(2.15)
$$
which is positive for $\sigma=1/2$. Lemma 2.3 implies that
$$
\sum_{\psi\in\Psi_1}\,\sum_{\rho\in\z(\psi)}\c^*\ss |\h\ss|^2\omega(\rho)\ge 0.
 \eqno(2.16)
 $$
for all functions $\h\sp$ defined on $\Omega$. Thus, in order to derive a contradiction from (A), it suffices to find  certain  functions $\h\sp$  and show that , under  Assumption  (A), (2.16) fails to hold.

\medskip\noindent
{\it   A heuristic argument}
\medskip\noindent

It is natural to consider the function $\h\sp$ which is of the form
$$
\sum_{n<P'}\frac{f(\log n/\log P')\mu\psi(n)}{n^s}
$$
for some  $P'<P$, where $f$ is a polynomial with $f(1)=0$.  It should be stressed, on assuming (A), that
  if $\h\sp$ is of the above form  with $P'$ slightly smaller than $P$ in the logarithmic scale,  the sum in (2.16) can be evaluated. This is analogous to the results of Conrey, Iwaniec and Soundararajan [5] and [6]. However, it seems that  such a choice of $\h\sp$ is not good enough for our purpose.  On the other hand, the left side of (2.16) could be ``close" to zero for certain non-zero functions $\h\sp$. For example, one may consider the choice that $\h\sp$ is a linear combination of $L(s+\be_j,\pc)$, $1\le j\le 3$. Thus, inspired by  the approximate formula for $L(s,\pc)$,  some new forms of $\h\sp$ could be introduced. However, our goal is not to show that (2.16) fails to hold for certain $\h\sp$ directly, instead, a variant of the argument will be adopted.
  
  The crucial part of the proof is to construct functions $H_1\sp$, $H_2\sp$, $J_1\sp$ and $J_2\sp$, each of which is of the form
  $$
\sum_{n<P'}\frac{f(\log n/\log P)\pc(n)}{n^s}
$$
with $\log P'/\log P$ close to $1/2$, and show that the left side of (2.16) is "small" if 
$$
\h\sp=H_1\sp+Z(s,\pc)\overline{H_2\sp}\quad\text{or}\quad \h\sp=J_1\sp+Z(s,\pc)\overline{J_2\sp}
$$
for $\sigma=1/2$. On the other hand, the sum
$$
\Xi_1^*:=\sum_{\psi\in\Psi_1}\,\sum_{\rho\in\z(\psi)}\c^*\ss\big( H_1\ss\overline{J_1\ss}+\overline{H_2\ss}J_2\ss\big)\omega(\rho)\eqno(2.17)
$$
is not too small in absolute value.

For $\psi\in\Psi_1$ and $\rho\in\tz(\psi)$, since
$$\aligned
 H_1\ss&\overline{J_1\ss}+\overline{H_2\ss}J_2\ss\\=&\big( H_1\ss+Z(\rho,\pc)\overline{H_2\ss}\big)\overline{J_1\ss}
 -\big(Z(\rho,\pc)\overline{J_1\ss}-J_2\ss\big)\overline{H_2\ss}
\endaligned
$$
and
$$
\big|Z(\rho,\pc)\overline{J_1\ss}-J_2\ss\big|=\big|J_1\ss-Z(\rho,\pc)\overline{J_2\ss}\big|,
$$
it follows that
$$\aligned
\big| H_1&\ss\overline{J_1\ss}+\overline{H_2\ss}J_2\ss\big|\\\le&\big| H_1\ss+Z(\rho,\pc)\overline{H_2\ss}\big||J_1\ss|
 +\big|J_1\ss-Z(\rho,\pc)\overline{J_2\ss}\big||H_2\ss|.
\endaligned
$$
This yields, by Lemma 2.3,
$$
|\Xi_1^*|\le \Xi_2^*+\Xi_3^*\eqno(2.18)
$$
where
$$\Xi_2^*=\sum_{\psi\in\Psi_1}\,\sum_{\rho\in\z(\psi)}\c^*\ss\big| H_1\ss+Z(\rho,\pc)\overline{H_2\ss}\big||J_1\ss|\omega(\rho),\eqno(2.19)
$$
$$
\Xi_3^*=\sum_{\psi\in\Psi_1}\,\sum_{\rho\in\z(\psi)}\c^*\ss\big|J_1\ss-Z(\rho,\pc)\overline{J_2\ss}\big||H_2\ss|\omega(\rho).
\eqno(2.20)$$

 \medskip\noindent
 {\it   Discrete mean estimates}
\medskip\noindent

 We  introduce some parameters and functions as follows. Let
$$
P_1=P^{0.504},\qquad P_2=P^{0.5}T^{-10},\qquad P_3=P^{0.498},\eqno(2.21)
$$
$$
\be_6=\frac{3i\al }{2},\qquad \be_7=\frac{5i\al }{2}.\eqno(2.22)
$$
(The definition of $\be_4$ and $\be_5$ will be given in Section 8.) Write
$$
H_{11}\sp=\sum_{n<P_1}\frac{\pc(n)}{n^s}\bigg(1-\frac{\log n}{\log P_1}\bigg)\bigg(\frac{P_1}{n}\bigg)^{\be_6},\eqno(2.23)
$$
$$
H_{12}\sp=\sum_{n<P_2}\frac{\pc(n)}{n^s}\bigg(1-\frac{\log n}{\log P_2}\bigg)\bigg(\frac{P_2}{n}\bigg)^{\be_7},\eqno(2.24)
$$
$$
H_{13}\sp=\sum_{n<P_3}\frac{\pc(n)}{n^s}\bigg(1-\frac{\log n}{\log P_3}\bigg)\bigg(\frac{P_3}{n}\bigg)^{\be_6}\eqno(2.25)
$$
(these sums are written as  $H_{1j}\sp$ rather than $H_{1j}(s,\pc)$). We also introduce the numerical constants
$$
\iota_2=0.94977-1.38995i,\quad\iota_3=-1.00635-0.22789i,\quad\iota_4=-0.68738+1.60688i.\eqno(2.26)
$$
We define
$$
H_1\sp=H_{11}\sp+\iota_2 H_{12}\sp,\quad H_2\sp=\bar{\iota}_3H_{13}\sp+\bar{\iota}_4 H_{12}\sp.\eqno(2.27)
$$ 
Let $\tilde{f}(z)$ be given by
$$
\tilde{f}(z)=
\begin{cases}
500(z-0.5)&\quad\text{if}\quad 0.5\le z\le 0.502\\
500(0.504-z)&\quad\text{if}\quad 0.502\le z\le 0.504\\
0&\quad\text{otherwise}.
\end{cases}\eqno(2.28)
$$
We define
$$
J_1\sp=\sum_n\frac{\pc(n)}{n^s}\tilde{f}\bigg(\frac{\log n}{\log P}\bigg),\eqno(2.29)
$$
$$
J_2\sp=\sum_n\frac{\pc(n)}{n^s}\tilde{f}\bigg(\frac{\log n}{\log P}+0.004-\tilde{\al}\bigg)\quad\text{with} \quad \tilde{\al}=\frac{\log D t_0}{\log P}.\eqno(2.30)
$$

Let
$$
\aa=\frac{6}{\pi^2}L'(1,\chi)^2\prod_{q|D}\frac{q}{q+1}.\eqno(2.31)
$$
It can be shown that the assumption (A) implies 
$$\mathfrak{a}\gg 1.
$$ 
(see Lemma 5.7).
\medskip\noindent

{\bf Proposition 2.4.}\,\, {\it
Assume that (A) holds. Then}
$$
|\Xi_1^*|>5\aa\p.
 $$
 
 \medskip\noindent
 
 {\bf Proposition 2.5.}\,\, {\it Assume that (A) holds. Then}
$$
\Xi_2^*<2\aa\p.
$$
\medskip\noindent

  {\bf Proposition 2.6.}\,\, {\it Assume that (A) holds. Then}
$$
\Xi_3^*=o(\aa\p).
$$

\medskip\noindent

Under Assumption (A), a contradiction is immediately derived from (2.18), Proposition 2.4, 2.5  and 2.6, This proves Theorem 1.

It should be remarked that Proposition 2.5 follows from the inequalities (on assuming (A))
$$\sum_{\psi\in\Psi_1}\,\sum_{\rho\in\z(\psi)}\c^*\ss\big| H_1\ss+Z(\rho,\pc)\overline{H_2\ss}\big|^2\omega(\rho)<0.001\aa\p,\eqno(2.32)
$$
and
$$\sum_{\psi\in\Psi_1}\,\sum_{\rho\in\z(\psi)}\c^*\ss |J_1\ss|^2\omega(\rho)<3000\aa\p,\eqno(2.33)
$$
by Cauchy's inequality (and Lemma 2.3). The proof of (2.32) involves some numerical calculations.

We conclude this section by proving Lemma 2.3.

\medskip\noindent

{\it  Proof of Lemma 2.3}.\,\,  By (2.) we see that $M(\rho+iv,\psi)\ne 0$ if $|\be_2|\le v\le |\be_3|$. This implies, by the mean-value theorem, that
 $$
M(\rho+\be_2,\psi)M(\rho+\be_3,\psi)>0,
$$  
since $M(1/2+it,\psi)$ is a real-valued continuous function in $t$. Similarly  we have
$$
\frac{M(\rho+\be_1,\psi)}{M(\rho+iv,\psi)}>0
$$
if $0<v\le|\be_1|$. Since
$$
\frac{M(\rho+\be_1,\psi)}{iM'(\rho,\psi)}=\lim_{v\to 0^+}\frac{vM(\rho+\be_1,\psi)}{M(\rho+iv,\psi)},
$$
it follows that
$$
\frac{M(\rho+\be_1,\psi)}{iM'(\rho,\psi)}\ge 0.
$$
 This completes the proof. $\Box$
\medskip\noindent

{\it Remark.}\,\,  It is implied in the proof of Lemma 2.3  that
$$
\bigg|\frac{L(\rho+\be_1,\psi)}{L'\ss}\bigg|=-i\frac{M(\rho+\be_1,\psi)}{M'\ss}\eqno(2.34)
$$
for $\psi\in\Psi_1$ and $\rho\in\tz(\psi)$.

 \section *{\centerline{3. The set $\Psi_1$}}
 
 \medskip\noindent
 
 Let $\nu(n)$ and $\upsilon(n)$ be given by
 $$
\zeta(s)L(s,\chi)= \sum_n\frac{\nu(n)}{n^s},\quad \zeta(s)^{-1}L(s,\chi)^{-1}= \sum_n\frac{\upsilon(n)}{n^s},\quad\sigma>1,
$$
respectively. It is easy to see that
$$
|\upsilon(n)|\le\nu(n)\le\tau_2(n).\eqno(3.1)
$$

{\bf Lemma 3.1.}\,\,  {\it Assume (A) holds. Then}
$$
\sum_{D^4<n\le P^2}\frac{\nu(n)^2}{n}\ll\l^{-2011}.\eqno(3.2)
$$
\medskip\noindent

{\it Proof.}\,\, Let
$$
\phi(s)=\zeta(s)^{-2}L(s,\chi)^{-2}\sum_n\frac{\nu(n)^2}{n^s}
$$
which has the Euler product representation 
$$
\phi(s)=\prod_p\phi_p(s)\qquad (\sigma>1).
$$ 
For $\sigma\ge\sigma_0>0$, by checking the cases  $\chi(p)=\pm 1$ and $\chi(p)=0$ respectively, it can be seen that   
$$\phi_p(s)=1+O(p^{-2\sigma})\quad\text{if}\quad p\nmid D,$$ 
 and  
$$\phi_p(s)=(1-p^{-s})(1+O(p^{-2\sigma}))\quad\text{if}\quad p|D,$$ 
the implied constant depending on $\sigma_0$. Thus $\phi(s)$ is analytic for $\sigma>1/2$ and it satisfies
$$
\phi(s)\ll\prod_{p|D}|1-p^{-s}|\eqno(3.3)
$$
for $\sigma\ge\sigma_1>1/2$, the implied constant depending on $\sigma_1$.  The left side of (3.2) is
$$
\aligned
\ll&\sum_n\frac{\nu(n)^2}{n}\big(\exp\{-n/P^2\}-\exp\{-n/D^4\}\big)\\
&=\tp\int_{(1)}\phi(1+s)\zeta(1+s)^2L(1+s,\chi)^2(P^{2s}-D^{4s})\Gamma(s)\,ds.
\endaligned
$$
Moving the line of integration to the left appropriately and applying  standard estimates, we see that the right side is equal to the residue of the integrand plus an acceptable error $O(D^{-c})$. The residue at $s=0$ can be written as
$$
\tp\int_{|s|=\al^*}\phi(1+s)\zeta(1+s)^2L(1+s,\chi)^2(P^{2s}-D^{4s})\Gamma(s)\,ds
$$
with $\al^*=\l^{-2024}$. If $|s|=\al^*$, then
 $\phi(1+s)=O(1)$ by (3.3),
$$
L(1+s,\chi)\ll \l^{-2022}
$$
by (A) and standard estimates, and
$$
\zeta(1+s)\ll\l^{2024},\qquad (P^{2s}-D^{4s})\Gamma(s)\ll\l^9.
$$
It follows that the residue at $s=0$ is $\ll \l^{-2011}$, whence the result follows. $\Box$

\medskip\noindent

{\bf Lemma 3.2.}\,\,  {\it Assume (A) holds. Then we have}
$$
\sum_{D^4<n\le D^8}\frac{\nu(n)^2\tau_2(n)^2}{n}\ll\l^{-2007}.
$$
\medskip\noindent

{\it Proof.}\,\, As the situation is analogous to Lemma 3.1 we give a sketch only. It can be verified that
the function
$$
\phi^*(s)=\zeta(s)^{-8}L(s,\chi)^{-8}\sum_{n}\frac{\nu(n)^2\tau_2(n)^2}{n^s}
$$
is analytic for $\sigma>1/2$ and it satisfies
$$
\phi^*(s)\ll\prod_{p|D}|1-p^{-4s}|
$$
for $\sigma\ge\sigma_1>1/2$, the implies constant depending on $\sigma_1$. Also, one can verify that
$$
\int_{|s|=\al^*}\phi^*(1+s)\zeta(1+s)^8L(1+s,\chi)^8(D^{8s}-D^{4s})\Gamma(s)\,ds\ll\l^{-2007}.
$$
This completes the proof. $\Box$
\medskip\noindent

{\bf Lemma 3.3.}\,\, {\it 
For any $s$ and any complex numbers $c(n)$ we have
$$
\sum_{\psi\in\Psi}\bigg|\sum_{n\le P}\frac{c(n)\psi(n)}{n^{s}}\bigg|^2\ll  \p
\sum_{n\le P}\frac{|c(n)|^2}{n^{2\sigma}}
$$
and}
$$
\sum_{\psi\in\Psi}\bigg|\sum_{n\le P^2}\frac{c(n)\psi(n)}{n^{s}}\bigg|^2\ll P^2
\sum_{n\le P^2}\frac{|c(n)|^2}{n^{2\sigma}}.
$$
\medskip\noindent

{\it Proof.}\,\,  The first assertion follows by the orthogonality relation;
the second assertion follows by the large sieve inequality. $\Box$
\medskip\noindent

Let
$$
F\sp=\sum_{n\le D^{4}}\frac{\nu(n)\psi(n)}{n^s},\qquad G\sp=\sum_{n\le D^{4}}\frac{\upsilon(n)\psi(n)}{n^s}.
$$
By (3.1) we may write
$$
F\sp^{20}=\sum_{n\le D^{80}}\frac{\nu_{20}(n)\psi(n)}{n^s},\qquad G\sp^{20}=\sum_{n\le D^{80}}\frac{\upsilon_{20}(n)\psi(n)}{n^s}
$$
with $|\nu_{20}{n}|\le\tau_{40}(n)$, $|\upsilon_{20}{n}|\le\tau_{40}(n)$.
Write
$$
X_1(x,\psi)=\sum_{n\le x}\frac{\nu_{20}(n)\psi(n)}{n^{s_0}},\qquad 
X_2(x,\psi)=\sum_{n\le x}\frac{\upsilon_{20}(n)\psi(n)}{n^{s_0}}.
$$
By Cauchy's inequality and  the first assertion of Lemma 3.3  we obtain
$$
\sum_{\psi\in\Psi}\bigg(|X_1(D^{80},\psi)|+|X_2(D^{80},\psi)|+\int_1^{D^{80}}\frac{|X_1(x,\psi)|+|X_2(x,\psi)|}{x}\,dx\bigg)^2\ll \p\l^{1602}.
$$
Thus we conclude
\medskip\noindent

{\bf Lemma 3.4.}\,\,{\it The  inequality
$$
|X_1(D^{80},\psi)|+|X_2(D^{80},\psi)|+\int_1^{D^{80}}\frac{|X_1(x,\psi)|+|X_2(x,\psi)|}{x}\,dx<\l^{1171}\eqno(3.4)
$$
holds for all but at most $O(\p\l^{-740})$ characters $\psi$ in $\Psi$.}
\medskip\noindent

Write 
$$
X_3(x,\psi)=\sum_{D^4<n\le x}\frac{\nu(n)\psi(n)}{n^{s_0}}\quad\text{for}\quad x>D^4.
$$
Assume that (A) holds. By Cauchy's inequality,  the second assertion of Lemma 3.2 and Lemma 3.1,
$$
\sum_{\psi\in\Psi}\bigg(|X_3(P^2,\psi)|+\int_{D^4}^{P^2}\frac{|X_3(x,\psi)|}{x}\,dx\bigg)^2\ll  P^2\l^{-1993}\ll \p\l^{-1909}.
$$
Thus we conclude
\medskip\noindent

{\bf Lemma 3.5}\,\,{\it Assume that (A) holds. The  inequality
$$
|X_3(P^2,\psi)|+\int_{D^4}^{P^2}\frac{|X_3(x,\psi)|}{x}\,dx<\l^{-585}\eqno(3.5)
$$
holds for all but at most $O(\p\l^{-739})$ characters $\psi$ in $\Psi$.}
\medskip\noindent

Let
$$
\varsigma(n)=\sum_{\substack{n=lm\\ l,m\le D^4}}\nu(l)\upsilon(m).
$$
We have $\varsigma(n)=0$ unless $n=1$ or $D^4<n\le D^8$.  Write
$$
X_4(x,\psi)=\sum_{D^4<n\le x}\frac{\varsigma(n)\psi(n)}{n^{s_0}}
$$
for $x>D^4$. It is direct to verify that
$$
|\varsigma(n)|\le\sum_{n=lm}\nu(l)|\upsilon(m)|\le \nu(n)\tau_2(n).
$$
Assume that (A) holds. By  Cauchy's inequality, the first assertion of Lemma 3.2 and Lemma 3.1,
$$
\sum_{\psi\in\Psi}\bigg(|X_4(D^8,\psi)|+\int_{D^4}^{D^8}\frac{|X_4(x,\psi)|}{x}\,dx\bigg)^2\ll \p\l^{-2005}.
$$
Thus we conclude
\medskip\noindent

{\bf Lemma 3.6.}\,\,{\it Assume that (A) holds. The  inequality
$$
|X_4(D^8,\psi)|+\int_{D^4}^{D^8}\frac{|X_4(x,\psi)|}{x}\,dx<\l^{-633}\eqno(3.6)
$$
holds for all but at most $O(\p\l^{-739})$ characters $\psi$ in $\Psi$.}
\medskip\noindent

We are now in a position to give the definition of $\Psi_1$: 
Let $\Psi_1$ be the subset of $\Psi$ 
such that $\psi\in\Psi_1$ if and only if the  inequalities (3.4), (3.5) and (3.6) simultaneously hold.

Proposition 2.1 follows from Lemma 3.4, 3.5 and 3.6 immediately.

\section  *{\centerline{4. Zeros of $L\sp L(s,\pc)$ in $\Omega$}}

\medskip\noindent

In this section we prove Proposition 2.2. We henceforth assume that  $\psi(\bmod\,p)\in\Psi_1$. This assumption will not be repeated in the statements of Lemma 4.1-4.8.

We begin by proving some consequences of the inequalities (3.4)-(3.6).

\medskip\noindent

{\bf Lemma 4.1.}\,\, {\it Let   
$$
\Omega_1=\big\{s:\quad 1/2-(100\l)^{-1}\log\l<\sigma<1+(100\l)^{-1}\log\l,\quad |t-2\pi t_0|<\l_1+5\big\}.
$$
If  $s\in\Omega_1$, then}
$$
|F\sp|+|G\sp|\ll\l^{79}.
$$
\medskip\noindent

{\it Proof.}\,\, By the Stieltjes integral we may write
$$
F\sp^{20}=1+\int_1^{D^{80}}\,x^{s_0-s}d\{X_1(x,\psi)\}.
$$
For $s\in\Omega_1$ and $1\le x\le D^{80}$ we have
$$
|x^{s_0-s}|\ll\l,\qquad\bigg|\frac{d}{dx}\big(x^{s_0-s}\big)\bigg|\ll x^{-1}\l^{406}.
$$
Hence, by partial integration,
$$
|F\sp|^{20}\ll 1+\l^{406}\bigg(|X_1(D^{80},\psi)|+\int_1^{D^{80}}\frac{|X_1(x,\psi)|}{x}\,dx\bigg).
$$
For $G\sp$ an entirely analogous bound is valid. The result now follows by (3.4). $\Box$

\medskip\noindent

{\bf Lemma 4.2.}\,\, {\it If  $s\in\Omega_1$, then}
$$
F\sp G\sp=1+O(\l^{-227}).
$$
\medskip\noindent

{\it Proof}.\,\, We have

$$
F\sp G\sp-1=\sum_{D^4<n\le D^8}\frac{\varsigma(n)\psi(n)}{n^s}
=\int_{D^4}^{D^8}x^{s_0-s}\,d\{X_4(x,\psi)\}.
$$
Thus, similar to (4.2), by partial integration we obtain
$$
F\sp G\sp-1\ll\l^{406}\bigg(|X_4(D^8,\psi)|+\int_{D^4}^{D^8}\frac{|X_4(x,\psi)|}{x}\,dx\bigg),
$$
 the right side being $O(\l^{-227})$ by (3.6). $\Box$
\medskip\noindent

\medskip\noindent

{\bf Lemma 4.3.}\,\, {\it Let
$$
\Omega_2=\big\{s:\quad 1/2-\l^{-1}<\sigma<1+\l^{-1},\quad |t-2\pi t_0|<\l_1+4\big\}.
$$ 
If  $s\in\Omega_2$, then}
$$
\frac{F'}{F}\sp=O(\l).
$$

{\it Proof}.\,\,  Assume  $|w|\le(200\l)^{-1}\log\l$,  so that $s+w\in\Omega_1$. By Lemma  4.1 and 4.2, 
$$
\l^{-88}\ll|F(s+w,\psi)|\ll\l^{88}.
$$
Thus the logarithm
$$
\mathfrak{l}(s,w):=\log\frac{F(s+w,\psi)}{F\sp},
$$
 which vanishes at $w=0$, is analytic in $w$, and it satisfies 
 $$
 \Re\{\mathfrak{l}(s,w)\}\ll\log\l.
 $$
Since
$$
\frac{F'}{F}\sp=\frac{\partial}{\partial w}\mathfrak{l}(s,w)|_{w=0},
$$
the result follows by Lemma 4 of  [15,  Chapter 2]. $\Box$
\medskip\noindent

We proceed to establish an approximate formula for $L\sp L(s,\pc)$. For this purpose we first introduce a weight $g(x)$ that will find application at various places. Let
$$
g(x)=\tp\int_{(c)}\frac{x^w\omega_1(w)}{w}\,dw\qquad (c>0)
$$
with
$$
\omega_1(w)=\exp\{w^2/(4\l^{30})\}.
$$
We may write
$$
g(x)=\tp\int_{(c)}\bigg(\int_{0}^{ x}y^{w-1}\,dy\bigg)\omega_1(w)\,dw.
$$
Since
$$
\tp\int_{(c)}\exp\{(\log y)w+w^2/(4\l^{30})\}\,dw=\frac{\l^{15}}{\sqrt{\pi}}\exp\{-\l^{30}(\log y)^2\},
$$
it follows, by changing the order of integration, that
$$
g(x)=\frac{\l^{15}}{\sqrt{\pi}}\int_{0}^{ x}\exp\{-\l^{30}(\log y)^2\}\,\frac{dy}{y}.
$$
This yields. by substituting $t=\l^{15}\log y$,
$$
g(x)=\frac{1}{\sqrt{\pi}}\int_{-\infty}^{\l^{15}\log x}\exp\{-t^2\}\,dt.\eqno(4.1)
$$
Thus the function $g(x)$ is increasing and it satisfies $0<g(x)<1$. Further we have
$$
g(x)=1+O(\exp\{-\l^{30}\log^2x\})\quad\text{if}\quad x\ge 1\eqno(4.2)
$$
and
$$
g(x)=O(\exp\{-\l^{30}\log^2x\})\quad\text{if}\quad x\le 1.\eqno(4.3)
$$

 Note that 
$\chi\psi$ is a primitive character $(\bmod Dp)$.  Write
$$
\tz\sp=Z\sp Z(s,\pc),
$$
so that
$$
L\sp L(s,\pc)=\tz\sp L(1-s,\bp)L(1-s,\chi\bp).\eqno(4.4)
$$
Recall that $s_0=1/2+2\pi it_0$. Assume that
$$
|\Re(s-s_0)|<100,\,\qquad|\Im(s-s_0)|<\l_1+3.
$$
 By (2.4) with $\theta=\psi$ and $\theta=\pc$ we have
$$
\tz\sp=\chi(-1)\tau(\psi)\tau(\pc)(Dp^2)^{-s}\vartheta(s)^2(1+O(e^{-\pi t})).
$$
This yields, by Stirling's formula,
$$
|\tz\sp|= (Dp^2t_0^2)^{1/2-\sigma}(1+o(1))\eqno(4.5)
$$
and
$$
\frac{\tz'}{\tz}\sp=-2\log P+O(\l).\eqno(4.6)
$$

{\bf Lemma 4.4.}\,\, {\it Let
$$
\Omega_3=\big\{s:\quad 1/2-\al<\sigma<1+\al,\quad |t-2\pi t_0|<\l_1+3\big\}.
$$ 
If  $s\in\Omega_3$, then}
$$
L\sp L(s,\pc)=F\sp+\tz\sp F(1-s,\bp)+O(\l^{-179}).
$$
\medskip\noindent

{\it Proof}.\,\, By the residue theorem,
$$
L\sp L(s,\pc)=\tp\bigg(\int_{(1)}-\int_{(-\sigma-1/2)}\bigg)L(s+w,\psi)L(s+w,\pc)P^{(9/5)w}\,\frac{\omega_1(w)\,dw}{w}.
$$
By (4.2) and (4.3),
$$\aligned
\tp\int_{(1)}&L(s+w,\psi)L(s+w,\pc)P^{(9/5)w}\,\frac{\omega_1(w)\,dw}{w}\\
&=F\sp+\sum_{D^4<n<P^2}\frac{\nu(n)\psi(n)}{n^s}g\bigg(\frac{P^{9/5}}{n}\bigg)+O(\ep)
\endaligned
$$
where
$$
\ep=\exp\{-c\l^{10}\}.
$$
By (3.5) and partial summation, the second sum on the right side above is 
$$
=\int_{D^4}^{P^2}g\bigg(\frac{P^{9/5}}{x}\bigg)x^{s_0-s}\,dX_3(x,\psi)\ll\l_1
\bigg(|X_3(P^2,\psi)|+\int_{D^4}^{P^2}\frac{|X_3(x,\psi)|}{x}\,dx\bigg)\ll\l^{-180}.
$$
On the other hand, by the functional equation, for $u=-\sigma-1/2$,
$$
L(s+w,\psi)L(s+w,\pc)=\tz(s+w,\psi)\sum_n\frac{\nu(n)\bp(n)}{n^{1-s-w}}.
$$
This sum is split into three sums according to $n\le D^4$, $D^4<n\le P^2$ and $P^2<n$. The proof is therefore reduced to showing that
$$\aligned
\tp\int_{(-\sigma-1/2)}&\tz(s+w,\psi)F(1-s-w,\bp)P^{(9/5)w}\,\frac{\omega_1(w)\,dw}{w}\\
&=-\tz(s,\psi)F(1-s,\bp)+O(P^{-c}),\endaligned\eqno(4.7)
$$
$$
\tp\int_{(-\sigma-1/2)}\tz(s+w,\psi)\bigg(\sum_{D^4<n\le P^2}\frac{\nu(n)\bp(n)}{n^{1-s-w}}\bigg)P^{(9/5)w}\,\frac{\omega_1(w)\,dw}{w}\ll\l^{-179}\eqno (4.8)
$$
and
$$
\tp\int_{(-\sigma-1/2)}\tz(s+w,\psi)\bigg(\sum_{n> P^2}\frac{\nu(n)\bp(n)}{n^{1-s-w}}\bigg)P^{(9/5)w}\,\frac{\omega_1(w)\,dw}{w}\ll P^{-c}.\eqno(4.9)
$$

 To prove (4.7) we move the contour of integration to the vertical segments
$$\begin{cases}
w=10+iv&\quad\text{with}\quad |v|<\l^{20},\\
 w=-\sigma-1/2+iv&\quad\text{with}\quad |v|\ge\l^{20},
\end{cases}$$
and to the two connecting horizontal segments
$$
w=u\pm i\l^{20}\quad\text{with}\quad -\sigma-1/2\le u\le 10.
$$
By  a trivial  bound for $\omega_1(w)$, (4.5) and the residue theorem we obtain (4.7).

To prove (4.8) we move the contour of integration to the vertical segments
$$\begin{cases}
w=-\al+iv&\quad\text{with}\quad |v|<\l^{20},\\
 w=-\sigma-1/2+iv&\quad\text{with}\quad |v|\ge\l^{20},
\end{cases}$$
and to the two connecting horizontal segments
$$
w=u\pm i\l^{20}\quad\text{with}\quad -\sigma-1/2\le u\le-\al.
$$
By  a trivial  bound for $\omega_1(w)$ and (4.5) we see that the left side of (4.8) is
$$
\ll P^{1-2\sigma}\int_{-\l^{20}}^{\l^{20}}\bigg|\sum_{D^4<n\le P^2}\frac{\nu(n)\psi(n)}{n^{s^*+iv}}\bigg|\,\frac{dv}{\al+iv}+\ep
$$
with $s^*=1+\al-\bar{s}$. By partial integration,
$$
\sum_{D^4<n\le P^2}\frac{\nu(n)\psi(n)}{n^{s^*+iv}}=\int_{D^4}^{P^2}x^{s_0-s^*-iv}\,dX_3(x,\psi)\ll P^{2\sigma-1}\l_1
\bigg(|X_3(P^2,\psi)|+\int_{D^4}^{P^2}\frac{|X_3(x,\psi)|}{x}\,dx\bigg)
$$
for $|v|\le\l^{20}$. From this and (3.5) we obtain (4.8).

The estimate (4.9) follows by  moving the contour of integration to the vertical segments
$$\begin{cases}
w=-\sigma-10+iv&\quad\text{with}\quad |v|<\l^{20},\\
 w=-\sigma-1/2+iv&\quad\text{with}\quad |v|\ge\l^{20},
\end{cases}$$
and to the two connecting horizontal segments
$$
w=u\pm i\l^{20}\quad\text{with}\quad -\sigma-10\le u\le-\sigma-1/2,
$$
and applying (4.5) and trivial bounds for $\omega_1(w)$ and the involved sum. $\Box$

  In order to prove Proposition 2.2,  it is appropriate to deal with the function 
$$
\a\sp=\frac{L\sp L(s,\pc)}{F\sp}.
$$ 
By Lemma 4.2, $\a\sp$ is analytic and it has the same zeros as 
$L(s,\psi)L(s,\psi\chi)$ in $\Omega_1$. Further, for $s\in\Omega_1$,  we have 
$$
F\sp^{-1}\ll\l^{79}
$$
by Lemma 4.1 and 4.2. This together with Lemma 4.4 implies that
$$
\a\sp=1+\b\sp+O(\l^{-100}),\eqno(4.10)
$$
for $s\in\Omega_3$, where
$$
\b\sp=\tz\sp\frac{F(1-s,\bp)}{F\sp}.
$$

The proof of Proposition 2.2 is reduced to proving three lemmas as follows.
\medskip\noindent

{\bf Lemma 4.5.} \,\,{\it If
$$
\frac12+\al^2<\sigma<1,\qquad |t-2\pi t_0|<\l_1+2,
$$
then}
$$A\sp\ne 0.$$ 
\medskip\noindent

{\it Proof.} \,\,We discuss in two cases.
\medskip\noindent

{\it Case 1.} $1/2+\l^{-1}\le\sigma<1$. 
\medskip\noindent

By Lemma 4.2 and trivial estimation,
$$\frac{F(1-s,\bp)}{F\sp}\ll D^c.$$ 
Hence, by (4.5),
$$
\b\sp\ll P^{1/2-\sigma}.
$$
 The result now follows by (4.10).
\medskip\noindent

{\it Case 2.} $1/2+\al^2\le\sigma<1/2+\l^{-1}$.
\medskip\noindent

Assume $1/2\le\sigma'\le\sigma$. Then both $\sigma'+it$ and $1-\overline{(\sigma'+it)}$ lie in $\Omega_2$. Hence,  by Lemma 4.3 and (4.6),
$$\aligned
\frac{\b'}{\b}(\sigma'+it,\psi)=&\frac{\tz'}{\tz}(\sigma'+it,\psi)-\frac{F'}{F}(\sigma'+it,\psi)-\frac{F'}{F}(1-\sigma'-it,\bp)\\
=&-2\log P+O(\l).
\endaligned\eqno(4.11)$$
Since
$|\b(1/2+it,\psi)|=1$,
it follows  that
 $$
 \log|\b\sp|=\Re\bigg\{\int_{1/2}^{\sigma}\frac{\b'}{\b}(\sigma'+it)\,d\sigma\bigg\}<(1/2-\sigma)\log P,
 $$
 Hence, by (4.10), 
$$ |\a\sp|>1-P^{1/2-\sigma}+O(\l^{-100})\gg\al. \qquad \Box
$$
\medskip\noindent

   {\bf Lemma 4.6.} {\it Suppose  $\rho=\beta+i\gamma$  is a zero of $\a\sp$ satisfying
$$
\frac12\le\beta<\frac12+\al^2,\qquad |\gamma-2\pi t_0|<\l_1+2.
$$   
Then $\beta=1/2$, $\a'\ss\ne 0$ and $\a(1/2+i\gamma+w,\psi)\ne 0$ if $0<|w|<\al(1-c'\al\l)$ where $c'$ is a sufficiently constant. }
\medskip\noindent

{\it Proof.}\,\, It suffices to show that the function $\a(1/2+i\gamma+w,\psi)$ has exactly one zero inside the circle
$|w|=\al(1-c'\al\l)$, counted with multiplicity. By the Rouch\'{e} theorem, this can be reduced  to proving that
$$
|\a(1/2+i\gamma+w,\psi)- (1-P^{-2w})|<|1-P^{-2w}|\quad\text{for}\quad |w|=\al(1-c'\al\l),
$$
since the function $1-P^{-2w}$
  has exactly one zero inside this circle  which is at $w=0$. In fact, we can prove that
  $$
  \a(1/2+i\gamma+w,\psi)- (1-P^{-2w})\ll\al\l\eqno(4.12)
$$
if $|w|<2\al$, the implied constant being independent of $c'$, and
$$
|1-P^{-2w}|>6c'\al\l\eqno(4.13)
$$
if $ |w|=\al(1-c'\al\l)$.

Assume $|w|<2\al$. By (4.10) we have
$$
 \a(1/2+i\gamma+w,\psi)- (1-P^{-2w})= \b(1/2+i\gamma+w,\psi)+P^{-2w}+O(\l^{-100}).
 $$

Noting that both $s$ and $1-s$ are in $\Omega_2$ if $s$ lies on the segment connecting $\rho$ and $1/2+i\gamma+w$, by (4.11) we have
$$\aligned
\frac{\b(1/2+i\gamma+w,\psi)}{\b\ss}=&\exp\bigg\{\int_{\rho}^{1/2+i\gamma+w}\frac{\b'}{\b}(s,\psi)\,ds\bigg\}
=\exp\big\{-2w\log P+O(\al\l)\big\}\\=&P^{-2w}+O(\al\l),
\endaligned
$$
Since $\b\ss=-1+O(\l^{-100})$ by (4.10), the estimate (4.12) follows.

Now assume $ |w|=\al(1-c'\al\l)$. Write $$w=\al(1-c'\al\l)(\cos\theta+i\sin\theta).$$ Then
$$
|P^{-2w}|=\exp\big\{-2\pi(1-c'\al\l)\cos\theta\big\},\quad \Im\{P^{-2w}\}=-|P^{-2w}|\sin\big\{2\pi(1-c'\al\l)\sin\theta\big\};
$$
$$
\text{if}\quad\cos\theta>c'\al\l,\quad\text{then}\quad |P^{-2w}|<1-6c'\al\l;
$$
$$
\text{if}\quad\cos\theta<-c'\al\l,\quad\text{then}\quad |P^{-2w}|>1+6c'\al\l;
$$
$$
\text{if}\quad |\cos\theta|\le c'\al\l, \quad\text{then}\quad |\Im\{P^{-2w}\}|>6c'\al\l,
$$
since
$$
\sqrt{1-(c'\al\l)^2}\le|\sin\theta|\le 1,
$$ so that
$$|\sin\big\{2\pi(1-c'\al\l)\sin\theta\big\}|=2\pi c'\al\l(1+o(1)).
$$
 In either case (4.16) holds. $\Box$
\medskip\noindent

Lemma 4.5 and 4.6 together imply the assertions (i) and  (ii) of Proposition 2.2. It is also proved that the gap between any distinct zeros of $\a\sp$ in $\Omega$ is $>\al(1-c'\al\l)$. To complete the proof of the gap assertion (iii), it now suffices to prove
\medskip\noindent

{\bf Lemma 4.7.} {\it Suppose  $\rho=1/2+i\gamma$ is a zero of $\a\sp$ satisfying $|\gamma-2\pi t_0|<\l_1+2$. Then the function $\a(\rho+w,\psi)$ has exactly three zeros inside the circle $|w|=\al(1+c'\al\l)$, counted with multiplicity.}
\medskip\noindent

{\it Proof.}\,\, In a way similar to the proof of Lemma 4.6, it is direct to verify that 
$$
|\a(\rho+w,\psi)- (1-P^{-2w})|<|1-P^{-2w}|
$$
if $ |w|=\al(1+c'\al\l)$. Hence, the functions $\a(\rho+w,\psi)$ and $1-P^{-2w}$
  have the same number of zeros inside this circle, while the later has exactly three zeros inside the same circle which are at
$w=0$, $w=i\al$ and $w=-i\al$. $\Box$
\medskip\noindent

We conclude this section by giving a result which is implied in the proof of Proposition 2.2.
\medskip\noindent

{\bf Lemma 4.8.}\,\,{\it Assume that $\rho$ is a zero of $L\sp L(s,\pc)$ in $\Omega$. Then we have}
$$
\tz\ss^{-1}=-G\ss F(1-\rho,\bp)+O(\l^{-100}).
$$

{\it Proof.}\,\,It  follows  from Lemma 4.4 that
$$
F\ss+\tz\ss F(1-\rho,\bp)\ll\l^{-179}.
$$
The result follows by multiplying both sides  by $\tz\ss^{-1}G\ss $ and applying Lemma 4.1.  $\Box$
\medskip\noindent

\section  *{\centerline{5. Some  analytic lemmas}}
\medskip\noindent

{\bf Lemma 5.1.}\,\,{\it Suppose $\psi(\bmod\,p)\in\Psi$, $|\sigma-1/2|
\le\al$, $|t-2\pi t_0|<\l_1+2$, $u=0$   and $|v|<\L^{20}$.  Then
$$
\frac{Z(s+w,\psi)-Z(s,\psi)(pt_0)^{-w}}{w}\ll\l^{-114}\eqno(5.1)
$$
$$
\frac{Z(s+w,\psi)-Z(s,\psi)(Pt_0)^{-w}}{w}\ll\l^{-68}\eqno(5.2)
$$
$$
\frac{Z(s+w,\pc)-Z(s,\pc)(Dpt_0)^{-w}}{w}\ll\l^{-114}\eqno(5.3)
$$
and}
$$
\frac{Z(s+w,\pc)-Z(s,\pc)(DPt_0)^{-w}}{w}\ll\l^{-68}\eqno(5.4)
$$

{\it Proof.}\,\, We have
$$
\frac{Z(s+w,\psi)}{Z\sp}=\exp\bigg\{\int_{0}^{w}\frac{Z'}{Z}(s+w',\psi)\,dw'\bigg\}.
$$
 Assume $|w'|\le |w|$. By (2.6) and the Stirling formula,
$$
\frac{Z'}{Z}(s+w',\psi)=-\log p+\frac{\vartheta'}{\vartheta}(s+w')+O(\ep)=-\log (pt/2\pi)+O(1/t_0).
$$
Hence
$$
Z(s+w,\psi)=Z\sp (pt/2\pi)^{-w}+O(|w|/t_0).
$$
This yields (5.1) and (5.2) since
$$
pt/2\pi=pt_0(1+O(\l^{-114})),\quad pt/2\pi=Pt_0(1+O(\l^{-68})).
$$
The proofs of (5.3) and (5.4) are similar. $\Box$
\medskip\noindent

{\bf Lemma 5.2.}\,\,{\it Let $\psi$ and $s$ be as in Lemma 5.1. Then}
$$
\frac{Y(s+\be_1,\psi)Y(s+\be_2,\psi)Y(s+\be_3,\psi)}{Y\sp}
=(pt_0)^{\be_3}Z(s,\psi)^{-1}(1+O(\l^{-123})).
$$
\medskip\noindent

{\it Proof.}\,\,The left side  is
$$
Z(s,\psi)^{-1}
\prod_{1\le j\le 3}\bigg(\frac{Y(s+\be_j,\psi)}{Y\sp}\bigg).
$$
By (2.6) and the Stirling formula, for $|w|<5\al$,
$$
\frac{Y'}{Y}(s+w,\psi)=-\frac12\frac{Z'}{Z}(s+w,\psi)=\frac12\log (pt_0)+O(\l^{-114}).
$$
 Hence, for $1\le j\le 3$,
$$
\frac{Y(s+\be_j,\psi)}{Y\sp}=\exp\bigg(\int_0^{\be_j}\frac{Y'}{Y}(s+w,\psi)\,dw\bigg)=(pt_0)^{\be_j/2}(1+O(\l^{-123})).
$$
The result now follows since
$$
\frac{\be_1+\be_2+\be_3}{2}=\be_3.\qquad \Box
$$
\medskip\noindent

Recall that $\vartheta(s)$ and  $\omega(s)$ are given by (2.3) and (2.15) respectively. It is known that
$$
\vartheta(1-s)=2(2\pi)^{-s}\Gamma(s)\cos(\pi s/2).
$$ 
For $t>1$ we have
$$
\vartheta(1-s)=\vartheta^*(1-s)\big(1+O(e^{-\pi t})\big)\eqno(5.5)
$$
where
$$ 
\vartheta^*(1-s)=(2\pi)^{-s}\Gamma(s)e(-s/4).
$$
Let
$$
\Delta_1(x)=\tp\int_{(3/2)} x^{-s}\vartheta^*(1-s)\omega(s)\,ds \eqno(5.6)
$$
and
$$
\Delta(x)=\Delta_1(x)e(x).\eqno(5.7)
$$
Note that
$$
\omega(1/2+2\pi ix)=\frac{\sqrt{\pi}}{\l_2}\exp\bigg\{-\bigg(\frac{\pi(x-t_0)}{\l_2}\bigg)^2\bigg\}>0.
$$

\medskip\noindent
{\bf Lemma 5.3.}\,\,
 {\it  If $x\le t_0^{1.02}$, then
 $$
\Delta(x)=\omega(1/2+2\pi ix)\big(1+O(\al)\big)+O(\ep);\eqno(5.8)
$$
 if $x>t_0^{1.02}$, then}
$$
\Delta(x)\ll\exp\{-(10^{-2}\l_2\log x)^2\}+\exp\{-x^{0.99}/\l_2\}.\eqno(5.9)
$$
\medskip\noindent

{\it Proof.}\,\,By the Mellin transform (see [1], Lemma 2) we have
$$
\Delta_1(x)=\int_0^\infty \exp\{(s_0-1)\log y-\l_2^2\log^2 y-2\pi ixy\}\,dy
$$
 where the logarithm vanishes at $y=1$. This yields, by substituting $y=e^u$,
 $$
\Delta(x)=\int_{-\infty}^\infty \exp\{s_0u-\l_2^2u^2-2\pi ix(e^u-1)\}\,du.\eqno(5.10)
$$

First assume $x\le t_0^{1.02}$. We may write 
$$
\exp\{s_0w-\l_2^2w^2-2\pi ix(e^w-1)\}=\exp\{2\pi i(t_0-x)w-\l_2^2w^2\}f^*(x,w)
$$
with
$$
f^*(x,w)=\exp\{(1/2)w-2\pi ix(e^w-1-w)\}.
$$
By the relation
$$
\int_{-\infty}^\infty \exp\{2\pi i(t_0-x)u-\l_2^2u^2\}\,du=\omega(1/2+2\pi ix)
$$
and Cauchy's theorem, the proof of (5.8) is reduced to showing that
$$
\int_{L_j} \exp\{2\pi i(t_0-x)w-\l+2^2w^2\}\big(f^*(x,w)-1\big)\,dw\ll \al \omega(1/2+2\pi ix)+\ep\eqno(5.11)
$$
for $1\le j\le 5$, where $L_j$ denote the segments
$$
L_1=(-\infty,\,-u^*],\quad L_2=[-u^*,\,-u^*+iv^*],\quad L_3=[-u^*+iv^*,\,u^*+iv^*],
$$
$$
L_4=\{u^*+iv^*,\,u^*],\quad L_5=[u^*,\,\infty)
$$
with
$$
u^*=\l_2^{-1}\l^{5},\qquad v^*=\frac{\pi(t_0-x)}{\l_2^2}.
$$

If $w\in L_1\cup L_5$,  then
$$
\exp\{2\pi i(t_0-x)w-\l_2^2w^2\}\big(f^*(x,w)-1\big)\ll\exp\{-\l_2^2u^2\}(1+e^{u/2}).
$$
Thus the left side of (5.11) is trivially $O(\ep)$ if $j=1,5$. By simple estimates,
$$
\int_{L_j} \big|\exp\{2\pi i(t_0-x)w-\l_2^2w^2\}\,dw\big|\ll \omega(1/2+2\pi ix)+\ep
$$
if $j=2,3,4$.
Since $t_0^{3.06}/\l_2^4\ll\al$, we have
$$
f^*(x,w)-1\ll |w|+x|w|^2\ll\al
$$
 for $w\in L_2\cup L_3\cup L_4$. These estimates together imply (5.11).
 
Now assume $x>t_0^{1.02}$. By Cauchy's theorem, the proof of (5.9) is reduced to showing that
$$
\int_{L_j'} \exp\{s_0w-\l_2^2w^2-2\pi ix(e^w-1)\}\,dw\ll\exp\{-(10^{-2}\l_2\log x)^2\}+\exp\{-x^{0.99}/\l_2\}\eqno(5.12)
$$
where $L'_j$, $1\le j\le 3$, denote the segments
$$
L'_1=\{-\infty,\,-10^{-2}\log x],\quad L'_2=[-10^{-2}\log x,\,-10^{-2}\log x-i/\l_2],$$
$$
 L'_3=[-10^{-2}\log x-i/\l_2,\,\infty-i/\l_2).
$$
We have
$$
\Re\{s_0w-\l_2^2w^2-2\pi ix(e^w-1)\}=\frac12u-\l_2^2(u^2-v^2)+2\pi v\bigg(\frac{xe^u\sin v}{v}-t_0\bigg).\eqno(5.13)
$$
If $w\in L'_1$, then the right side of (5.13)  is
$$
\le\frac12u-(10^{-2}\l_2\log x)^2.
$$
This yields (5.12) with $j=1$. For  $w\in L'_2\cup L_3'$ we have 
$$xe^u\ge x^{0.99}>2t_0\quad\text{and}\quad  -1/\l_2\le v\le 0, 
$$
so that
$$
2\pi \bigg(\frac{xe^u\sin v}{v}-t_0\bigg)>x^{0.99}.
$$
If $w\in L_2'$, then $(\l_2u)^2=(10^{-2}\l_2\log x)^2$, so the right side of (5.13) is
$$
\le -(10^{-2}\l_2\log x)^2+O(1).
$$
This yields (5.12) with $j=2$ . If $w\in L'_3$, then $v=-1/\l_2$, and the right side of (5.13) is
$$
<1+\frac12u-(\l_2u)^2-x^{0.99}/\l_2.
$$
This yields (5.12) with $j=3$ . $\Box$
\medskip\noindent

\medskip\noindent

As a  consequence of Lemma 5.3, the Mellin transform
$$
\delta(s):=\int_{0}^{\infty}\Delta(x)x^{s-1}\,dx\eqno(5.14)
$$
 is analytic for $\sigma>0$.
\medskip\noindent

\medskip\noindent

{\bf Lemma 5.4.}
(i). {\it  If $1/2\le\sigma\le 2$, then}
$$
\delta(s)\ll\l^c|s|^{-2}.
$$
\medskip\noindent

(ii).\,\,{\it If $|s-1|<10\al$, then}
$$
\delta(s)=1+O(\al\log\l). 
$$
\medskip\noindent

{\it Proof.}\,\, (i). Using partial integration twice we obtain
$$
\delta(s)=\frac{1}{s(s+1)}\int_0^\infty\,\Delta''(x)x^{s+1}\,ds.
$$
By (5.10) we have
 $$
\Delta^{''}(x)=-4\pi^2\int_{-\infty}^\infty( e^u-1)^2 \exp\{s_0u-B^2u^2-2\pi ix(e^u-1)\}\,du.
$$
Thus some upper bounds for $\Delta''(x)$ analogous to Lemma 5.3 can be obtained, and (i)
 follows. 

(ii). Assume $|s-1|<10\al$. By Lemma 5.3, on the right side of (5.14), the integral on the part $|x-t_0|\ge \l_1$ contributes $O(\ep)$. For $|x-t_0|<\l_1$ we have
$$
x^{s-1}=1+O(\al\log \l).
$$
Since
$$
\int_{0}^{\infty}\omega(1/2+2\pi ix)\,dx=1+O(\ep),
$$
 (ii) follows. $\Box$
\medskip\noindent

{\it Throughout the rest of this paper we assume that (A) holds. This assumption will not be repeated in the statements of the lemmas and propositions in the sequel. }

\medskip\noindent

  The next  two lemmas are weaker forms of the Deuring-Heillbronn Phenomenon.
 \medskip\noindent
 
  {\bf Lemma 5.5.}\,\, {\it  The function $L(s,\chi)$ has a simple real zero $\tilde{\rho}$ such that 
  $$
  1-\tilde{\rho}=O(\l^{-2022})\eqno(5.15)
  $$
  and $L(s,\chi)$ has no other zeros in the region}
  $$
  \sigma>1-2\l^{-1},\quad |t|<2D.
  $$
  \medskip\noindent

  {\bf Lemma 5.6.}\,\, {\it  For any primitive character $\theta(\bmod\,r)$ with $r<T$ and $\theta\ne\chi$,
  $$
  \sum_{p\sim P}\theta(p)p^{1+it}\ll\p\exp\{-\l^{9/2}\}
  $$  
   if $|t|\le D$.}
   \medskip\noindent

\medskip\noindent

 {\bf Lemma 5.7.}\,\, {\it  We have}
 $$
 L'(1,\chi)\gg\frac{D}{\varphi(D)}.
 $$

{\it Proof.}\,\, The right side of the equality
$$
\tp\int_{(1)}\zeta(1+s)L(1+s,\chi)\frac{D^{4s}\omega_1(s)}{s}\,ds=\sum_n\frac{\nu(n)}{n}g\bigg(\frac{D^4}{n}\bigg).
$$
 is
$$
\gg\sum_{n|D}\frac{1}{n}\gg\frac{D}{\varphi(D)},
$$
since $\nu(n)=1$ if $n|D$, while the left side is, by moving the line of integration to the left appropriately and applying (A), equal to $L'(1,\chi)+o(1)$. $\Box$

\medskip\noindent

{\bf Lemma 5.8.}\,\, {\it If
$$
\al \le|s-1|\le 10\al,
$$
then
$$
L(s,\chi)=L'(1,\chi)(s-1)+O(\al_2)
$$
where}
$$
\al_2=\l^{-15}.
$$

{\it Proof.}\,\, This follows from the relation
$$
L(s,\chi)=L(1,\chi)+L'(1,\chi)(s-1)+\int_1^s(s-w)L^{''}(w,\chi)\,dw,
$$
(A) and a simple bound for $L^{''}(w,\chi)$. $\Box$

\medskip\noindent

{\bf Lemma 5.9.}\,\,{\it Suppose $\psi\in\Psi_1$, $|\sigma-1/2|\le\al$, $|t-2\pi t_0|\le\l_1+10$ and $|s-\rho|\gg\al$ for any zero $\rho$ of $L\sp$. Then}
$$
\frac{L(s+\be_1,\psi)}{L\sp}\ll\log P. 
$$

{\it Proof.}\,\,It is known that
$$
\frac{L'}{L}(s',\psi)=\sum_{|\rho-s'|<1}\frac{1}{s'-\rho}+O(1/\al)\eqno(5.16)
$$
for $|\Re\{s'\}-1/2|\le\al$ and $|\Im\{s'-2\pi t_0\}|<\l_1+10$, where $\rho$ runs through the zeros of $L(s',\psi)$.

We can assume $L(s+\be_1,\psi)\ne 0$. Suppose $\sigma\ge 1/2$. By (5.16),
$$
-\Re\bigg\{\frac{L'}{L}(u+it+\be_1,\psi)\bigg\}<O(1/\al)\quad\text{for}\quad \sigma<u<1/2+\al,
$$
so that 
$$
\log\frac{|L(s+\be_1,\psi)|}{|L(1/2+\al+it+\be_1,\psi)|}<O(1).
$$
By (9.1) and the condition $|s-\rho|\gg\al$ for any $\rho$,
$$
\Re\bigg\{\frac{L'}{L}(u+it,\psi)\bigg\}=O(1/\al)\quad\text{for}\quad \sigma<u<1/2+\al,
$$
so that 
$$
\log\frac{|L(1/2+\al+it,\psi)|}{|L(s,\psi)|}<O(1).
$$ 
Further, by (5.16) and Proposition 2.2 (iii),
$$
\bigg|\frac{L'}{L}(1/2+\al+it',\psi)\bigg|<\frac{1}{\al}\sum_{k<1/al}\frac{1}{k}+O(1/\al)<\frac{\log\log P}{\al}+O(1/\al)
$$
for $t<t'<t+|\be_1|$, so that
$$
\log\frac{|L(1/2+\al+it+\be_1,\psi)|}{|L(1/2+\al+it,\psi)|}<\log\log P+O(1).
$$
Combining theses estimates we obtain the result. In the case $\sigma<1/2$ the proof is analogous. $\Box$

 \medskip\noindent

\section *{\centerline{6. Approximate formula for $L(s,\psi)$}}
\medskip\noindent

Write
$$
g^*(y)=
\begin{cases}
g(y)&\quad\text{if}\quad y>1/2,\\
0&\quad\text{otherwise}.
\end{cases}
$$
Let
$$
P_4=PT^{-2}t_0\quad\text{with}\quad T=\exp\{\l^{1.1}\}.
$$

{\bf Lemma 6.1.}\,\,{\it Suppose   $\psi(\bmod\,p)\in\Psi$, $|\sigma-1/2|<2\al$ and $|t-2\pi t_0|<\l_1+2$. Then 
$$
L(s,\psi)=K(s,\psi)+Z\sp N(1-s,\bp)+O(E_1\sp),
$$
where
$$
K(s,\psi)=\sum_{n}\frac{\psi(n)}{n^s}\,g^*\bigg(\frac{P_4}{n}\bigg)$$
$$
N(1-s,\bp)=\sum_{n}\frac{\bp(n)}{n^{1-s}}\,g^*\bigg(\frac{T^2}{n}\bigg)$$
and where }
$$\aligned
 E_1\sp=\l^{-68}\int_{-\l^{20}}^{\l^{20}}\bigg|\sum_{n<T^3}\frac{\psi(n)}{n^{s+iv}}\bigg|\omega_1(iv)dv+\ep.
\endaligned
$$
{\it Proof.}\,\, By (4.3) we have
$$
\tp\int_{(1)}L(s+w,\psi)\frac{P_4^w\omega_1(w)\,dw}{w}=K(s,\psi)+O(\ep).
$$
The left side above is, by moving the line of integration to $u=-1$, equal to
$$
L(s,\psi)+\tp\int_{(-1)}L(s+w,\psi)\frac{P_4^w\omega_1(w)\,dw}{w}.
$$
It therefore suffices to show that
$$
\tp\int_{(-1)}L(s+w,\psi)\frac{P_4^w\omega_1(w)\,dw}{w}=-Z\sp N(1-s,\bp)+O(E_1\sp).\eqno(6.1)
$$

For $u=-1$ we have, by the functional equation (2.2) with $\theta=\psi$, 
$$
L(s+w,\psi)=Z(s+w,\psi)\bigg(\sum_{n<T^3}\frac{\bp(n)}{n^{1-s-w}}+\sum_{n
\ge T^3}\frac{\bp(n)}{n^{1-s-w}}\bigg)
$$
We first show that 
$$
\int_{(-1)}Z(s+w,\psi)\bigg(\sum_{n
\ge T^3}\frac{\bp(n)}{n^{1-s-w}}\bigg)\frac{P_4^w\omega_1(w)\,dw}{w}\ll\ep.\eqno(6.2)
$$

 We move the contour of integration in (6.2) to the vertical segments
$$
u=-1,\quad |v|>\l^{20}
$$ and 
$$
u=-\l^9,\quad |v|\le\l^{20}
$$
with the horizontal connecting segments
$$
-\l^9\le u\le -1\quad |v|=\l^{20}.
$$
By a trivial bound for $\omega_1(w)$ and simple estimates, , the  integrals on the segments $u=-1$, $|v|>\l^{20}$  contribute  $\ll\ep$. On the other hand, in the case  $-\l^9\le u\le -1$, $|v|\le\l^{20}$, by (4.) we have
$$
Z(s+w,\psi)P_4^w\ll(2T^2)^{-u}
$$
and
$$
\sum_{n\ge T^3}\frac{\bp(n)}{n^{1-s-w}}\ll T^{3u+1/2},
$$
whence (6.2) follows.  The proof of (6.1) is therefore reduced to showing that
$$\aligned
\tp\int_{-1-i\l^{20}}^{-1+i\l^{20}}Z(s+w,\psi)&\bigg(\sum_{n
<T^3}\frac{\bp(n)}{n^{1-s-w}}\bigg)\frac{P_4^w\omega_1(w)\,dw}{w}\\
=&-Z\sp N(1-s,\bp)+O(E_1\sp).
\endaligned\eqno(6.3)
$$

Note that $Pt_0/P_4=T^2$. The left side of (6.3) is equal to
$I'+I{''}$
with
$$
I'=Z\sp\cdot\tp\int_{-1-i\l^{20}}^{-1+i\l^{20}}\bigg(\sum_{n
<T^3}\frac{\bp(n)}{n^{1-s-w}}\bigg)\frac{T^{-2w}\omega_1(w)\,dw}{w},\eqno(6.4)
$$
$$
I^{''}=\tp\int_{(-1)}\big(Z(s+w,\psi)-Z\sp(Pt_0)^{-w}\big)\bigg(\sum_{n
<T^2}\frac{\bp(n)}{n^{1-s-w}}\bigg)\frac{P_4^w\omega_1(w)\,dw}{w}.\eqno(6.5)
$$
Replacing the segment $u=-1$, $|v|\le\l^{20}$ by $u=-1$ and using the change of variable $w\to-w$, we obtain
$$
I'=-Z\sp N(1-s,\bp)+O(\ep).
$$
On the other hand, moving the segment   $u=-1$, $|v|\le\l^{20}$ to $u=0$, $|v|\le\l^{20}$ and applying Lemma 5.1, we find the 
the right side of (6.5) is $\ll E_1\sp$. $\Box$

\medskip\noindent

\section  *{\centerline{7. Mean-value formula I}}
\medskip\noindent

\medskip\noindent

Let $\n(d)$ denote the set of positive integers such that $h\in\n(d)$ if and only if every prime factor of $h$ divides $d$ (note that $1\in\n(d)$ for every $d$ and $\n(1)=\{1\}$). Assume $1\le j\le 3$ in what follows. Write
$$
\tk(d;m,s)=\sum_{\substack{h\in\n(d)\\(h,m)=1}}\frac{\kappa(dh)}{h^s},
$$
and
$$
 \la(m,s)=\prod_{q|m}\frac{(1-q^{-s-\be_1})(1-q^{-s-\be_2})(1-q^{-s-\be_3})}{1-q^{-s}}.
$$
For notational simplicity we write
$$
\tk(d;m,1-\be_j)=\tk_{0j}(d;m),\quad \la(m,1-\be_j)=\la_{0j}(m).
$$
Let
$$
\xi_{0j}(n;d,r)=\tilde{\la}_{0j}(n,dr)\sum_{\substack{n=d_1k\\(k,r)=1}}\frac{\tk_{0j}(d_1;drk)\mu(k)k^{1-\be_j}}{\varphi(k)}
$$ 
with
 $$
 \tilde{\la}_{0j}(n,dr)=\prod_{\substack{q|n\\(q,dr)=1}} \la_{0j}(q),
 $$
For $\psi(\bmod\,p)\in\Psi$ write
$$
\c\sp=-i(pt_0)^{\be_3}Z\sp^{-1}\frac{L(s+\be_1,\psi)L(s+\be_2,\psi)L(s+\be_3,\psi)}{L\sp}.\eqno(7.1)
$$
Let $\bold{a}=\{a(n)\}$ denote  a sequence of complex numbers satisfying
$$
a(n)\ll 1,\quad a(n)=0\quad\text{if}\quad n\ge PT^{-2}.\eqno(7.2)
$$
Write
$$
\bar{\bold{a}}=\{\bar{a}(n)\},\quad \bar{a}(n)=\overline{a(n)},
$$
$$
A(\bold{a};s,\psi)=\sum_n\frac{a(n)\psi(n))}{n^s},\quad A(\bold{a};1-s,\bp)=\sum_n\frac{a(n)\bp(n)}{n^{1-s}}.
$$
Let $\j(z)$ denote the segment $[s_0+z-i\l_1,\,s_0+z+i\l_1]$.

The goal of this section is to prove
\medskip\noindent

{\bf Proposition 7.1.}\,\,{\it Assume $\bold{a}_1$ and $\bold{a}_2$ satisfy (7.2).
We have
$$\aligned
\Theta_1(\bold{a}_1,\bold{a}_2):&=\sum_{\psi\in\Psi_1}\,\frac{1}{2\pi i}\int_{\j(1)}\c\sp A(\bold{a}_1;s,\psi)A(\bold{a}_2;1-s,\bp)\omega(s)\,ds\\=&\frac{1}{\al}\bigg(
\frac12S_{1}(\bold{a}_1,\bold{a}_2)
+2S_{2}(\bold{a}_1,\bold{a}_2)+\frac32S_{3}(\bold{a}_1,\bold{a}_2)\bigg)\p\\
&+O(E(\bold{a}_1, \bold{a}_2))+o(\p),
\endaligned
$$
where
$$
\p=\sum_{p\sim P}p,
$$
$$
S_{j}(\bold{a}_1,\bold{a}_2)=\sum_{d}\sum_{r}\frac{|\mu(r)|\la_{0j}(dr)}{dr\varphi(r)}\bigg(\sum_m\frac{a_1(drm)}{m^{1-\be_j}}\bigg)\bigg(\sum_n\frac{a_2(drn)\xi_{0j}(n;d,r)}{n}\bigg),
$$ 
and}
$$
E(\bold{a}_1, \bold{a}_2)=\p\l^2\sum_{ 1\le j\le 3}|S_{j}(\bold{a}_1,\bold{a}_2)|.
$$

 In this and the next two sections we assume that  $1\le j\le 3$.
\medskip\noindent

{\it Proof of Proposition 7.1: Initial steps}
\medskip\noindent

 For notational simplicity we write $\Theta_1$ for $\Theta_1(\bold{a}_1,\bold{a}_2)$ . Note that for any $\psi\in\Psi$, $\c\sp$ is analytic if $\sigma>1$ and $|t-2\pi t_0|\le\l_1$.  We need to show that the sum  over $\Psi_1$ in  the expression for $\Theta_1$ can  be extend to the sum over $\Psi$ with an acceptable error. Namely we prove
$$
\sum_{\psi\in\Psi_2}\int_{\j(1)}\c\sp A(\bold{a}_1;s,\psi)A(\bold{a}_2;1-s,\bp)\omega(s)\,ds=o(\p).\eqno(7.3)
$$
Here Proposition 2.1 is crucial.

 Let $\kappa(n)$ be given by
$$
\sum_n\frac{\kappa(n)}{n^s}=\frac{\zeta(s+\be_1)\zeta(s+\be_2)\zeta(s+\be_3)}{\zeta(s)},
$$
 and regard $a_1$ as an arithmetic function. For $\sigma>1$
$$
\frac{L(s+\be_1,\psi)L(s+\be_2,\psi)L(s+\be_3,\psi)}{L\sp}A(\bold{a}_1;s,\psi)=\sum_{m}\frac{(\kappa*a_1)(m)\psi(m)}{m^s}.
$$
Thus, for $\psi(\bmod\,p)\in\Psi$ and $s\in\j(1)$,
$$
\c\sp A(\bold{a}_1;s,\psi)=-i(pt_0)^{\be_3}Z\sp^{-1}\sum_{m}\frac{(\kappa*a_1)(m)\psi(m)}{m^s}.
$$
The sum over $m$ is split into two sums according to $m<P^2$ and $m\ge P^2$. To handle the second one we appeal to the trivial bounds
$$
\quad\sum_{m\ge P^2}\frac{(\kappa*a_1)(m)\psi(m)}{m^s}\ll P^{2(1-\sigma)}\l^c,\quad A(\bold{a}_2;1-s,\bp)\ll (PT^{-2})^\sigma
$$
for $\sigma\ge 3/2$, and
$$
Z\sp^{-1}\ll (pt_0)^{\sigma-1/2}
$$
for $3/2\le\sigma\le\l^9$, $|t-2\pi t_0|\le \l_1$. Thus, moving the segment to $\j(\l^9)$ and applying the simple estimate
$$
\int_{\j(z)}|\omega(s)\,ds|\ll 1,\quad |z|\le\l^9, \eqno(7.4)
$$
 we obtain
$$
\int_{\j(1)}Z\sp^{-1}\bigg(\sum_{m\ge P^2}\frac{(\kappa*a_1)(m)\psi(m)}{m^s}\bigg)A(\bold{a}_2;1-s,\bp)\omega(s)\,ds\ll \ep.
$$
To handle the sum over $m<P^2$  we move the path of integration to $\j(0)$. Hence
$$\aligned
\bigg|\int_{\j(1)}&\tilde{\c}\sp A(\bold{a}_1;s,\psi)A(\bold{a}_2;1-s,\bp)\sp\omega(s)\,ds\bigg|\\
&\le\int_{\j(0)}\bigg|\sum_{m<P^2}\frac{(\kappa*a_1)(m)\psi(m)}{m^s}\bigg||A(\bold{a}_2;1-s,\bp)\omega(s)\,ds|+O(\ep).
\endaligned
$$
 By (7.4), the proof of (7.3) is reduced to showing that
$$
\sum_{\psi\in\Psi_2}\bigg|\sum_{m<P^2}\frac{(\kappa*a_1)(m)\psi(m)}{m^s}\bigg||A(\bold{a}_2;1-s,\bp)|=o(\p),\quad \sigma=1/2.
\eqno(7.5)$$
Note that $(\kappa*a_1)(m)=O(\tau_5(m))$.  The left side above is, by Cauchy's inequality, 
$$\aligned
\le&\bigg(\sum_{\psi\in\Psi}\bigg|\sum_{m<P^2}\frac{(\kappa*a_1)(m)\psi(m)}{m^s}\bigg|^2\bigg)^{1/2}\bigg(\sum_{\psi\in\Psi}|A(\bold{a}_2;1-s,\bp)|^4\bigg)^{1/4}\bigg(\sum_{\psi\in\Psi_2}1\bigg)^{1/4}\\
&\ll\bigg(P^2\sum_{m<P^2}\frac{\tau_5(m)^2}{m}\bigg)^{1/2}\bigg(P^2\sum_{m<P^2}\frac{\tau_2(m)^2}{m}\bigg)^{1/4}
\bigg(\sum_{\psi\in\Psi_2}1\bigg)^{1/4}.
\endaligned
$$
This yields  (7.5) by Proposition 2.1 and (2.9).

By (7.3)  we may write
$$
\Theta_1=-i\sum_{p\sim P}(pt_0)^{\be_3}\,\sideset{}{^*}\sum_{\psi(\bmod\,p)}I_1(\psi)+o(\p)\eqno(7.6)
$$
where
$$
I_1(\psi)=\tp\int_{\j(1)}Z\sp^{-1}\bigg(\sum_{m}\frac{(\kappa*a_1)(m)\psi(m)}{m^s}\bigg)A(\bold{a}_2;1-s,\bp)\omega(s)\,ds.
$$
Assume $\psi(\bmod\,p)\in\Psi$.  By (2.4),  in the integral $I_1(\psi)$, the factor $Z(s,\psi)^{-1}$  can be replaced by $$\tau(\bp)p^{s-1}\vartheta^*(1-s),$$ and then, by a trivial bound for $\vartheta^*(1-s)\omega(s)$ on $\sigma=3/2$,  the segment $\j(1)$ can be replaced by the line $\sigma=3/2$, with negligible errors. Thus,  integration term-by--term gives
$$
I_1(\psi)=\frac{\tau(\bp)}{p}\sum_m\,\sum_n\frac{(\kappa*a_1)(m)a_2(n)\psi(m)\bp(n)}{n}\Delta_1\bigg(
\frac{m}{pn}\bigg)+O(\ep).
 $$
 This yields
$$
\sideset{}{^*}\sum_{\psi(\bmod\,p)}I_1(\psi)=\frac{1}{p}\sum_m\,\sum_n\frac{(\kappa*a_1)(m)a_2(n)}{n}\Delta_1\bigg(
\frac{m}{pn}\bigg)\bigg( \,\,\sideset{}{^*}\sum_{\psi(\bmod\,p)}\tau(\bp)\psi(m)\bp(n)\bigg)+O(\ep).
 $$
 Note that $(n,p)=1$ if $n<PT^{-2}$. On substituting $d=(m,n)$, $m=dl$, $n=dk$, we find that the main term on the right side above is 
 $$
 \frac{1}{p}\sum_d\,\frac{1}{d}\sum_l\,\sum_{(k,l)=1}\frac{(\kappa*a_1)(dl)a_2(dk)}{k}
\bigg( \,\,\sideset{}{^*}\sum_{\psi(\bmod\,p)}\tau(\bp)\psi(l)\bp(k)\bigg)\Delta_1\bigg(\frac{l}{pk}\bigg).
$$
 Since $\tau(\psi_p^0)=-1$, for $(kl,p)=1$ we have
  $$
\sideset{}{^*}\sum_{\psi (\bmod p)}\tau(\bp)\psi(l)\bp(k)=pe\bigg(\frac{l\bar{k}}{p}\bigg)+O(1)
$$
with $\bar{k}k\equiv 1(\bmod\,p)$. Hence, by Lemma 5.4,
$$
\sideset{}{^*}\sum_{\psi(\bmod\,p)}I(\psi)=\sum_{d}\frac{1}{d}\sum_{(l,p)=1}\,\sum_{(k,l)=1}\frac{(\kappa*a_1)(dl)b(dk)}{k}e\bigg(\frac{l\bar{k}}{p}\bigg)
\Delta_1\bigg(\frac{l}{pk}\bigg)+O( PT^{-c}).\eqno(7.7)
$$
 By trivial estimation, this  remains valid if  the constraint $(l,p)=1$ is removed. Further, by  the relation
$$
\frac{\bar{k}}{p}\equiv-\frac{\bar{p}}{k}+\frac{1}{pk}\qquad(\bmod \,1)
$$
we have
$$
e\bigg(\frac{l\bar{k}}{p}\bigg)\Delta_1\bigg(\frac{l}{pk}\bigg)=e\bigg(-\frac{l\bar{p}}{k}\bigg)\Delta\bigg(\frac{l}{pk}\bigg)
$$
Thus  the right side of (7.7) is
$$
\sum_{d}\frac{1}{d}\sum_l\,\sum_{(k,l)=1}\frac{(\kappa*a_1)(dl)b(dk)}{k}e\bigg(-\frac{l\bar{p}}{k}\bigg)\Delta\bigg(\frac{l}{pk}\bigg)
+O( PT^{-c}).\eqno(7.8)$$
For  $(l,k)=1$ we have
$$
e\bigg(\frac{-l\bar{p}}{k}\bigg)=\frac{\mu(k)}{\varphi(k)}+\frac{1}{\varphi(k)}\sideset{}{'}\sum_{\theta (\bmod k)}\tau(\bt)\theta(l)\bt(-p).$$
Inserting this into  (7.8) we deduce that
$$
\sideset{}{^*}\sum_{\psi(\bmod\,p)}I(\psi)=\t_{11}(p)+\t_{12}(p)+O(PT^{-c})\eqno(7.9)
$$
where
$$
\t_{11}(p)=\sum_{d}\frac{1}{d}\sum_l\,\sum_{(k,l)=1}\frac{(\kappa*a_1)(dl)a_2(dk)\mu(k)}{k\varphi(k)}\Delta\bigg(\frac{l}{pk}\bigg)
$$
and
$$
\t_{12}(p)=\sum_{d}\frac{1}{d}\sum_l\,\sum_{(k,l)=1}\frac{(\kappa*a_1)(dl)a_2(dk)}{k\varphi(k)}\bigg(\,\sideset{}{'}\sum_{\theta (\bmod k)}\tau(\bt)\theta(l)\bt(-p)\bigg)\Delta\bigg(\frac{l}{pk}\bigg).
$$

Note that $(pt_0)^{\be_3}=-1+O(\al_1)$. By (7.6) and (7.9), the proof of Proposition 7.1 is reduced to showing that
$$\aligned
i\t_{11}(p)=&\frac{p}{\al}\bigg(
\frac12S_{1}(\bold{a}_1,\bold{a}_2)
+2S_{2}(\bold{a}_1,\bold{a}_2)+\frac32S_{3}(\bold{a}_1,\bold{a}_2)\bigg)\\
&+O\bigg(P\l^2\sum_{1\le j\le 3}|S_{j}(\bold{a}_1,\bold{a}_2)|\bigg)+o(P)\endaligned\eqno(7.10)
$$
for $p\sim P$, and
$$
\sum_{p\sim P}p^{\be_3}\t_{12}(p)=o(\p).\eqno(7.11)
$$

{\it   Proof of Proposition 7.1: The error term}

\medskip\noindent

In this subsection we prove (7.11).

Changing the order of summation gives
$$\aligned
\sum_{p\sim P}p^{\be_3} \t_{12}(p)=&\sum_{d}\frac{1}{d}\sum_{k}\,\frac{a_2(dk)}{k\varphi(k)}\sideset{}{'}\sum_{\theta (\bmod k)}\tau(\bt)\sum_{l}(\kappa*a_1)(dl)\theta(l)\\
&\times\sum_{p\sim P}p^{\be_3}\bt(-p)\,\Delta\bigg(\frac{l}{pk}\bigg).
\endaligned
\eqno(7.12)
$$
 If $k=hr$, $r>1$ and $\theta(\bmod\, k)$ is induced by a primitive character $\theta^*(\bmod\,r)$, then $|\tau(\bt)|\le r^{1/2}$. Thus, the inner sum over the non-principal $\theta(\bmod\,k)$  above is
$$
\ll \sum_{\substack{hr=k\\r>1}}r^{1/2}\,\,\sideset{}{^*}\sum_{\theta (\bmod r)}|\s(r,h,d;\theta)|,$$
where 
$$
\s(r,h,d;\theta)=\sum_{(l,h)=1}(\kappa*a_1)(dl)\theta(l)\sum_{p\sim P}p^{\be_3}\bt(p)\Delta\bigg(\frac{l}{phr}\bigg).
$$
Inserting this into (7.12)  we obtain
$$
\sum_{p\sim P}p^{\be_3} \t_{12}(p)\ll\sideset{}{'}\sum_{\{d,h,r\}}\frac{1}{dh\varphi(hr)\sqrt{r}}\,\sideset{}{^*}\sum_{\theta (\bmod r)}|\s(r,h,d;\theta)|,\eqno(7.13)
$$
where $\sideset{}{'}\sum_{\{d,h,r\}}$ denotes a sum over the triples $\{d,h,r\}$ satisfying  $dhr<P_1$ and $r>1$. By Lemma 5.3, for $dhr<P_1$, the terms in $\s(r,h,d;\theta)$ with $l>P^2$ make a negligible error. Also, by the Mellin transform and Lemma 5.4 (i),
$$
\sum_{p\sim P}p^{\be_3}\bt(p)\Delta\bigg(\frac{l}{phr}\bigg)\ll\l^chrl^{-1}\int_{-\infty}^{\infty}\bigg|\sum_{p\sim P}p^{1+it+\be_3}\bt(p)\bigg|\frac{dt}{1+t^2}\eqno(7.14)
$$
for $l\le P^2$.

Assume $1<r<D$ and $\theta$ is a primitive character $(\bmod\,r)$. By Lemma 5.6,  the right side of (7.14) is
$$
\ll\frac{hr}{l}P^2D^{-c}.
$$
Hence
$$
\s(r,h,d;\theta)\ll\tau_5(d)hrP^2D^{-c}.
$$
Thus, on the right side of (7.13), the total contribution from the terms with $1<r<D$ is $O(P^2D^{-c})$. The proof of 
(7.11) is therefore reduced to showing that
  $$
 \frac{1}{R^{3/2}}\sum_{R\le r<2R}\,\,\sideset{}{^*}\sum_{\theta (\bmod r)}|\s(r,h,d;\theta)|\ll \tau_5(d)hP^2D^{-c} \eqno(7.15)
 $$
 for 
 $$dh<P_1\quad\text{and}\quad D\le R<(dh)^{-1}P_1
 $$
  which are henceforth assumed .

 Let $\mathfrak{I}(y)$ denote the interval
 $$
 [(1/3)Pt_0y,\,4Pt_0y].
 $$
 By Lemma 5.1,  for $R\le r<2R$,  the terms with 
 $l\notin\mathfrak{I}(Rh)$ 
 in $\s(r,h,d;\theta)$  make a negligible contribution.  Thus, it suffices to prove  (7.15) with   $\s(r,h,d;\theta)$ replaced by
 $$
\s^*(R,r,h,d;\theta):=\sum_{\substack{l\in\mathfrak{I}(Rh)\\(l,h)=1}}(\kappa*a_1)(dl)\theta(l)\sum_{p\simeq P}\bt(p)p^{\be_3}\Delta\bigg(\frac{l}{phr}\bigg).
$$
By the Mellin transform and Lemma 5.4 (i),
$$
\s^*(R,r,h,d;\theta)\ll  \l^c\int_{-\infty}^{\infty}\bigg|\sum_{\substack{l\in\mathfrak{I}(Rh)\\(l,h)=1}}(\kappa*a_1)(dl)\theta(l)l^{-1-it}\bigg|\bigg|\sum_{p\simeq P}\bt(p)p^{1+it+\be_3}\bigg|\frac{dt}{1+t^2}.
$$
For $\sigma=1$, by the large sieve inequality we have
$$
\sum_{R\le r<2R}\,\,\sideset{}{^*}\sum_{\theta (\bmod r)}\bigg|\sum_{\substack{l\in\mathfrak{I}(Rh)\\(l,h)=1}}(\kappa*a_1)(dl)\theta(l)l^{-s}\bigg|^2\ll \tau_5(d)^2\l^c(R^2+PRht_0)(PRht_0)^{-1}\ll \tau_5(d)^2\l^c
$$
and
$$
\sum_{R\le r<2R}\,\,\sideset{}{^*}\sum_{\theta (\bmod r)}\bigg|\sum_{p\simeq P}\bt(p)p^{s}\bigg|^2\ll(R^2+P)P^3.
$$
It follows by Cauchy's inequality that 
  $$
 \frac{1}{R^{3/2}}\sum_{R\le r<2R}\,\,\sideset{}{^*}\sum_{\theta (\bmod r)}|\s^*(R,r,h,d;\theta)|\ll \tau_5(d)h\l^c(R^{1/2}P^{3/2}+R^{-1/2}P^2)\ll  \tau_5(d)hP^2D^{-c}.
$$
 This yields (7.15).
\medskip\noindent

{\it   Proof of Proposition 7.1: The main term}
\medskip\noindent

In this subsection we prove (7.10). 

Assume $p\sim P$.  We may write
$$
\t_{11}(p)=\sum_{d}\frac{1}{d}\sum_k\frac{a_2(dk)\mu(k)}{k\varphi(k)}\,\sum_{(l,k)=1}(\kappa*a_1)(dl)\Delta\bigg(\frac{l}{pk}\bigg).
\eqno(7.16)$$
Since
$$
(\kappa*a_1)(dl)=\sum_{m_1m_2=dl}\kappa(m_1)a_1(m_2)=\sum_{d=d_1d_2}\,\sum_{\substack{m_1m_2=dl\\(m_1,d)=d_1}}\kappa(m_1)a_1(m_2),
$$
 it follows, by substituting $m_1=d_1l_1$, that
$$
(\kappa*a_1)(dl)=\sum_{d=d_1d_2}\,\sum_{\substack{l=l_1l_2\\(l_1,d_2)=1}}\kappa(d_1l_1)a_1(d_2l_2).\eqno(7.17)
$$
Hence,
$$
\sum_{(l,k)=1}(\kappa*a_1)(dl)\Delta\bigg(\frac{l}{pk}\bigg)=\sum_{d=d_1d_2}\,\sum_{(l_2,k)=1}a_1(d_2l_2)\,\sum_{(l_1,kd_2)=1}
\kappa(d_1l_1)\Delta\bigg(\frac{l_1l_2}{pk}\bigg).
$$
Inserting this into (7.16) we obtain
$$
\t_{11}(p)=\sum_{d_1}\,\sum_{d_2}\frac{1}{d_1d_2}\sum_k\frac{a_2(d_1d_2k)\mu(k)}{k\varphi(k)}\,
\sum_{(l_2,k)=1}a_1(d_2l_2)\,\sum_{(l_1,d_2k)=1}
\kappa(d_1l_1)\Delta\bigg(\frac{l_1l_2}{pk}\bigg).
\eqno(7.18)$$
 The innermost sum  is, by the Mellin transform, equal to
$$
\tp\int_{(3/2)}\bigg(\sum_{(l,d_2k)=1}\frac{\kappa(d_1l)}{l^s}\bigg)\bigg(\frac{pk}{l_2}\bigg)^s\delta(s)\,ds\eqno(7.19)
$$
(here we have rewritten $l$ for $l_1$). In view of (7.2), we can assume that
$$
d_2l_2<PT^{-2},\quad d_1d_2k<PT^{-2}.
$$
 Every $l$ with $(l, d_2k)=1$ can be uniquely written as $l=hr$ such that $h\in\n(d_1)$, $(h,d_2k)=1$)and $(r, d_1d_2k)=1$.  Hence, 
 $$
 \sum_{(l,kd_2)=1}\frac{\kappa(d_1l)}{l^s}=\tk(d_1;d_2k,s)\la(d_1d_2k,s)\frac{\zeta(s+\be_1)\zeta(s+\be_2)\zeta(s+\be_3)}{\zeta(s)},
$$
By the simple bounds
$$
|\tk(d_1;d_2k,s)|\le\tau_5(d_1)\prod_{q|d_1}\bigg|1+\frac{c}{q^\sigma}\bigg|,\quad |\la(m,s)|\le\prod_{q|m}\bigg|1+\frac{c}{q^\sigma}\bigg|
$$
for $\sigma>9/10$, we can move the contour of integration in (7.19) to the vertical segments
$$\aligned
&s=1+\al+it\qquad&\text{with}\quad  |t|\ge D,\\
&s=1-\l^{-1}+it\quad&\text{with}\quad  |t|\le D
\endaligned
$$
and to the two connecting horizontal segments
$$
s=\sigma\pm iD\quad\text{with}\quad 1-\l^{-1}\le\sigma\le 1+\al.
$$
This yields, by Lemma 5.2 (i) and standard estimates,  that the integral (7.18) is equal to the sum of the residues of the integrand at $s=1-\be_j$, $1\le j\le 3$, plus an acceptable error. Namely we have
$$
\sum_{(l,kd_2)=1}\kappa(d_1l)\Delta\bigg(\frac{ll_2}{pk}\bigg)=\sum_{ 1\le j\le 3}\r_j\tk_{0j}(d_1;d_2k)\la_{0j}(d_1d_2k)\bigg(\frac{pk}{l_2}\bigg)^{1-\be_j}+O\bigg(\frac{pk\ep_1}{l_2}\bigg),
$$
where $\r_j$ is the residue of the function
$$
\frac{\zeta(s+\be_1)\zeta(s+\be_2)\zeta(s+\be_3)}{\zeta(s)}\delta(s)
$$
at $s=1-\be_j$, and where
$$
\ep_1=\exp\{-c\l^{1/10}\}.
$$
This yields
$$\aligned
\sum_{(l_2,k)=1}&a_1(d_2l_2)\,\sum_{(l_1,kd_2)=1}\kappa(d_1l_1)\Delta\bigg(\frac{l_1l_2}{pk}\bigg)\\
=&\sum_{1\le j\le 3}\r_j(pk)^{1-\be_j}\sum_{d=d_1d_2}\tk_{0j}(d_1;d_2k)\la_{0j}(d_1d_2k)
\sum_{(l,k)=1}\frac{a_1(d_2l)}{l^{1-\be_j}}+O(Pk\ep_1).
\endaligned
$$
(here we have rewritten $l$ for $l_2$). Inserting this into (7.18) and rearranging the terms we obtain
$$\aligned
\t_{11}(p)=&\sum_{1\le j\le 3}\r_jp^{1-\be_j}S_j^*(\bold{a}_1,\bold{a}_2)+O(P\ep_1)
\endaligned\eqno(7.20)
$$
where
$$
S_j^*(\bold{a}_1,\bold{a}_2)=\sum_d\,\sum_{d_1}\,\sum_k\frac{a_2(d_1dk)\tk_{0j}(d_1,dk)\la_{0j}(d_1dk)\mu(k)}{d_1d\varphi(k)k^{\be_j}}\sum_{(l,k)=1}\frac{a_1(dl)}{l^{1-\be_j}}.$$
By the M\"{o}bius inversion,
  $$\aligned
S_{j}^*(\bold{a}_1,\bold{a}_2)=&\sum_{d}\sum_{d_1}\sum_k\frac{a_2(d_1dk)\tk_{0j}(d_1;dk)\la_{0j}(d_1dk)\mu(k)}{d_1d\varphi(k)k^{\be_j}}\\&
\times\sum_{l}\frac{a_1(dl)}{l^{1-\be_j}}\bigg(\sum_{r|(k,l)}\mu(r)\bigg).
\endaligned
$$
This yields, by substituting $k=rk_1$,  $l=rm$, $n=d_1k_1$ and changing the order of summation,
 $$\aligned
S_{j}^*(\bold{a}_1,\bold{a}_2)=\sum_r\,\sum_d\frac{|\mu(r)|}{dr\varphi(r)}\sum_{m}\frac{a_1(drm)}{m^{1-\be_j}}\sum_n\frac{a_2(drn)\la_{0j}(drn)}{n}\sum_{\substack{n=d_1k_1\\(k_1,r)=1}}\frac{\tk_{0j}(d_1,drk_1)\mu(k_1)k_1^{1-\be_j}}{\varphi(k_1)}.
\endaligned$$
Since
   $$
   \la_{0j}(drn)=\la_{0j}(dr)\tilde{\la}_{0j}(n,dr),
 $$
 it follows that
 $$
 \la_{0j}(drn)\sum_{\substack{n=d_1k_1\\(k_1,r)=1}}\frac{\tk_{0j}(d_1,drk_1)\mu(k_1)k_1^{1-\be_j}}{\varphi(k_1)}=\la_{0j}(dr)
 \xi_{0j}(n;d,r).
$$
Hence
 $$
 S_{j}^*(\bold{a}_1,\bold{a}_2)=S_{j}(\bold{a}_1,\bold{a}_2).\eqno(7.21)
 $$
 On the other hand, by Lemma 5.2 (ii) and direct calculation we have
 $$
i\r_1p^{-\be_1}=\frac{1}{2\al}+O(\l),
$$
$$
i\r_2p^{-\be_2}=\frac{2}{\al}+O(\l),
$$
$$
i\r_3p^{-\be_3}=\frac{3}{2\al}+O(\l).
$$
 Combining  these with (7.20) and (7.21) we obtain (7.10), and complete the proof of Proposition 7.1. $\Box$

  \section  *{\centerline{8. Evaluation of $\Xi_{11}$}}
  
  \medskip\noindent
  
  We first prove a general result as follows.
  
  {\bf Lemma 8.1.}\,\,{\it For any $\bold{a}_1$ and $\bold{a}_2$ satisfying (7.2),}
 $$
 \sum_{\psi\in\Psi_1}\,\sum_{\rho\in\tz(\psi)}\c^*\ss A(\bold{a}_1;\rho,\psi)A(\bold{a}_2,1-\rho,\bp)\omega(\rho)=\Theta_1(\bold{a}_1,\bold{a}_2)+\overline{\Theta_1(\bar{\bold{a}}_2, \bar{\bold{a}}_1)}+o(\p).
 $$
  \medskip\noindent
  
  {\it Proof.}\,\,Write
  $$
  \tilde{\c}\sp=-i\frac{M(s+\be_1,\psi)M(s+\be_2,\psi)M(s+\be_3,\psi)}{M(s,\psi)}.
  $$
  Assume $\psi\in\Psi_1$. By Proposition 2.2, we can choose a rectangle  $\mathfrak{R}$ with vertices at
$$
s_0\pm\al+i\l_1^\pm
$$
such that $\l_1^\pm=\pm\l_1+O(\al)$ and such that the set of zeros of $L\sp$ inside $\mathfrak{R}$ is exactly $\tz(\psi)$. Further, without loss of generality, we can assume that $|s-\rho|\gg\al$  if $s\in\mathfrak{R}$ and  $\rho$ is a zero of $L\sp$. By Lemma 5.9, the residue theorem and a simple bound for $\omega(s)$,
$$\aligned
\sum_{\rho\in\tz(\psi)}&\c^*\ss A(\bold{a}_1;\rho,\psi)A(\bold{a}_2,1-\rho,\bp)\omega(\rho)\\&=\tp\int_{\mathfrak{R}}\tilde{\c}\sp A(\bold{a}_1;s,\psi)A(\bold{a}_2,1-s,\bp)\omega(s)\,ds
\\&=\tilde{I}_1^+(\bold{a}_1,\bold{a}_2;\psi)-\tilde{I}_1^-(\bold{a}_1,\bold{a}_2;\psi)+O(\ep)
\endaligned\eqno(8.1)
$$
where
$$
\tilde{I}_1^\pm(\bold{a}_1,\bold{a}_2;\psi)=\frac{1}{2\pi i} \int_{\j(\pm \al)}\tilde{\c}\sp A(\bold{a}_1;s,\psi)A(\bold{a}_2,1-s,\bp)\omega(s)\,ds.
$$
Write $s'=1-\bar{s}$. If
$$
s=\al+s_0+iv\in\j(\al).
$$
then
$$
s'=-\al+s_0+iv\in\j(-\al).
$$
By (2.11) and analytic continuation, for $s\in\j(\al)$  we have
$$
\overline{\tilde{\c}\sp}=-\tilde{\c}(s',\psi)
$$
and
$$
\overline{A(\bold{a}_1;s,\psi)A(\bold{a}_2,1-s,\bp)\omega(s)}=A(\bar{\bold{a}}_2;s',\psi)A(\bar{\bold{a}}_1,1-s',\bp)\omega(s').
$$
These together imply that
$$
-\tilde{I}_1^-(\bold{a}_1,\bold{a}_2;\psi)=\overline{\tilde{I}_1^+(\bar{\bold{a}}_2,\bar{\bold{a}}_1;\psi)}.
$$
 Hence, by (8.6),  
$$
\sum_{\rho\in\tz(\psi)}\c^*\ss A(\bold{a}_1;\rho,\psi)A(\bold{a}_2,1-\rho,\bp)\omega(\rho)=\tilde{I}_1^+(\bold{a}_1,\bold{a}_2;\psi)
+\overline{\tilde{I}_1^+(\bar{\bold{a}}_2,\bar{\bold{a}}_1;\psi)}+O(\ep).
$$
Write
$$
I_1^+(\bold{a}_1,\bold{a}_2;\psi)=\frac{1}{2\pi i} \int_{\j( \al)}\c\sp A(\bold{a}_1;s,\psi)A(\bold{a}_2,1-s,\bp)\omega(s)\,ds.
$$
 By Lemma 5.2 and 5.9,
 $$
 \tilde{I}_1^+(\bold{a}_1,\bold{a}_2;\psi)-I_1^+(\bold{a}_1,\bold{a}_2;\psi)\ll\l^{-114} \int_{\j( \al)}\bigg|L(s+\be_2,\psi)L(s+\be_3,\psi) A(\bold{a}_1;s,\psi)A(\bold{a}_2,1-s,\bp)\omega(s)\,ds\bigg|.
 $$
  By Cauchy's inequality, Lemma 6.1, and the second assertion of Lemma 3.3, for $s\in\j(\al)$,
  $$
  \sum_{\psi\in\Psi_1}|L(s+\be_2,\psi)L(s+\be_3,\psi)|^2\ll P^2\l^{36},
  $$
   $$
 \sum_{\psi\in\Psi_1}|A(\bold{a}_1;s,\psi)A(\bold{a}_2,1-s,\bp)|^2\ll P^2\l^{36}.
  $$
  These estimates together with (7.4) and (2.9)  imply
  $$
  \sum_{\psi\in\Psi_1}\big( \tilde{I}_1^+(\bold{a}_1,\bold{a}_2;\psi)-I_1^+(\bold{a}_1,\bold{a}_2;\psi)\big)=o(\p).
  $$
  On the other hand, moving the segment $\j(\al)$ to $\j(1)$ gives
  $$
   \sum_{\psi\in\Psi_1}I_1^+(\bold{a}_1,\bold{a}_2;\psi)=\Theta_1(\bold{a}_1,\bold{a}_2)+O(\ep).
   $$
   Hence
   $$
   \sum_{\psi\in\Psi_1}\tilde{I}_1^+(\bold{a}_1,\bold{a}_2;\psi)=\Theta_1(\bold{a}_1,\bold{a}_2)+o(\p).
   $$
  For the sum of $\overline{\tilde{I}_1^+(\bar{\bold{a}}_2,\bar{\bold{a}}_1;\psi)}$ a similar result holds. This completes the proof. $\Box$

   Recall that $H_1\sp$ and $H_2\sp$ are given by (2.27). Write
$$
\Xi_1=\sum_{\psi\in\Psi_1}\,\sum_{\rho\in\z(\psi)}\c^*\ss \big|H_1\ss+Z(\rho,\pc)\overline{H_2\sp}\big|^2\omega(\rho).
$$
Since $|Z(s,\pc)|=1$ if $\sigma=1/2$, it follows that
$$
\Xi_1=\Xi_{11}+\Xi_{12}+2\Re\{\Xi_{13}\}\eqno(8.2)
$$
where
$$
\Xi_{11}=\sum_{\psi\in\Psi_1}\,\sum_{\rho\in\z(\psi)}\c^*\ss|H_1\ss|^2\omega(\rho),\eqno(8.3)
$$
$$
\Xi_{12}=\sum_{\psi\in\Psi_1}\,\sum_{\rho\in\z(\psi)}\c^*\ss |{H_2\ss}|^2\omega(\rho),\eqno(8.4)
$$
$$
\Xi_{13}=\sum_{\psi\in\Psi_1}\,\sum_{\rho\in\z(\psi)}\c^*\ss Z(\rho,\pc)^{-1}H_1\ss H_2\ss\omega(\rho).\eqno(8.5)
$$

 We may write
  $$
H_{1j}\sp=\sum_n\frac{\vk_j(n)\pc(n)}{n^s}\eqno(8.6)
$$
where
$$
\vk_1(n)=
\begin{cases}
(1-\log n/\log P_1)(P_1/n)^{\be_6}&\quad\text{if}\quad n<P_1,\\
0&\quad\text{otherwise},
\end{cases}
$$
$$
\vk_2(n)=
\begin{cases}
(1-\log n/\log P_2)(P_2/n)^{\be_7}&\quad\text{if}\quad n<P_2,\\
0&\quad\text{otherwise},
\end{cases}
$$
$$
\vk_3(n)=
\begin{cases}
(1-\log n/\log P_3)(P_3/n)^{\be_6}&\quad\text{if}\quad n<P_3,\\
0&\quad\text{otherwise}.
\end{cases}
$$

  In this section we evaluate $\Xi_{11}$. By Lemma 8.1,
$$
\Xi_{11}=2\Re\{\Theta_1(\bold{a}_{11},\bold{a}_{21})\}+o(\p)\eqno(8.7)
$$
with 
$$
a_{11}(n)=\chi(n)\big(\vk_1(n)+\iota_2\vk_2(n)\big),\quad a_{21}(n)=\overline{a_{11}(n)}.\eqno(8.8)
$$
By Proposition 7.1, our goal  is reduced to evaluating the sum
  $$\aligned
S_{j}(\bold{a}_{11},\bold{a}_{21})=\sum_d\,\sum_{r}&\frac{|\chi(d)||\mu\chi(r)|\la_{0j}(dr)}{dr\varphi(r)}\bigg(\sum_m\frac{\chi(m)\big(\vk_1(drm)+\iota_2\vk_2(drm)\big)}{m^{1-\be_j}}\bigg)\\
&\times\bigg(\sum_n\frac{\chi(n)\big(\bar{\vk}_1(drn)+\bar{\iota}_2\bar{\vk}_2(drn)\big)\xi_{0j}(n;d,r)}{n}\bigg).
\endaligned
$$

 Write
   $$
   \be_4=\be_1,\quad \be_5=\be_2,
   $$
   so that
   $$
   \{\be_j,\,\be_{j+1},\,\be_{j+2}\}=\{\be_1,\,\be_2,\,\be_3\}.
   $$
 \medskip\noindent
 
 {\bf Lemma 8.2.}\,\,{\it Suppose $T<x<P$. Then for $\mu=6, 7$,
 $$
 \sum_{m<x}\frac{\chi(m)}{m^{1-\be_j}}\bigg(\frac{x}{m}\bigg)^{\be_\mu}\log\frac{x}{m}=L'(1,\chi)\f_{j\mu}(x)+O(\l^{-6})
 $$
 where}
 $$
\f_{j\mu}(x)=\big(1+(\be_\mu-\be_j)\log x \big)x^{\be_\mu}.
$$
 
 {\it Proof.}\,\, The sum is equal to
 $$
 \tp\int_{(1)}L(1-\be_j+s,\chi)\frac{x^s\,ds}{(s-\be_\mu)^2}.
 $$
 We move the contour of integration to the vertical segments
$$\aligned
&s=\al+it\qquad&\text{with}\quad  |t|\ge D,\\
&s=-\l^{-1}+it\quad&\text{with}\quad  |t|\le D
\endaligned
$$
and to the two connecting horizontal segments
$$
s=\sigma+\pm iD,\quad\text{with}\quad -\l^{-1}\le\sigma\le \al.
$$
It follows by Lemma 5.6 that
$$\aligned
\sum_{m<x}\frac{\chi(m)}{m^{1-\be_j}}\bigg(\frac{x}{m}\bigg)^{\be_\mu}\log\frac{x}{m}=& \tp\int_{|s|=5\al}L(1-\be_j+s,\chi)\frac{x^s\,ds}{(s-\be_\mu)^2}+O(\ep_1)\\
=&L'(1,\chi)\cdot \tp\int_{|s|=5\al}\frac{x^s(s-\be_j)}{(s-\be_\mu)^2}\,ds+O(\l^{-6}).
\endaligned
$$
The result now follows by direct calculation. $\Box$
 
 The sum involving $\xi_{0j}(n;d,r)$ is more involved. We need the following lemma that will be proved in Appendix A.

\medskip\noindent

{\bf Lemma 8.3.}\,\,{\it Suppose  $dr<PT^{-2}$. The function
 $$\aligned
\u_{j}(d,r;s):=\frac{L(s,\chi)}{L(s+\be_{j+1},\chi)L(s+\be_{j+2},\chi)}\sum_{n}\frac{\chi(n)\xi_{0j}(n;d,r)}{n^s}
\endaligned
$$
 is analytic  and it satisfies
$$
|\u_{j}(d,r;s)|<c\prod_{q|dr}\big(1+cq^{-\sigma}\big)
$$
 for $\sigma>9/10$. Further, if  $|s-1|\le 5\al$, then
$$
\u_{j}(d,r;s)=\Pi(d,r)+O(\l^{-8})
$$
 where}
$$
\Pi(d,r)=\prod_{q|dr}\frac{1}{1-\chi(q)q^{-1}}\prod_{\substack{(q,r)=1\\q|d}}\frac{1-q^{-1}-\chi(q)q^{-1}}{1-q^{-1}}.
$$
 
 {\bf Lemma 8.4.}\,\,{\it Suppose $dr<PT^{-2}$ and $T<x<P$. Then for $\mu=6,7$,
 $$
\sum_{n<x}\frac{\chi(n)\xi_{0j}(n;d,r)}{n}\bigg(\frac{x}{n}\bigg)^{-\be_\mu}\log \frac{x}{n}=L'(1,\chi)\Pi(d,r)\g_{j\mu}(x)+O(\l^{-6})
$$
where}
$$
\g_{j\mu}(x)=\frac{\be_{j+1}\be_{j+2}}{\be_\mu^2}+\bigg(1-\frac{\be_{j+1}\be_{j+2}}{\be_\mu^2}-\frac{(\be_{j+1}-\be_\mu)(\be_{j+2}-\be_\mu)}{\be_\mu}\log x\bigg)x^{-\be_\mu}.
$$

{\it Proof.}\,\, The sum is equal to
$$
\tp\int_{(1)}\bigg(\sum_n\frac{\chi(n)\xi_{0j}(n;d,r)}{n^{1+s}}\bigg)\frac{x^s\,ds}{(s+\be_\mu)^2}.\eqno(8.9)
$$
The contour of integration is moved in the same way as in the proof of Lemma 8.2. Recall that 
$$
\sum_n\frac{\chi(n)\xi_{0j}(n;d,r)}{n^{1+s}}=\frac{L(1+s+\be_{j+1},\chi)L(1+s+\be_{j+2},\chi)}{L(1+s,\chi)}\u_j(d,r,1+s).
$$
 By Lemma 5.5 and 5.6, for $|s|=5\al$ we have
$$
\frac{L(1+s+\be_{j+1},\chi)L(1+s+\be_{j+2},\chi)}{L(1+s,\chi)}=L'(1,\chi)\frac{(s+\be_{j+1})(s+\be_{j+2})}{s}+O(\l^{-15}).
$$
Combining these results with Lemma 8.3, we find that the integral (8.9) is equal to
$$
L'(1,\chi)\Pi(d,r)\cdot \tp\int_{|s|=5\al}\frac{(s+\be_{j+1})(s+\be_{j+2})}{s}\frac{x^s\,ds}{(s+\be_\mu)^2}+O(\l^{-6}).
$$
The result now follows by direct calculation. $\Box$

  By Lemma 8.2 with $x=P_1/dr$ and  $x=P_2/dr$ respectively we have
$$
\sum_m\frac{\chi(m)\vk_1(drm)}{m^{1-\be_j}}=\frac{L'(1,\chi)}{\log P_1}\f_{j6}(P_1/dr)+O(\l^{-15})
$$
 if $dr<P_1/T$, and 
$$
\sum_m\frac{\chi(m)\vk_2(drm)}{m^{1-\be_j}}=\frac{L'(1,\chi)}{\log P_2}\f_{j7}(P_2/dr)+O(\l^{-15})
$$
 if $dr<P_2/T$.  On the other hand, by Lemma 8.3 with $x=P_1/dr$ and  $x=P_2/dr$ respectively we have
$$
\sum_n\frac{\chi(n)\bar{\vk}_1(drn)\xi_{0j}(n;d,r)}{n}=\frac{L'(1,\chi)\Pi(d,r)}{\log P_1}\g_{j6}(P_1/dr)+O(\l^{-15})
$$
if $dr<P_1/T$, and
$$
\sum_n\frac{\chi(n)\bar{\vk}_2(drn)\xi_{0j}(n;d,r)}{n}=\frac{L'(1,\chi)\Pi(d,r)}{\log P_2}\g_{j7}(P_2/dr)+O(\l^{-15})
$$
if $dr<P_2/T$. Gathering these results together we conclude, by simple approximation, that
$$\aligned
S_{j}(\bold{a}_{11},\bold{a}_{21})=&L'(1,\chi)^2\sum_{dr<P_2}\frac{|\mu(r)\chi(dr)|}{dr\varphi(r)}\la_{0j}(dr)\Pi(d,r)\\
&\quad\times\bigg(\frac{\f_{j6}(P_1/dr)}{\log P_1}+\iota_2\frac{\f_{j7}(P_2/dr)}{\log P_2}\bigg)\bigg(\frac{\g_{j6}(P_1/dr)}{\log P_1}+\bar{\iota}_2\frac{\g_{j7}(P_2/dr)}{\log P_2}\bigg)\\
&+L'(1,\chi)^2\sum_{P_2\le dr<P_1}\frac{|\mu(r)\chi(dr)|}{dr\varphi(r)}\la_{0j}(dr)\Pi(d,r)\frac{\f_{j6}(P_1/dr)\g_{j6}(P_1/dr)}{(\log P_1)^2}+o(\al).
\endaligned
$$
It can be shown, by verifying the case $n=q^k$, that
$$
\sum_{n=dr}\frac{|\mu(r)|}{r}\prod_{\substack{(q,r)=1\\q|d}}(1-q^{-1}-\chi(q)q^{-1})=\prod_{q|n}(1-\chi(q)q^{-1}),
$$
so that
$$
\sum_{n=dr}\frac{|\mu(r)|}{\varphi(r)}\Pi(d,r)=\frac{n}{\varphi(n)}.\eqno(8.10)
$$
It follows, by substituting $n=dr$,  that
$$\aligned
S_{j}(\bold{a}_{11},\bold{a}_{21})=&L'(1,\chi)^2\sum_{n<P_2}\frac{|\chi(n)|\la_{0j}(n)}{\varphi(n)}\\
&\quad\times\bigg(\frac{\f_{j6}(P_1/n)}{\log P_1}+\iota_2\frac{\f_{j7}(P_2/n)}{\log P_2}\bigg)\bigg(\frac{\g_{j6}(P_1/n)}{\log P_1}+\bar{\iota}_2\frac{\g_{j7}(P_2/n)}{\log P_2}\bigg)\\
&+L'(1,\chi)^2\sum_{P_2\le n<P_1}\frac{|\chi(n)|\la_{0j}(n)}{\varphi(n)}\frac{\f_{j6}(P_1/n)\g_{j6}(P_1/n)}{(\log P_1)^2}+o(\al).
\endaligned
$$
Since for $n<P$,
$$
\la_{0j}(n)=\prod_{q|n}(1-q^{-1})^2(1+O(\al\log q/q))=\frac{\varphi(n)^2}{n^2}+O(\al_1)
$$
and, for $x<P$, 
$$
\sum_{n<x}\frac{|\chi(n)|\varphi(n)}{n^2}=\frac{\varphi(D)}{D}\bigg(\prod_{(q,D)=1}(1-q^{-2})\bigg)\log x+O(\log\l)=\frac{6}{\pi^2}\bigg(\prod_{q|D}\frac{q}{q+1}\bigg)\log x+O(\log\l),
$$
(see [T, 1.2.12]), it follows by partial integration that
$$\aligned
S_{j}(\bold{a}_{11},\bold{a}_{21})=&\aa\int_{1}^{P_2}\bigg(\frac{\f_{j6}(P_1/x)}{\log P_1}+\iota_2\frac{\f_{j7}(P_2/x)}{\log P_2}\bigg)\bigg(\frac{\g_{j6}(P_1/x)}{\log P_1}+\bar{\iota}_2\frac{\g_{j7}(P_2/x)}{\log P_2}\bigg)\,\frac{dx}{x}\\
&+\frac{\aa}{(\log P_1)^2}\int_{P_2}^{P_1}\frac{\f_{j6}(P_1/x)\g_{j6}(P_1/x)}{x}\,dx+o(\al).
\endaligned\eqno(8.11)
$$
This yields, by the change of variable $x\to P_1/x$ or $x\to P_2/x$,
$$\aligned
S_{j}(\bold{a}_{11},\bold{a}_{21})=&\frac{\aa}{(\log P_1)^2}\int_{1}^{P_1}\frac{\f_{j6}(x)\g_{j6}(x)}{x}\,dx+\frac{\aa|\iota_2|^2}{(\log P_2)^2}\int_{1}^{P_2}\frac{\f_{j7}(x)\g_{j7}(x)}{x}\,dx\\
&+\frac{\aa\iota_2}{(\log P_1)(\log P_2)}\int_{1}^{P_2}\frac{\f_{j7}(x)\g_{j6}(P^{0.004}T^{10}x)}{x}\,dx\\
&+\frac{\aa\bar{\iota}_2}{(\log P_1)(\log P_2)}\int_{1}^{P_2}\frac{\f_{j6}(P^{0.004}T^{10}x)\g_{j7}(x)}{x}\,dx+o(\al),
\endaligned\eqno(8.12)
$$
since $P_1/P_2=P^{0.004}T^{10}$. By direct calculation, for  $0\le z\le 1$ we have
$$
\f_{j6}(P^z)=\ff_{j6}(z)+O(\l^{-8}),\qquad \g_{j6}(P^z)=\gh_{j6}(z)+O(\l^{-8}),
$$
$$
\f_{j7}(P^z)=\ff_{j7}(z)+O(\l^{-8}),\qquad \g_{j7}(P^z)=\gh_{j7}(z)+O(\l^{-8})
$$
where
$$
\ff_{16}(z)=\bigg(1+\frac{\pi i}{2}z\bigg)e^{(3\pi i/2)z},\quad \gh_{16}(z)=\frac83+\bigg(-\frac53-\frac{\pi i}{2}z\bigg)e^{(-3\pi i/2)z},
\eqno(8.13)$$
$$
\ff_{26}(z)=\bigg(1-\frac{\pi i}{2}z\bigg)e^{(3\pi i/2)z},\quad \gh_{26}(z)=\frac43+\bigg(-\frac13+\frac{\pi i}{2}z\bigg)e^{(-3\pi i/2)z},
\eqno(8.14)$$
$$
\ff_{36}(z)=\bigg(1-\frac{3\pi i}{2}z\bigg)e^{(3\pi i/2)z},\quad \gh_{36}(z)=\frac89+\bigg(\frac19+\frac{\pi i}{6}z\bigg)e^{(-3\pi i/2)z}.
\eqno(8.15)$$
$$
\ff_{17}(z)=\bigg(1+\frac{3\pi i}{2}z\bigg)e^{(5\pi i/2)z},\quad \gh_{17}(z)=\frac{24}{25}+\bigg(\frac{1}{25}+\frac{\pi i}{10}z\bigg)e^{(-5\pi i/2)z},
\eqno(8.16)$$
$$
\ff_{27}(z)=\bigg(1+\frac{\pi i}{2}z\bigg)e^{(5\pi i/2)z},\quad \gh_{27}(z)=\frac{12}{25}+\bigg(\frac{13}{25}+\frac{3\pi i}{10}z\bigg)e^{(-5\pi i/2)z},
\eqno(8.17)$$
$$
\ff_{37}(z)=\bigg(1-\frac{\pi i}{2}z\bigg)e^{(5\pi i/2)z},\quad \gh_{37}(z)=\frac{8}{25}+\bigg(\frac{17}{25}-\frac{3\pi i}{10}z\bigg)e^{(-5\pi i/2)z}.
\eqno(8.18)$$
Substituting $x=P^z$ we obtain
$$
\frac{1}{(\log P_1)^2}\int_{1}^{P_1}\frac{\f_{j6}(x)\g_{j6}(x)}{x}\,dx=\frac{1}{(0.504)^2\log P}\int_{0}^{0.504}\ff_{j6}(z)\mathfrak{g}_{j6}(z)\,dz+o(\al),
$$
$$
\frac{1}{(\log P_2)^2}\int_{1}^{P_2}\frac{\f_{j7}(x)\g_{j7}(x)}{x}\,dx=\frac{1}{(0.5)^2\log P}\int_{0}^{0.5}\ff_{j7}(z)\mathfrak{g}_{j7}(z)\,dz+o(\al),
$$
$$\aligned
\frac{1}{(\log P_1)(\log P_2)}&\int_{1}^{P_2}\frac{\f_{j7}(x)\g_{j6}(P^{0.004}T^{10}x)}{x}\,dx\\
&=\frac{1}{(0.5)(0.504)\log P}\int_{0}^{0.5}\ff_{j7}(z)\mathfrak{g}_{j6}(z+0.004)\,dz+o(\al),
\endaligned
$$
$$\aligned
\frac{1}{(\log P_1)(\log P_2)}&\int_{1}^{P_2}\frac{\f_{j6}(P^{0.004}T^{10}x)\g_{j7}(x))}{x}\,dx\\
&=\frac{1}{(0.5)(0.504)\log P}\int_{0}^{0.5}\ff_{j6}(z+0.004)\mathfrak{g}_{j7}(z)\,dz+o(\al).
\endaligned
$$
Inserting these into (8.13) we obtain
$$
\frac{1}{2\al}S_1(\bold{a}_{11},\bold{a}_{21})+\frac{2}{\al}S_2(\bold{a}_{11},\bold{a}_{21})+\frac{3}{2\al}S_3(\bold{a}_{11},\bold{a}_{21})=\aa
\big(b_{11}+\iota_2b_{21}+\bar{\iota}_2b_{12}+|\iota_2|^2b_{22}\big)+o(1)
$$
where
$$
b_{11}=\frac{1}{0.504^2\pi}\int_0^{0.504}\bigg(\frac12\ff_{16}(z)\gh_{16}(z)+2\ff_{26}(z)\gh_{26}(z)+\frac32\ff_{36}(z)\gh_{36}(z)\bigg)\,dz,\eqno(8.19)
$$
$$
b_{22}=\frac{1}{0.5^2\pi}\int_0^{0.5}\bigg(\frac12\ff_{17}(z)\gh_{17}(z)+2\ff_{27}(z)\gh_{27}(z)+\frac32\ff_{37}(z)\gh_{37}(z)\bigg)\,dz,\eqno(8.20)
$$
$$\aligned
b_{21}=&\frac{1}{(0.504)(0.5)\pi}\\
&\times\int_{0}^{0.5}\bigg(\frac12\ff_{17}(z)\gh_{16}(z+0.004)+2\ff_{27}(z)\gh_{26}(z+0.004)
+\frac32\ff_{37}(z)\gh_{36}(z+0.004)\bigg)\,dz.
\endaligned\eqno(8.21)
$$
$$\aligned
b_{12}=&\frac{1}{(0.504)(0.5)\pi}\\
&\times\int_{0}^{0.5}\bigg(\frac12\ff_{16}(z+0.004)\gh_{17}(z)+2\ff_{26}(z+0.004)\gh_{27}(z)
+\frac32\ff_{36}(z+0.004)\gh_{37}(z)\bigg)\,dz.\endaligned\eqno(8.22)
$$
It follows by Proposition 7.1 that
$$
\Theta_1(\bold{a}_{11},\bold{a}_{21})=\big(b_{11}+\iota_2b_{21}+\overline{\iota_2}b_{12}+|\iota_2|^2b_{22}\big)\aa\p+o(\p).
$$
This yields, by (8.7),
$$
\Xi_{11}=\mathfrak{c}_1\aa\p+o(\p)\eqno(8.23)
$$
where
$$
\mathfrak{c}_1=c_{11}+\iota_2c_{21}+\overline{\iota_2}c_{12}+|\iota_2|^2c_{22}
$$
with
$$
c_{11}=b_{11}+\overline{b_{11}},
$$
$$
c_{22}=b_{22}+\overline{b_{22}},
$$
$$
c_{12}=b_{12}+\overline{b_{21}},
$$
$$
c_{21}=\overline{c_{12}}.
$$
Let $\epsilon$ be a complex number satisfying $|\epsilon|<10^{-5}$, not necessarily the same in each occurrence.  
Numerical calculation shows that
$$
c_{11}=3.61226+\epsilon/2,
$$
$$
c_{22}=1.32215+\epsilon/2,
$$
$$
c_{12}=-0.45757-0.18179i+\epsilon/\sqrt{2}.
$$
It follows that
$$
\mathfrak{c}_1<6.9955.\eqno(8.24)
$$

 \section  *{\centerline{9. Evaluation of $\Xi_{12}$}}

  In this section we evaluate $\Xi_{12}$. By Lemma 8.1,
  $$
  \Xi_{12}=2\Re\{\Theta_1(\bold{a}_{12},\bold{a}_{22})\}+o(\p)\eqno(9.1)
  $$
 with
 with 
$$
a_{12}(n)=\chi(n)\big(\bar{\iota}_3{\vk}_3(n)+\bar{\iota}_4\vk_2(n)\big),\quad a_{22}(n)=\overline{a_{21}(n)}.\eqno(9.2)
$$
Hence
  $$\aligned
S_{j}(\bold{a}_{12},\bold{a}_{22})=&\sum_{d}\,\frac{|\chi(d)||\mu\chi(r)|\la_{0j}(dr)}{dr\varphi(r)}\bigg(\sum_m\frac{\chi(m)\big(
\bar{\iota}_3\vk_3(drm)+\bar{\iota}_4\vk_2(drm)\big)}{m^{1-\be_j}}\bigg)\\
&\times\bigg(\sum_n\frac{\chi(n)\big(\iota_3\bar{\vk}_3(drn)+\iota_4\bar{\vk}_2(drn)\big)\xi_{0j}(n;d,r)}{n}\bigg).
\endaligned
$$

By Lemma 8.2 and 8.4, for $dr<P_3/T$,
$$
\sum_m\frac{\chi(m)\vk_3(drm)}{m^{1-\be_j}}=\frac{L'(1,\chi)}{\log P_3}\f_{j6}(P_3/dr)+O(\l^{-6})
$$
and
$$
\sum_n\frac{\chi(n)\bar{\vk}_3(drn)\xi_{0j}(n;d,r)}{n}=\frac{L'(1,\chi)\Pi(d,r)}{\log P_3}\g_{j6}(P_3/dr)+O(\l^{-6}).
$$
Thus, in a way similar to the proof of (8.12), we deduce that
$$\aligned
S_{j}(\bold{a}_{12},\bold{a}_{22})=&\aa\sum_{r<D}\int_{1}^{P_3}\bigg(\bar{\iota}_3\frac{\f_{j6}(P_3/x)}{\log P_3}+\bar{\iota}_4\frac{\f_{j7}(P_2/x)}{\log P_2}\bigg)\bigg(\iota_3\frac{\g_{j6}(P_3/x)}{\log P_3}+\iota_4\frac{\g_{j7}(P_2/x)}{\log P_2}\bigg)\,\frac{dx}{x}\\
&+\frac{|\iota_4|^2\aa}{(\log P_2)^2} \int_{P_3}^{P_2}\f_{j7}(P_2/x)\g_{j7}(P_2/x)\,\frac{dx}{x}+o(\al).
\endaligned
$$
It follows, by the change of variables $x\to P_2/x$ or $x\to P_3/x$, that
$$\aligned
S_{j}(\bold{a}_{12},\bold{a}_{22})=&\frac{\aa|\iota_3|^2}{(\log P_3)^2}\int_{1}^{P_3}\frac{\f_{j6}(x)\g_{j6}(x)}{x}\,dx+\frac{\aa|\iota_4|^2}{(\log P_2)^2}\int_{1}^{P_2}\frac{\f_{j7}(x)\g_{j7}(x)}{x}\,dx\\
&+\frac{\aa\bar{\iota}_3\iota_4}{(\log P_2)(\log P_3)}\int_{1}^{P_3}\frac{\f_{j6}(x)\g_{j7}(P^{0.002}T^{-10}x)}{x}\,dx\\
&+\frac{\aa\bar{\iota}_4\iota_3}{(\log P_2)(\log P_3)}\int_{1}^{P_3}\frac{\f_{j7}(P^{0.002}T^{-10}x)\g_{j6}(x)}{x}\,dx+o(\al),
\endaligned
$$
since $P_2/P_3=P^{0.002}T^{-10}$. This yields, by substituting $x=P^z$,
$$
\frac{1}{2\al}S_1(\bold{a}_{12},\bold{a}_{22})+\frac{2}{\al}S_2(\bold{a}_{12},\bold{a}_{22})+\frac{3}{2\al}S_3(\bold{a}_{12},\bold{a}_{22})=\aa
\big(|\iota_3|^2b_{33}+\iota_3\bar{\iota}_4b_{34}+\iota_4\bar{\iota_3}b_{43}+|\iota_4|^2b_{44}\big)+o(1)
$$
where
$$
b_{33}=\frac{1}{0.498^2\pi}\int_0^{0.498}\bigg(\frac12\ff_{16}(z)\gh_{16}(z)+2\ff_{26}(z)\gh_{26}(z)+\frac32\ff_{36}(z)\gh_{36}(z)\bigg)\,dz,\eqno(9.3)
$$
$$
b_{44}=b_{22},\eqno(9.4)
$$
$$
b_{34}=\frac{1}{(0.504)(0.498)\pi}\int_{0}^{0.498}\bigg(\frac12\ff_{17}(z+0.002)\gh_{16}(z)+2\ff_{27}(z+0.002)\gh_{26}(z)
+\frac32\ff_{37}(z+0.002)\gh_{36}(z)\bigg)\,dz,\eqno(9.5)
$$
$$
b_{43}=\frac{1}{(0.504)(0.498)\pi}\int_{0}^{0.498}\bigg(\frac12\ff_{16}(z)\gh_{17}(z+0.002)+2\ff_{26}(z)\gh_{27}(z+0.002)
+\frac32\ff_{36}(z)\gh_{37}(z+0.002)\bigg)\,dz.\eqno(9.6)
$$
It follows by Proposition 7.1 that
$$
\Theta_1(\bold{a}_{12},\bold{a}_{22})=\big(|\iota_3|^2b_{33}+\iota_3\bar{\iota}_4b_{34}+\iota_4\bar{\iota}_3b_{43}+|\iota_4|^2b_{44}\big)\aa\p+o(\p)
$$
This yields
$$\Xi_{12}=\mathfrak{c}_2\aa\p+o(\p)\eqno(9.7)
$$
where
$$
\mathfrak{c}_2=|\iota_3|^2c_{33}+\iota_3\bar{\iota}_4c_{34}+\iota_4\bar{\iota}_3c_{43}+|\iota_4|^2c_{44}
$$
with
$$
c_{33}=b_{33}+\overline{b_{33}},
$$
$$
c_{44}=c_{22},
$$
$$
c_{34}=b_{34}+\overline{b_{43}},
$$
$$
c_{43}=\overline{c_{34}}.
$$
Numerical calculation shows that
$$
c_{33}=3.69507+\epsilon/2,
$$
$$
c_{34}=-0.4526+0.19474i+\epsilon/\sqrt{2}.
$$
It follows that
$$
\mathfrak{c}_2<6.9955.\eqno(9.8)
$$

 \section  *{\centerline{10. Proof of Proposition 2.4}}

The goal of this section is to prove Proposition 2.4. We continue to assume $1\le j\le 3$.

 Note that $\tilde{f}(z)\in\mathbf{R}$. If $\sigma=1/2$, then
 $$
 \overline{J_1\sp}=\sum_n\frac{\chi\bp(n)}{n^{1-s}}\tilde{f}\bigg(\frac{\log n}{\log P}\bigg)=J_1(1-s,\bp),
 $$
  $$
 \overline{J_2\sp}=\sum_n\frac{\chi\bp(n)}{n^{1-s}}\tilde{f}\bigg(\frac{\log n}{\log P}+0.004-\tilde{\al}\bigg)=J_2(1-s,\bp).
 $$
 By Lemma 8.1 we have
 $$
\sum_{\psi\in\Psi_1}\, \sum_{\rho\in\z(\psi)}\c^*\ss H_1\ss\overline{J_1\ss}\omega(\rho)=\Theta_1(\bold{a}_{11},\bold{a}_{13})+\overline{\Theta_1(\bold{a}_{13},\bold{a}_{21})}+o(\p)
 $$
 with $a_{11}(n)$ and $a_{21}(n)$ given by (8.8), and with
 $$
  a_{13}(n)=\chi(n)\tilde{f}\bigg(\frac{\log n}{\log P}\bigg),
   $$
 and
 $$
 \sum_{\psi\in\Psi_1}\,\sum_{\rho\in\z(\psi)}\c^*\ss\overline{H_2\ss}J_2\ss\omega(\rho)=\Theta_1(\bold{a}_{14},\bold{a}_{22})+\overline{\Theta_1(\bold{a}_{12},\bold{a}_{14})}+o(\p) $$
 with $a_{12}(n)$ and $a_{22}(n)$ given by (9.2), and with
  $$
  a_{14}(n)=\chi(n)\tilde{f}\bigg(\frac{\log n}{\log P}+0.004-\tilde{\al}\bigg),  $$
 Hence
$$
\Xi_1^*=\Theta_1(\bold{a}_{11},\bold{a}_{13})+\overline{\Theta_1(\bold{a}_{13},\bold{a}_{21})}
+\Theta_1(\bold{a}_{14},\bold{a}_{22})+\overline{\Theta_1(\bold{a}_{12},\bold{a}_{14})}+o(\p) .\eqno(10.1)
$$

 \medskip\noindent

{\bf Lemma 10.1.}\,\, {\it 
Write
$$
\v_{1j}(y)=\sum_m\frac{\chi(m)\tilde{f}\big(\log(ym)/\log P\big)}{m^{1-\be_j}}.
$$
If $1\le y\le P^{0.5}/T$, then
$$
\v_{1j}(y)\ll T^{-c};\eqno(10.2)
$$
if $P^{0.5}<y\le P^{0.502}/T$, then
$$
\v_{1j}(y)=\frac{500L'(1,\chi)}{\log P}\big(-1-\be_j\log(y/P^{0.5})\big)+O(\l^{-15});\eqno(10.3)
$$
if $P^{0.502}< y\le P^{0.504}/T$, then
$$
\v_{1j}(y)=\frac{500L'(1,\chi)}{\log P}\big(1-\be_j\log(P^{0.504}/y)\big)+O(\l^{-15});\eqno(10.4)
$$
if
$$
y\in(P^{0.5}/T,\,P^{0.5}]\cup(P^{0.502}/T,\,P^{0.502}]\cup(P^{0.504}/T,\,P^{0.504}),
$$
then}
$$
\v_{1j}(y)\ll \l^{-7}.\eqno(10.5)
$$

\medskip\noindent

{\it Proof.}\,\,The inequalities (10.2) follow by the Polya-Vinogradov inequality and partial summation.
Write
$$
P_1'=P^{0.504},\quad P_2'=P^{0.502},\quad P_3'=P^{0.5}.
$$
Then
$$
\tilde{f}\big(\log y/\log P\big)=\frac{500}{\log P}\cdot\tp\int_{(1)}\frac{(P_1')^s-2(P_2')^s+(P_3')^s}{y^s}\,\frac{ds}{s^2}.\eqno(10.6)
$$

In the case  $P^{0.5}< y\le P^{0.502}/T$,  
$$
\v_{1j}(y)=\frac{500}{\log P}\cdot\tp\int_{(1)}L(1-\be_j+s,\chi)
\frac{(P_1')^s-2(P_2')^s}{y^s}\,\frac{ds}{s^2}.\eqno(10.7)
$$
by (10.6). The contour of integration is moved in the same way as in the proof of lemma 8.1. Hence, by Lemma 5.8, the right side above is equal to
$$
\frac{500L'(1,\chi)}{\log P}\cdot\tp\int_{|s|=5\al}
\frac{(P_1')^s-2(P_2')^s}{y^s}\,\frac{(s-\be_j)ds}{s^2}+o(\l^{-15}).
$$
The yields (11.3) since the residue of the function
$$
\frac{(P_1')^s-2(P_2')^s}{y^s}\frac{s-\be_j}{s^2}
$$
at $s=0$ is equal to $-1-\be_j\log(y/P^{0.5})$.

Similarly, in the case  $P^{0.502}<y\le P^{0.504}/T$ we have
$$
\v_{1j}(y)=\frac{500}{\log P}\cdot\tp\int_{(1)}L(1-\be_j+s,\chi)
\bigg(\frac{P_1'}{y}\bigg)^s\,\frac{ds}{s^2}.
$$
This yields (10.4) in a way similar to the proof of (10.3). 

In the case $y\in (P^{0.5}/T,\,P^{0.5}]\cup(P^{0.504}/T,\,P^{0.504})$ and $m<T$, 
we have $$\tilde{f}\big(\log(ym)/\log P\big)\ll\al_1,$$ so (10.5)
follows by the Polya-Vinogradov inequality and partial summation; in the case $y\in (P^{0.502}/T,\,P^{0.502}]$, we use  (10.7) and a simple estimate for the integral on the new segments to obtain (10.5). $\Box$
\medskip\noindent

{\bf Lemma 10.2.}\,\, {\it Write
$$
\v_{2j}(d,r)=\sum_n\frac{\chi(n)\tilde{f}\big(\log(drn)/\log P\big)\xi_{0j}(n;d,r)}{n}.
$$
If $dr\le P^{0.5}/T$, then
$$
\v_{2j}(d,r)=\frac{L'(1,\chi)\Pi(d,r)}{500}\be_{j+1}\be_{j+2}\log P+O(\l^{-15});\eqno(10.8)
$$
 if $P^{0.5}< dr\le P^{0.502}/T$,  then
$$
\v_{2j}(d,r)=\frac{500L'(1,\chi)\Pi(d,r)}{\log P}\big(-1+\y_{1j}(dr)\big)+O(\l^{-15})\eqno(10.9)
$$
with
$$
\y_{1j}(y)=(\be_{j+1}+\be_{j+2})\log(y/P^{0.5})+\frac{\be_{j+1}\be_{j+2}}{2}\bigg(\log^2\frac{P^{0.504}}{y}-2\log^2\frac{P^{0.502}}{y}\bigg);
$$
if  $P^{0.502}<dr\le P^{0.504}/T$, then
$$
\v_{2j}(d,r)=\frac{500L'(1,\chi)\Pi(d,r)}{\log P}\big(1+\y_{2j}(dr)\big)+O(\l^{-15})\eqno(10.10)
$$
with
$$
\y_{2j}(y)=(\be_{j+1}+\be_{j+2})\log(P^{0.504}/y)+\frac{\be_{j+1}\be_{j+2}}{2}\log^2(P^{0.504}/y);
$$
if
 $$dr\in(P^{0.5}/T,\,P^{0.5}]\cup(P^{0.502}/T,\,P^{0.502}]\cup(P^{0.504}/T,\,P^{0.504}),
 $$
 then}
$$
\v_{2j}(d,r)\ll \l^{-7}.\eqno(10.11)
$$

\medskip\noindent

{\it Proof.}\,\,First assume $dr\le P^{0.5}/T$. In view of (8.9) and (10.6),  we have
$$\aligned
\v_{2j}(d,r)=\frac{500}{\log P}&\cdot\int_{(1)}\frac{L(1+\be_{j+1}+s,\chi)L(1+\be_{j+2}+s,\chi)\u_{0j}(s;d,r)
}{L(1+s,\chi)}\\
&\times\frac{(P_1')^s-2(P_2')^s+(P_3')^s}{(dr)^s}\,\frac{ds}{s^2}.
\endaligned
$$
The contour of integration is moved in the same way as the proof of Lemma 8.1. Thus, by Lemma 5.8 and 8.2,
$$\aligned
\v_{2j}(d,r)=\frac{500L'(1,\chi)\Pi(d,r)}{\log P}\cdot\int_{|s|=10\al}\frac{(\be_{j+1}+s)(\be_{j+2}+s)}{s}\frac{(P_1')^s-2(P_2')^s+(P_3')^s}{(dr)^s}\,\frac{ds}{s^2}+O(\l^{-15}).
\endaligned
$$
This yields (10.8) since the function
$$
\frac{(P_1')^s-2(P_2')^s+(P_3')^s}{s^2}
$$
is analytic and equal to $(\log P/500)^2$ at $s=0$.

Similarly, in the case $P^{0.5}/T<dr\le P^{0.502}$ we have
$$\aligned
\v_{2j}(d,r)=\frac{500L'(1,\chi)\Pi(d,r)}{\log P}\cdot\tp\int_{|s|=10\al}\frac{(\be_{j+1}+s)(\be_{j+2}+s)}{s}\frac{(P_1')^s-2(P_2')^s}{(dr)^s}\,\frac{ds}{s^2}+O(\l^{-15});
\endaligned
$$
 in the case  $P^{0.502}<dr\le P^{0.504}/T$ we have
$$\aligned
\v_{2j}(d,r)=\frac{500L'(1,\chi)\Pi(d,r)}{\log P}\cdot\tp\int_{|s|=10\al}\frac{(\be_{j+1}+s)(\be_{j+2}+s)}{s}\bigg(\frac{P_1'}{dr}\bigg)^s\,\frac{ds}{s^2}+O(\l^{-15}).
\endaligned
$$
These yield (10.9) and (10.10). The proof of (10.11) is similar to that of (10.5) in the case $y\in (P^{0.502}/T,\,P^{0.502}]$. $\Box$

\medskip\noindent

By direct calculation, for $0\le z\le 1$ we have
$$
\be_j\log P^{z-0.5}=\pi ij(z-0.5)+O(\l^{-8}),\quad\be_j\log P^{0.504-z}=\pi ij(0.504-z)+O(\l^{-8}),
$$
$$
\y_{1j}(P^z)=\yy_{1j}(z)+O(\l^{-8}),\quad \y_{2j}(P^z)=\yy_{2j}(z)+O(\l^{-8}),
$$
where
$$
\yy_{11}(z)=5\pi i(z-0.5)-3\pi^2\big((0.504-z)^2-2(0.502-z)^2\big),
$$
$$ \yy_{21}(z)=5\pi i(0.504-z)-3\pi^2(0.504-z)^2,
$$
$$
\yy_{12}(z)=4\pi i(z-0.5)-\frac{3\pi^2}{2}\big((0.504-z)^2-2(0.502-z)^2\big),
$$
$$
 \yy_{22}(z)=4\pi i(0.504-z)-\frac{3\pi^2}{2}(0.504-z)^2,
$$
$$
\yy_{13}(z)=3\pi i(z-0.5)-\pi^2\big((0.504-z)^2-2(0.502-z)^2\big),
$$
$$
 \yy_{23}(z)=3\pi i(0.504-z)-\pi^2\big((0.504-z)^2.
$$

{\it Remark.}\,\,Similar to (10.6),
$$
\tilde{f}\big(\log y/\log P+0.004-\tilde{\al}\big)=\frac{500}{\log P}\cdot\tp\int_{(1)}\bigg(\frac{Dt_0}{P^{0.004}}\bigg)^s\frac{(P_1')^s-2(P_2')^s+(P_3')^s}{y^s}\,\frac{ds}{s^2}.
$$
Thus, with simple modification, Lemma 10.1 and 10.2 apply to the sums 
$$
\sum_m\frac{\chi(m)\tilde{f}\big(\log(ym)/\log P+0.004-\tilde{\al}\big)}{m^{1-\be_j}}
$$
and
$$
\sum_n\frac{\chi(n)\tilde{f}\big(\log(drn)/\log P+0.004-\tilde{\al}\big)\xi_{0j}(n;d,r)}{n}.
$$

\medskip\noindent

{\it Evaluation of $\Theta_1(\bold{a}_{11},\bold{a}_{13})$.}
\medskip\noindent

We have

$$\aligned
S_{j}(\bold{a}_{11},\bold{a}_{13})=&\sum_{r}\,\sum_d\frac{|\chi(d)||\mu\chi(r)|\la_{0j}(dr)}{dr\varphi(r)}\bigg(\sum_m\frac{\chi(m)\big(
\vk_1(drm)+\iota_2\vk_2(drm)\big)}{m^{1-\be_j}}\bigg)\\
&\times\bigg(\sum_n\frac{\chi(n)\tilde{f}\big(\log(drn)/\log P\big)\xi_{0j}(n;d,r)}{n}\bigg).
\endaligned
$$
The right side is split into three sums according to
 $$dr<P^{0.5},\quad P^{0.5}\le dr<P^{0.502},\quad
P^{0.502}\le dr<P^{0.504}.$$ By Lemma 10.2 and the results in Section 8,  the sum over $dr<P^{0.5}$ is equal to
$$
\aligned
\frac{L'(1,\chi)^2}{500}&\be_{j+1}\be_{j+2}\sum_{n<P^{0.5}}\,\frac{|\chi(n)|\la_{0j}(n)}{\varphi(n)}
\bigg(\frac{\f_{j6}(P^{0.504}/n)}{0.504}+\iota_2\frac{\f_{j7}(P^{0.5}/n)}{0.5}\bigg)+o(\al)
\\=&\frac{\aa\be_{j+1}\be_{j+2}}{500}\int_{1}^{P^{0.5}}\bigg(\frac{\f_{j6}(P^{0.504}/x)}{0.504}+\iota_2\frac{\f_{j7}(P^{0.5}/x)}{0.5}\bigg)\,\frac{dx}{x}+o(\al)\\
=&\frac{\aa\be_{j+1}\be_{j+2}\log P}{500}\int_{0}^{0.5}\bigg(\frac{\ff_{j6}(0.004+z)}{0.504}+\iota_2\frac{\ff_{j7}(z)}{0.5}\bigg)\,dz+o(\al).
\endaligned
$$
 Similarly, the sum over $P^{0.5}\le dr<P^{0.502}$ is equal to
$$\aligned
\frac{500L'(1,\chi)^2}{0.504\log^2P}&
\sum_{P^{0.5}\le n<P^{0.502}}\,\frac{|\chi(n)|\la_{0j}(n)}{\varphi(n)}
\f_{j6}(P^{0.504}/n)\big(-1+\y_{1j}(n)\big)+o(\al)\\
=&\frac{500\aa}{0.504\log P}
\int_{0.5}^{0.502}
\ff_{j6}(0.504-z)\big(-1+\yy_{1j}(z)\big)\,dz+o(\al);
\endaligned
$$
the sum over $P^{0.502}\le dr<P^{0.504}$ is equal to
$$\aligned
\frac{500L'(1,\chi)^2}{0.504\log^2P}&
\sum_{P^{0.502}\le n<P^{0.504}}\frac{|\chi(n)|\la_{0j}(n)}{\varphi(n)}
\f_{j6}(P^{0.504}/n)\big(1+\y_{2j}(n)\big)+o(\al)\\
=&\frac{500\aa}{0.504\log P}
\int_{0.502}^{0.504}
\ff_{j6}(0.504-z)\big(1+\yy_{2j}(z)\big)\,dz+o(\al).
\endaligned
$$
Noting that
$$
\be_{j+1}\be_{j+2}\log P=-(11-6j+j^2)\pi\al+o(\al)
 $$
and  gathering the above results together we conclude
 $$
 \frac{1}{\al}S_{j}(\bold{a}_{11},\bold{a}_{13})=d_{3j}\aa+o(1)
 $$
 with
 $$\aligned
d_{3j} =&-\frac{(11-6j+j^2)\pi}{500}\int_0^{0.5}\bigg(\frac{1}{0.504}\ff_{j6}(0.004+z)+\frac{\iota_2}{0.5}\ff_{j7}(z)\bigg)\,dz\\
 &+\frac{500}{0.504\pi}\int_0^{0.002}\big(\ff_{j6}(z)-\ff_{j6}(0.002+z)\big)\,dz\\
 &+\frac{500}{0.504\pi}\bigg(\int_{0.5}^{0.502}\ff_{j6}(0.504-z)\yy_{1j}(z)\,dz+\int_{0.502}^{0.504}\ff_{j6}(0.504-z)\yy_{2j}(z)\,dz\bigg).
 \endaligned
 $$
Hence, by Proposition 7.1,
 $$
 \Theta_1(\bold{a}_{11},\bold{a}_{13})=\bigg(\frac12d_{31}+2d_{32}+\frac32d_{33}\bigg)\aa\p+o(\p).\eqno(10.12)
 $$
\medskip\noindent

{\it Evaluation of $\Theta_1(\bold{a}_{13},\bold{a}_{21})$.}
\medskip\noindent
\medskip\noindent

We have
$$\aligned
S_{j}(\bold{a}_{13},\bold{a}_{21})=&\sum_{r}\,\sum_d\frac{|\chi(d)||\mu\chi(r)|\la_{0j}(dr)}{dr\varphi(r)}\bigg(\sum_m\frac{\chi(m)\tilde{f}\big(\log(drm)/\log P\big)}{m^{1-\be_j}}\bigg)\\
&\times\bigg(\sum_n\frac{\chi(n)\big(
\bar{\vk}_1(drn)+\bar{\iota}_2\bar{\vk}_2(drn)\big)\xi_{0j}(n;d,r)}{n}\bigg).
\endaligned
$$
The right side is split into three sums in the same way as in the last subsection.
  By Lemma 10.1 and the results in Section 8, the sum over $dr<P^{0.5}$ is $o(\al)$;
   the sum over $P^{0.5}\le dr<P^{0.502}$ is equal to
  $$\aligned
&\frac{500L'(1,\chi)^2}{0.504\log^2P}
\sum_{P^{0.5}\le n<P^{0.502}}\frac{|\chi(n)|\la_{0j}(n)}{\varphi(n)}
\big(-1-\be_j\log(y/P^{0.5})\big)\g_{j6}(P^{0.504}/dr)+o(\al)\\
&=\frac{500\aa}{0.504\log P}
\int_{0.5}^{0.502}
\big(-1-\pi ij(z-0.5)\big)\gh_{j6}(0.504-z)\,dz+o(\al);
\endaligned
$$
 the sum over $P^{0.502}<dr\le P^{0.504}$ is equal to
$$\aligned
&\frac{500L'(1,\chi)^2}{0.504\log^2P}
\sum_{P^{0.502}\le n<P^{0.504}}\frac{|\chi(n)|\la_{0j}(n)}{\varphi(n)}
\big(1-\be_j\log(P^{0.504}/n)\big)\g_{j6}(P^{0.504}/n)+o(\al)\\
&=\frac{500\aa}{0.504\log P}
\int_{0.502}^{0.504}
\big(1-\pi ij(0.504-z)\big)\gh_{j6}(0.504-z)\,dz+o(\al),
\endaligned
$$
Gathering these results together we conclude
  $$
 \frac{1}{\al}S_{j}(\bold{a}_{13},\bold{a}_{21})=d_{4j}\aa+o(1)
 $$
 with
 $$\aligned
 d_{4j}=
 &\frac{500}{0.504\pi}\int_0^{0.002}\big(\gh_{j6}(z)-\gh_{j6}(0.002+z)\big)\,dz\\
 -&\frac{500ij}{0.504}\int_{0}^{0.002}\big(\gh_{j6}(0.002+z)(0.002-z)+\gh_{j6}(z)z\big)\,dz.\endaligned
 $$
Hence, by Proposition 7.1,
 $$
 \Theta_1(\bold{a}_{13},\bold{a}_{21})=\bigg(\frac12d_{41}+2d_{42}+\frac32d_{43}\bigg)\aa\p+o(\p).\eqno(10.13)
 $$
 \medskip\noindent
 
{\it Evaluation of $\Theta_1(\bold{a}_{14},\bold{a}_{22})$.}
\medskip\noindent

We have
$$\aligned
S_{j}(\bold{a}_{14},\bold{a}_{22})=&\sum_{r}\,\sum_d\frac{|\chi(d)||\mu\chi(r)|\la_{0j}(dr)}{dr\varphi(r)}\bigg(\sum_m\frac{\chi(m)\tilde{f}\big(\log(drm)/\log P+0.004-\tilde{\al}\big)}{m^{1-\be_j}}\bigg)\\
&\times\bigg(\sum_n\frac{\chi(n)\big(\iota_3
\bar{\vk}_3(drn)+\iota_4\bar{\vk}_2(drn)\big)\xi_{0j}(n;d,r)}{n}\bigg).
\endaligned
$$
The right side is split into three sums according to
 $$dr<P^{0.496},\quad P^{0.496}\le dr<P^{0.498},\quad
P^{0.498}\le dr<P^{0.5}. $$
 By a result  similar to Lemma 10.1 and the results in Section 8, the sum over $dr<P^{0.496}$ is $o(\al)$; the sum over $P^{0.496}\le dr<P^{0.498}$ is equal to
 $$\aligned
&\frac{500L'(1,\chi)^2}{\log^2P}
\sum_{P^{0.496}\le n<P^{0.498}}\frac{|\chi(n)|\la_{0j}(n)}{\varphi(n)}
\big(-1-\be_j\log(P^{-0.496}n)\big)\\
&\,\,\times\bigg(\frac{\iota_3}{0.498}\g_{j6}(P^{0.498}/n)+\frac{\iota_4}{0.5}\g_{j7}(P^{0.5}/n)\bigg)+o(\al)\\
&=\frac{500\aa}{\log P}
\int_{0.496}^{0.498}
\big(-1-\pi ij(z-0.496)\big)\\
&\times\bigg(\frac{\iota_3}{0.498}\gh_{j6}(0.498-z)+\frac{\iota_4}{0.5}\gh_{j7}(0.5-z)\big)\,dz+o(\al);
\endaligned
$$
the sum over $P^{0.498}\le dr<P^{0.5}$ is equal to
$$\aligned
&\frac{500L'(1,\chi)^2}{\log^2P}\frac{\iota_4}{0.5}
\sum_{P^{0.498}\le n<P^{0.5}}\frac{|\chi(n)|\la_{0j}(n)}{\varphi(n)}
\big(1-\be_j\log(P^{0.5}/n)\big)\g_{j7}(P^{0.5}/n)+o(\al)\\
=&\frac{1000\iota_4\aa}{\log P}
\int_{0.498}^{0.5}
\big(1-\pi ij(0.5-z)\big)\gh_{j7}(0.5-z)\,dz+o(\al).
\endaligned
$$
Gathering the above results together  we conclude
  $$
 \frac{1}{\al}S_{j}(\bold{a}_{14},\bold{a}_{22})=(d_{5j}'+d_{5j})\aa+o(1) $$
 with
 $$
 d_{5j}'=-\frac{500\iota_3}{0.498\pi}\int_0^{0.002}\gh_{j6}(z)\,dz,
 $$
 $$\aligned
 d_{5j}=  &\frac{1000\iota_4}{\pi}\int_{0}^{0.002}\big(\gh_{j7}(z)-\gh_{j7}(0.002+z)\big) \,dz\\
 -&\frac{500ij\iota_3}{0.498}\int_{0}^{0.002}(0.002-z)\gh_{j6}(z) \,dz\\
 -&1000ij\iota_4\int_{0}^{0.002}\big((0.002-z)\gh_{j7}(0.002+z)+z\gh_{j7}(z)\big)\,dz.
 \endaligned
 $$
Hence, by Proposition 7.1,
 $$
\Theta_1(\bold{a}_{14},\bold{a}_{22})=\bigg(\frac12(d_{51}'+d_{51})+2(d_{52}'+d_{52})+\frac32(d_{53}'+d_{53})\bigg)\aa\p+o(\p).
 \eqno(10.14)$$

  \medskip\noindent
  
{\it Evaluation of $\Theta_1(\bold{a}_{12},\bold{a}_{14})$.}
\medskip\noindent

We have
$$\aligned
S_{j}(\bold{a}_{12},\bold{a}_{14})=&\sum_{r}\,\sum_d\frac{|\chi(d)||\mu\chi(r)|\la_{0j}(dr)}{dr\varphi(r)}\bigg(\sum_m\frac{\chi(m)\big(\bar{\iota}_3
\vk_3(drm)+\bar{\iota}_4\vk_2(drm)\big)}{m^{1-\be_j}}\bigg)\\
&\times\bigg(\sum_n\frac{\chi(n)\tilde{f}\big(\log(drn)/\log P+0.004-\tilde{\al}\big)\xi_{0j}(n;d,r)}{n}\bigg).
\endaligned
$$
The right side is split into three sums in the same way as  in the last subsection.
By a result similar to Lemma 10.2 and the results in Section 8, the sum over $dr<P^{0.496}$ is equal to
$$\aligned
\frac{L'(1,\chi)^2}{500}&\be_{j+1}\be_{j+2}\sum_{n<P^{0.496}}\frac{|\chi(n)|\la_{0j}(n)}{\varphi(n)}
\bigg(\frac{\bar{\iota}_3\f_{j6}(P^{0.498}/n)}{0.498}+\frac{\bar{\iota}_4\f_{j7}(P^{0.5}/n)}{0.5}\bigg)+o(\al)\\
=&\frac{\aa\be_{j+1}\be_{j+2}}{500}\int_{1}^{P^{0.496}}\bigg(\frac{\bar{\iota}_3\f_{j6}(P^{0.498}/x)}{0.498}+\frac{\bar{\iota}_4\f_{j7}(P^{0.5}/x)}{0.5}\bigg)\,\frac{dx}{x}+o(\al)\\
=&\frac{\aa\be_{j+1}\be_{j+2}\log P}{500}\int_{0}^{0.496}\bigg(\frac{\bar{\iota}_3\ff_{j6}(0.002+z)}{0.498}+\frac{\bar{\iota}_4\ff_{j7}(0.004+z)}{0.5}\bigg)\,dz+o(\al);
\endaligned
$$
the sum over $P^{0.496}\le dr<P^{0.498}$ is equal to
$$\aligned
&\frac{500L'(1,\chi)^2}{\log^2P}
\sum_{p^{0.496}\le n<P^{0.498}}\frac{|\chi(n)|\la_{0j}(n)}{\varphi(n)}
\bigg(\frac{\bar{\iota}_3\f_{j6}(P^{0.498}/n)}{0.498}+\frac{\bar{\iota}_4\f_{j7}(P^{0.5}/n)}{0.5}\bigg)\\
&\qquad\qquad\times\big(-1+\y_{1j}(P^{0.004}n)\big)+o(\al)\\
&=\frac{500\aa}{\log P}
\int_{0}^{0.002}
\bigg(\frac{\bar{\iota}_3\ff_{j6}(z)}{0.498}+\frac{\bar{\iota}_4\ff_{j7}(0.002+z)}{0.5}\bigg)
\big(-1+\yy_{1j}(0.502-z)\big)\,dz+o(\al);
\endaligned
$$
the sum over $P^{0.498}\le dr<P^{0.5}$ is equal to
$$\aligned
&\frac{1000\bar{\iota_4}L'(1,\chi)^2}{\log^2P}\
\sum_{p^{0.498}\le n<P^{0.5}}\frac{|\chi(n)|\la_{0j}(n)}{\varphi(n)}
\f_{j7}(P^{0.5}/n)\big(1+\y_{2j}(P^{0.004}n)\big)+o(\al)\\
&=\frac{1000\bar{\iota_4}\aa}{\log P}
\int_{0}^{0.002}
\ff_{j7}(z)\big(1+\yy_{2j}(0.504-z)\big)\,dz+o(\al).
\endaligned
$$
Gathering the above results together  we conclude
  $$
 \frac{1}{\al}S_{j}(\bold{a}_{12},\bold{a}_{14})=(d_{6j}'+d_{6j})\aa+o(1) $$
 with
 $$
 d_{6j}'=-\frac{500\bar{\iota}_3}{0.498\pi}\int_0^{0.002}\ff_{j6}(z)\,dz,\eqno(10.15)
 $$
 $$\aligned
 d_{6j}=&-\frac{(11-6j+j^2)\pi}{500}\int_{0}^{0.496}\bigg(\frac{\bar{\iota}_3\ff_{j6}(0.002+z)}{0.498}+\frac{\bar{\iota}_4\ff_{j7}(0.004+z)}{0.5}\bigg)\,dz\\
&+ \frac{1000\bar{\iota}_4}{\pi}\int_{0}^{0.002}\big(\ff_{j7}(z)-\ff_{j7}(0.002+z)\big)\,dz\\
+&\frac{500}{\pi}
\int_{0}^{0.002}
\bigg(\frac{\bar{\iota}_3\ff_{j6}(z)}{0.498}+\frac{\bar{\iota}_4\ff_{j7}(0.002+z)}{0.5}\bigg)
\yy_{1j}(0.502-z)\,dz\\
+&\frac{1000\bar{\iota}_4}{\pi}
\int_{0}^{0.002}
\ff_{j7}(z)\yy_{2j}(0.504-z)\,dz.
\endaligned
$$
Hence, by Proposition 7.1,
 $$
 \Theta_1(\bold{a}_{12},\bold{a}_{14})=\bigg(\frac12(d_{61}'+d_{61})+2(d_{62}'+d_{62})+\frac32d_{63}'+d_{63})\bigg)\aa\p+o(\p).
 \eqno(10.16)$$

 It follows from (10.1) and (10.12)-(10.16) that
$$
\Xi_1^*=(\mathfrak{d}'+\mathfrak{d})\aa\p+o(\p)\eqno(10.17)
$$
where
$$
\mathfrak{d}'=\frac12\big(d_{51}'+\overline{d'_{61}}\big)+2\big(d_{52}'+\overline{d'_{62}}\big)+\frac32\big(d_{53}'+\overline{d'_{63}}\big),
$$
$$
\mathfrak{d}=\frac12\big(d_{31}+d_{51}+\overline{d_{41}+d_{61}}\big)+2\big(d_{32}+d_{52}+\overline{d_{42}+d_{62}}\big)+\frac32\big(d_{33}+d_{53}+\overline{d_{43}+d_{63}}\big)
$$
For $\mu=6, 7$,
$$\ff_{j\mu}(0)=\gh_{j\mu}(0)=1.
$$
Hence
$$
\mathfrak{d}'\simeq -\frac{8\iota_3}{0.498\pi}.
$$
We only need a crude lower bound for the real part of $\mathfrak{d}'$. By direct calculation we have 
$$
\Re\{\mathfrak{d}'\}>-\frac{8}{0.498\pi}\Re\{\iota_3\}-0.04>5.1.
$$
The contribution from $\mathfrak{d}$ is minor since if $z$ is close to $1/2$, then for $\mu=1, 2$,
$$
\yy_{\mu j}(z)\simeq 0.
$$
Thus, by direct calculation we have the crude bound
$$
|\Re\{\mathfrak{d}\}|<0.1.
$$
Combining these bounds with (10.16) we complete the proof of Proposition 2.4. $\Box$

\medskip\noindent

{\it Remark.} Numerical calculation actually shows that
$$
\Re\{\mathfrak{d}\}>0.
$$

\medskip\noindent

 \section  *{\centerline{11.  Proof of Proposition 2.6}}
\medskip\noindent 

By the result of Section 9,
$$
\Xi_{12}\ll\aa\p.
$$
Hence, by Cauchy's inequality, the proof of Proposition 2.6 is reduced to showing that
$$
\sum_{\psi\in\Psi_1}\,\sum_{\rho\in\tz(\psi)}\c^*\ss|J_1\ss-Z(\rho,\pc)J_2(1-\rho,\bp)|^2\omega(\rho)=o(\aa\p).\eqno(11.1)
$$

Let $\tg_1(y)$ and $\tg_2(y)$ be given by
$$
\tg_1(y)=-500\int_{0.5}^{0.502}g\bigg(\frac{P^z}{y}\bigg)\,dz+500\int_{0.502}^{0.504}g\bigg(\frac{P^z}{y}\bigg)\,dz,
$$
$$
\tg_2(y)=-500\int_{0.496}^{0.498}g\bigg(\frac{P^zDt_0}{y}\bigg)\,dz+500\int_{0.498}^{0.5}g\bigg(\frac{P^zDt_0}{y}\bigg)\,dz.
$$
Write
$$
\tilde{J}_\mu\sp=\sum_n\frac{\pc(n)\tg_\mu(n)}{n^s},\quad \mu=1,2,
$$
and
$$
\eta_\pm=\exp\{\pm\l^{-10}\}. 
$$
By (5) and (5), for $\sigma=1/2$, the terms with $n\le P^{0.5}\eta_-$ or $n\ge P^{0.504}\eta_+$ in $\tilde{J}_1\sp$ contribute
 $\ll\ep$.
 \medskip\noindent

 {\bf Lemma 11.1.}\,\,{\it If
 $$
 y\in[P^{0.5}\eta_+,\,P^{0.502}\eta_-]\cup[P^{0.502}\eta_+,\,P^{0.504}\eta_-],
 $$
 then
 $$
 \tilde{f}\bigg(\frac{\log y}{\log P}\bigg)-\tg_1(y)\ll\ep;\eqno(11.2)
 $$
 if
  $$
 y\in(P^{0.5}\eta_-,\,,P^{0.5}\eta_+)\cup(P^{0.502}\eta_-,\,P^{0.502}\eta_+)\cup(P^{0.504}\eta_-,\,P^{0.504}\eta_+),
 $$
 then}
 $$
 \tilde{f}\bigg(\frac{\log y}{\log P}\bigg)-\tg_1(y)\ll\l^{-10}.\eqno(11.3)
 $$
\medskip\noindent

{\it Proof.}\,\, 
Since
$$
g(1/x)=\frac{1}{\sqrt{\pi}}\int_{-\infty}^{-\l\log x}\exp(-t^2)\,dt=\frac{1}{\sqrt{\pi}}\int_{\l\log x}^{\infty}\exp(-t^2)\,dt,
$$
it follows that
$$
g(x)+g(1/x)=1.\eqno(11.4)
$$
Write
$$
u=\frac{\log y}{\log P}.
$$

First assume  $T^{0.502}\eta_+\le y\le P^{0.503}$. By (11.4),
$$
\int_{0.502}^{2u-0.502}g\bigg(\frac{P^z}{y}\bigg)\,dz=u-0.502;
$$
by (6.2),
$$
\int_{2u-0.502}^{0.504}g\bigg(\frac{P^z}{y}\bigg)\,dz=1.006-2u+O(\ep);
$$
by (6.3),
$$
\int_{0.5}^{0.502}g\bigg(\frac{P^z}{y}\bigg)\,dz\ll\ep.
$$
These together imply (11.2).

Now assume  $P^{0.503}\le y\le P^{0.504}\eta_-$. By (11.4),
$$
\int_{2u-0.504}^{0.504}g\bigg(\frac{P^z}{y}\bigg)\,dz=0.504-u;
$$
by (6.2),
$$
\int_{0.502}^{2u-0.504}g\bigg(\frac{P^z}{y}\bigg)\,dz=2u-1.006+O(\ep);
$$
by (6.3),
$$
\int_{0.5}^{0.502}g\bigg(\frac{P^z}{y}\bigg)\,dz\ll\ep.
$$
These together imply (11.2). The proof of (11.2) with $ y\in[P^{0.5}\eta_+,\,P^{0.502}\eta_-]$ is similar.

We give the  proof of   (11.3) with  $P^{0.502}\eta_-<y<P^{0.502}\eta_+$ only; the other cases are similar.  We have
$$
\int_{0.502}^{0.504}g\bigg(\frac{P^z}{y}\bigg)\,dz=\frac{1}{500}+O(\l^{-10}),\quad
\int_{0.5}^{0.502}g\bigg(\frac{P^z}{y}\bigg)\,dz=O(\l^{-10}),
$$
by  (5.2) and (5.3), and 
$$
\tilde{f}\bigg(\frac{\log y}{\log P}\bigg)=1+O(\l^{-10}).\quad\Box
$$

It follows from Lemma 11.1 that
$$
J_1\sp-\tilde{J}_1\sp=\sum_{n\in\mathfrak{I}_1}\frac{\pc(n)}{n^s}\big(\tilde{f}(\log n/\log P)-\tilde{g}_1(n)\big)+O(\ep)\eqno(11.5)
$$
for $\sigma=1/2$, where
$$
\mathfrak{I}_1=(P^{0.5}\eta_-,\,P^{0.5}\eta_+)\cup (P^{0.502}\eta_-,\,P^{0.502}\eta_+)\cup(P^{0.504}\eta_-,\,P^{0.504}\eta_+).
$$
By (11.3), (8.25) and (8.26),
$$
\sum_{\psi\in\Psi_1}\,\sum_{\rho\in\tz(\psi)}\c^*\ss\bigg|\sum_{n\in\mathfrak{I}_1}\frac{\pc(n)}{n^\rho}\big(\tilde{f}(\log n/\log P)-\tilde{g}_1(n)\big)\bigg|^2\omega(\rho)=o(\aa\p).
$$
Hence
$$
\sum_{\psi\in\Psi_1}\,\sum_{\rho\in\tz(\psi)}\c^*\ss\big|J_1\ss-\tilde{J}_1\ss\big|^2\omega(\rho)=o(\aa\p).
$$
Similarly,
$$
\sum_{\psi\in\Psi_1}\,\sum_{\rho\in\tz(\psi)}\c^*\ss\big|J_2\ss-\tilde{J}_2\ss\big|^2\omega(\rho)=o(\aa\p).
$$
The proof of (11.1) is therefore reduced to showing that
$$
\sum_{\psi\in\Psi_1}\,\sum_{\rho\in\tz(\psi)}\c^*\ss\big|\tilde{J}_1\ss-Z(\rho,\pc)\tilde{J}_2(1-\rho,\bp)\big|^2\omega(\rho)=o(\aa\p).
\eqno(11.6)$$

\medskip\noindent

{\bf Lemma 11.2.}\,\, {\it Suppose $\sigma=1/2$ and $|t-2\pi t_0|<\l_1$. Then
$$
\tilde{J}_1\sp=Z(s,\pc)\tilde{J}_2(1-s,\bp)+O(E_2\sp)
$$
where}
$$
E\sp=\l^{-68}\int_{-\l^{20}}^{\l^{20}}\bigg|\sum_{n<P_1}\frac{\pc(n)}{n^{s+iv}}\bigg|\omega_1(iv)\,dv.
$$
\medskip\noindent

{\it Proof.}\,\, Suppose $0.5\le z\le 0.504$. In a way similar to the proof of Lemma 6.1, it can be deduced that 
$$
\sum_n\frac{\pc(n)}{n^s}g\bigg(\frac{P^z}{n}\bigg)=L(s,\pc)-Z(s,\pc)\sum_n\frac{\chi\bp(n)}{n^{1-s}}g\bigg(\frac{P^{1-z}Dt_0}{n}\bigg)+O(E_2\sp).
$$
Hence
$$\aligned
\int_{0.5}^{0.502}&\bigg\{\sum_n\frac{\pc(n)}{n^s}g\bigg(\frac{P^z}{n}\bigg)\bigg\}\,dz\\&=\frac{1}{500}L(s,\pc)
-Z(s,\pc)\int_{0.498}^{0.5}\bigg\{\sum_n\frac{\chi\bp(n)}{n^{1-s}}g\bigg(\frac{P^zDt_0}{n}\bigg)\bigg\}\,dz+O(E_2\sp),
\endaligned
$$
and
$$\aligned
\int_{0.502}^{0.504}&\bigg\{\sum_n\frac{\pc(n)}{n^s}g\bigg(\frac{P^z}{n}\bigg)\bigg\}\,dz\\
&=\frac{1}{500}L(s,\pc)-Z(s,\pc)\int_{0.496}^{0.498}\bigg\{\sum_n\frac{\chi\bp(n)}{n^{1-s}}g\bigg(\frac{P^zDt_0}{n}\bigg)\bigg\}\,dz+O(E_2\sp).
\endaligned
$$
This completes the proof. $\Box$
\medskip\noindent

By Lemma 11.2, the proof of (11.6) is reduced to showing that
$$
\sum_{\psi\in\Psi_1}\,\sum_{\rho\in\tz(\psi)}\c^*\ss E_2\ss^2\omega(\rho)=o(\aa\p).
$$
This follows by (8.25), (8.26) and simple estimates. 

 \section  *{\centerline{12.  Evaluation of $\Xi_{15}$}}
 \medskip\noindent

The sum $H_{11}\sp$ is split into
$$
H_{11}\sp=\sum_{n<P^{1/2}}+\sum_{P^{1/2}\le n<P_1}=H_{14}\sp+H_{15}\sp,\quad\text{say}.\eqno(12.1)
$$
We can write
$$
H_1\sp H_2\sp=B\sp+H_{15}\sp H_2\sp
$$
with
$$
B\sp=(H_{14}\sp+\iota_2H_{12}\sp)H_2\sp.\eqno(12.2)
$$
Accordingly we have
$$
\Xi_{13}=\Xi_{14}+\Xi_{15}\eqno(12.3)
$$
with
$$
\Xi_{14}=\sum_{\psi\in\Psi_1}\,\sum_{\rho\in\tz(\psi)}\c^*\ss Z(\rho,\pc)^{-1}B\ss\omega(\rho),\eqno(12.4)
$$
$$
\Xi_{15}=\sum_{\psi\in\Psi_1}\,\sum_{\rho\in\tz(\psi)}\c^*\ss Z(\rho,\pc)^{-1}H_{15}\ss H_2\ss\omega(\rho).\eqno(12.5)
$$
The goal of this section is to evaluate $\Xi_{15}$.

Write
$$
\vk_{12}(y)=\frac{1}{0.504}\bigg(\frac{P_1}{y}\bigg)^{\be_6}\int_{0.5}^{0.504}\bigg\{g\bigg(\frac{P^z}{y}\bigg)-g\bigg(\frac{P^{0.5}}{y}\bigg)\bigg\}\,dz
$$
and
$$
\tilde{H}_{15}\sp=\sum_{P^{0.5}\eta_-<n<P_1\eta_+}\frac{\vk_{12}(n)\pc(n)}{n^s}.
$$
We first claim that
$$
\sum_{\psi\in\Psi_1}\,\sum_{\rho\in\tz(\psi)}\c^*\ss |H_{15}\ss-\tilde{H}_{15}\ss|^2\omega(\rho)=o(\aa\p).\eqno(12.6)
$$

By the argument in the last section, for $P^{0.5}\eta_+<n<P_1\eta_-$,
$$
\int_{0.5}^{0.504}g\bigg(\frac{P^z}{y}\bigg)\,dz=\frac{\log P_1-\log n}{\log P}+O(\ep),
$$
so that
$$
\vk_{12}(n)-\vk_1(n)\ll\ep.
$$
Also we have
$$
\vk_{12}(n)-\vk_1(n)\ll
\begin{cases}
\l^{-10}&\quad\text{if}\quad P_1\eta_-\le n< P_1\eta_+,\\
1&\quad\text{if}\quad P^{0.5}\eta_-< n\le P^{0.5}\eta_+.
\end{cases}
$$
 These bounds  together with (8.25) and (8.26) imply  (12.6). We briefly describe the argument as follows. Let $\mathfrak{t}(y)$
denote the characteristic function of the interval $(P^{0.5}\eta_-, P^{0.5}\eta_+]$. If $d\le P^{0.5}(\eta_+-\eta_-)$, then
$$
\sum_l\frac{\mathfrak{t}(dl)}{l}\ll\l^{-10}.
$$
If $d>P^{0.5}(\eta_+-\eta_-)$, then there exists at most one $l$ such that $\mathfrak{t}(dl)=1$. Hence
$$
\sum_{d>P^{0.5}(\eta_+-\eta_-)}\frac{1}{d}\bigg(\sum_l\frac{\mathfrak{t}(dl)}{l}\bigg)\le \sum_l\frac{1}{l^2}\sum_d\frac
{\mathfrak{t}(dl)}{d}\ll\l^{-10}.
$$

Assume $\sigma=1/2$, $|t-2\pi t_0|<\l_1$ and $0.5\le z\le 0.504$. Using the proof of Lemma 11.2 with $s+\be_6$ in place of $s$ we deduce that
$$\aligned
\tilde{H}_{15}\sp=&\frac{Z(s+\be_6,\pc)P_1^{\be_6}}{0.504}\sum _n\frac{\chi\bp(n)}{n^{1-s-\be_6}}\int_{0.496}^{0.5}\bigg\{g\bigg(\frac{P^{0.5}Dt_0}{n}\bigg)-g\bigg(\frac{P^{z}Dt_0}{n}\bigg)\bigg\}\,dz\\&+O(E(s+\be_6,\psi)).
\endaligned$$
Write
$$
P''_1=P^{0.496}Dt_0,\qquad P''_2=P^{0.5}Dt_0.
$$
By (4) and (4), the terms with $n\le P''_1\eta_-$ or $n\ge P''_2\eta_+$ above contribute $\ll\ep$. If 
$$P''_1\eta_+<n<P''_2\eta_-,$$ 
then, in a way similar to the proof of (10.),
$$
\int_{0.496}^{0.5}\bigg\{g\bigg(\frac{P^{0.5}Dt_0}{n}\bigg)-g\bigg(\frac{P^{z}Dt_0}{n}\bigg)\bigg\}\,dz=\frac{\log n}{\log P}-\tilde{\al}-0.496+O(\ep).
$$
Thus, in a way similar to the proof of Proposition 2.6, by the above discussion we deduce that
$$
\sum_{\psi\in\Psi_1}\,\sum_{\rho\in\tz(\psi)}\c^*\ss |Z(\rho,\pc)^{-1}\tilde{H}_{15}\ss-\bar{H}_{16}(1-s,\bp)|^2\omega(\rho)=o(\aa\p).\eqno(12.8)
$$
where
$$
\bar{H}_{16}(1-s,\bp)=\frac{1}{\log P_1}\sum_{P''_1<n<P''_2}\frac{\chi\bp(n)}{n^{1-s}}\big(n/P''_1\big)^{\be_6}\log\big(n/P''_1\big).
$$
By (12.7) and (12.8) we obtain
$$
\Xi_{15}=\sum_{\psi\in\Psi_1}\,\sum_{\rho\in\tz(\psi)}\c^*\ss \bar{H}_{16}(1-\rho,\bp)H_2\ss\omega(\rho)+o(\p).
$$
This implies, by Lemma 8.1,
$$
\Xi_{15}=\Theta_1(\bold{a}_{12}, \bold{a}_{25})+\overline{\Theta_1(\bold{a}_{15}, \bold{a}_{22})}+o(\p)\eqno(12.9)
$$ 
where
$$
a_{15}(n)=\chi(n)\vk_{13}(n)
$$
with
$$
\vk_{13}(n)=
\begin{cases}
(\log P_1)^{-1}(n/P''_1)^{-\be_6}\log(n/P''_1)&\quad\text{if}\quad P''_1<n<P''_2,\\
0&\quad\text{otherwise},
\end{cases}
$$
$$
a_{25}(n)=\overline{a_{15}(n)}.
$$

Assume $1\le j\le 3$ in what follows.

\medskip\noindent

{\bf Lemma 12.1.}\,\,{\it We have
$$
\sum_l\frac{\chi(l)\vk_{13}(dl)}{l^{1-\be_j}}\ll
\begin{cases} T^{-c}&\quad\text{if}\quad d\le P''_1/T,\\
\al_1&\quad\text{if}\quad P''_1/T<d\le P''_1,
\end{cases}
$$
and
$$
\sum_l\frac{\chi(l)\vk_{13}(dl)}{l^{1-\be_j}}=\frac{L'(1,\chi)}{\log P_1}\big(-1+(2\be_6-\be_j)\log(d/P''_1)+\epsilon_{1j}(d)\big),\quad |\epsilon_{1j}(d)|<10^{-5},
$$
if} $P''_1<d<P_2$,. 
\medskip\noindent

{\it Proof.} Note that $\vk_{13}(n)\ll\al_1$ if $P''_1/T\le n\le P''_1T$. In the case $d\le P''_1$ the results follow by the Polya-Vinogradov inequality and partial summation. On the other hand, in the case $P''_1<d<P_2$ the sum becomes
$$\aligned
\frac{(d/P''_1)^{-\be_6}\log(d/P''_1)}{\log P_1}&\sum_{l<P''_2/d}\frac{\chi(l)}{l^{1+\be_6-\be_j}}+\frac{(d/P''_1)^{-\be_6}}{\log P_1}\sum_{l<P''_2/d}\frac{\chi(l)\log l}{l^{1+\be_6-\be_j}}+O(T^{-c})\\
=&\frac{L'(1,\chi)(d/P''_1)^{-\be_6}}{\log P_1}\big(-1+(\be_6-\be_j)\log(d/P''_1)\big)+O(\l^{-15}).
\endaligned
$$
Since $P''_2/d>T^{10}$ and $(d/P''_1)^{-\be_6}=1-\be_6\log(d/P''_1)+...$,
 the result follows by Lemma 5.8 and a simple estimates. $\Box$
\medskip\noindent

{\bf Lemma 12.2.}\,\,{\it If $dr<P''_1/T$, then
$$
\sum_l\frac{\chi(l)\bar{\vk}_{13}(drl)\xi_{j}(l;d,r)}{l}=b^*L'(1,\chi)\Pi(d,r)(\log P)\be_{j+1}\be_{j+2}+O\eqno(12.10)
$$
with
$$
b^*=\frac{1}{0.504}\int_0^{0.004}ze^{3\pi iz/2}\,dz.
$$
If $P''_1/T\le n\le P''_1$, then}
$$
\sum_l\frac{\chi(l)\bar{\vk}_{13}(drl)\xi_{j}(l;d,r)}{l}\ll\al_1;\eqno(12.11)
$$

{\it Proof.}\,\,If $P''_1<n<P''_2$, then
$$
\bar{\vk}_{13}(n)=-\frac{1}{\log P_1}\frac{\partial}{\partial w}(n/P''_1)^{\be_6-w}|_{w=0}.
$$
Hence
$$
\sum_l\frac{\chi(l)\bar{\vk}_{13}(drl)\xi_{j}(l;d,r)}{l}=-\frac{1}{\log P_1}\frac{\partial}{\partial w}\bigg((dr/P''_1)^{\be_6-w}
\sum_{P''_1/dr<l<P''_2/dr}\frac{\chi(l)\xi_{j}(l;d,r)}{l^{1-\be_6+w}}\bigg)|_{w=0}.
$$
Suppose $dr\le P''_1/T$. By (4) and (4), for $|w|=\al$ we have
$$\aligned
\sum_{P''_1/dr<l<P''_2/dr}&\frac{\chi(l)\xi_{j}(l;d,r)}{l^{1-\be_6+w}}=\sum_l\frac{\chi(l)\xi_{j}(l;d,r)}{l^{1-\be_6+w}}\bigg\{
g\bigg(\frac{P''_2}{drl}\bigg)-g\bigg(\frac{P''_1}{drl}\bigg)\bigg\}+O(\l^{-15})\\
=&\tp\int_{(1)}\bigg(\sum_l\frac{\chi(l)\xi_{j}(l;d,r)}{l^{1-\be_6+w+s}}\bigg)\big((P''_2/dr)^s-(P''_1/dr)^s\big)\frac{\omega_1(s)\,ds}{s}+O(\l^{-15}).
\endaligned$$
In a way similar to the proof of Lemma 8.4, by lemma 8.2 and 5.8, we find that the right side above is equal to 
$$\aligned
L'(1,\chi)\Pi(d,r)&\cdot\tp\int_{|s|=5\al}\frac{(s+w-\be_6+\be_{j+1})(s+w-\be_6+\be_{j+2})}{s+w-\be_6}\,\frac{(P''_2/dr)^s-(P''_1/dr)^s}{s}\,ds\\
&+O(\l^{-15})\\
=&L'(1,\chi)\Pi(d,r)\be_{j+1}\be_{j+2}\frac{(P''_2/dr)^{\be_6-w}-(P''_1/dr)^{\be_6-w}}{\be_6-w}+O(\l^{-15}).
\endaligned
$$
It follows by Cauchy' integral formula that
$$
\frac{\partial}{\partial w}\bigg((dr/P''_1)^{\be_6-w}
\sum_{P''_1/dr<l<P''_2/dr}\frac{\chi(l)\xi_{j}(l;d,r)}{l^{1-\be_6+w}}\bigg)|_{w=0}=\frac{\partial}{\partial w}\bigg(\int_{1}^{P^{0.004}}y^{\be_6-w-1}\,dy\bigg)|_{w=0}+O(\l^{-6})
$$
The main term on the right side is, by substituting $y=P^z$, equal to
$$
-(\log P)^2\int_0^{0.004}ze^{3\pi iz/2}\,dz.
$$
Gathering these results together we obtain (12.10). The proof of (12.11) is similar to that of .

{\bf Lemma 12.3.}\,\,{\it If $P''_1<dr<P_2$, then}
$$\aligned
\sum_l\frac{\chi(l)\bar{\vk}_{13}(drl)\xi_{j}(l;d,r)}{l}=&\frac{L'(1,\chi)\Pi(d,r)}{\log P_1}\big(-1+(-2\be_6+\be_{j+1}+\be_{j+2})\log(dr/P''_1)+\epsilon_{2j}(dr)\big),\\
&\quad |\epsilon_{2j}(dr)|<10^{-5}.
\endaligned
$$

{\it Proof.}\,\, The left side is equal to
$$
-\frac{1}{\log P_1}\frac{\partial}{\partial w}\bigg((dr/P''_1)^{\be_6-w}
\sum_{l<P''_2/dr}\frac{\chi(l)\xi_{j}(l;d,r)}{l^{1-\be_6+w}}\bigg)|_{w=0}.
$$
Assume $|w|=\al$. In a way similar to the proof of Lemma 12.1, we deduce that
$$\aligned
\sum_{l<P''_2/dr}&\frac{\chi(l)\xi_{j}(l;d,r)}{l^{1-\be_6+w}}=L'(1,\chi)\Pi(d,r)\\
&\times\cdot\tp\int_{|s|=5\al}\frac{(s+w-\be_6+\be_{j+1})(s+w-\be_6+\be_{j+2})}{s+w-\be_6}\,\frac{(P''_2/dr)^s}{s}\,ds
+O(\l^{-15}).
\endaligned
$$
By direct calculation,
$$\aligned
\tp\int_{|s|=5\al}&\frac{(s+w-\be_6+\be_{j+1})(s+w-\be_6+\be_{j+2})}{s+w-\be_6}\,\frac{(P''_2/dr)^s}{s}\,ds\\
&=w-\be_6+\be_{j+1}+\be_{j+2}+\be_{j+1}\be_{j+2}\int_{1}^{P''/dr}y^{\be_6-w-1}\,dy,
\endaligned
$$
and the derivative of 
$$
(dr/P''_1)^{\be_6-w}\bigg(w-\be_6+\be_{j+1}+\be_{j+2}+\be_{j+1}\be_{j+2}\int_{1}^{P_2''/dr}y^{\be_6-w-1}\,dy\bigg)
$$
at $w=0$ is equal to
$$\aligned
(dr/P''_1)^{\be_6}&\bigg(1-\be_{j+1}\be_{j+2}\int_{1}^{P_2''/dr}y^{\be_6-1}\log y\,dy\bigg)\\
&-(dr/P''_1)^{\be_6}\log(dr/P''_1)\bigg(-\be_6+\be_{j+1}+\be_{j+2}+\be_{j+1}\be_{j+2}\int_{1}^{P_2''/dr}y^{\be_6-1}\,dy\bigg).
\endaligned
$$
This can be written as the form
$$
1-(-2\be_6+\be_{j+1}+\be_{j+2})\log(dr/P''_1)-\epsilon_{2j}(dr),
$$
since $(dr/P''_1)^{\be_6}=1+\be_6\log(dr/P''_1)+...$. $\Box$

\medskip\noindent

{\it Evaluation of $\Theta_1(\bold{a}_{12},\bold{a}_{25})$.}

\medskip\noindent

We have
$$\aligned
S_{j}(\bold{a}_{12},\bold{a}_{25})=&\sum_{r}\,\sum_d\frac{|\chi(d)||\mu\chi(r)|\la_{0j}(dr)}{dr\varphi(r)}\bigg(\sum_m\frac{\chi(m)\big(\bar{\iota}_3
\vk_3(drm)+\bar{\iota}_4\vk_2(drm)\big)}{m^{1-\be_j}}\bigg)\\
&\times\bigg(\sum_n\frac{\chi(n)\bar{\vk}_{13}(n)\xi_{0j}(n;d,r)}{n}\bigg).
\endaligned
$$
The sum is split into three sums according to $dr\le P_1''/T$, $P_1''/T< dr\le P''_1$ and $P_1''< dr<P_2$. By lemma 8.2, 8.3 and 12.1,the sum over $P_1''/T< dr\le P''_1$ contributes $o(\al)$; the sum over $dr\le P_1''/T$ is equal to
$$\aligned
L'(1,\chi)^2&b^*\be_{j+1}\be_{j+2}\sum_{n<P^{0.496}}\frac{|\chi(n)|\la_{0j}(n)}{\varphi(n)}
\bigg(\frac{\bar{\iota}_3\f_{j6}(P^{0.498}/n)}{0.498}+\frac{\bar{\iota}_4\f_{j7}(P^{0.5}/n)}{0.5}\bigg)+o(\al)\\
=&\aa b^*(\log P)\be_{j+1}\be_{j+2}\int_0^{0.496}\bigg(\frac{\bar{\iota}_3\ff_{j6}(0.498-z)}{0.498}+\frac{\bar{\iota}_4\ff_{j7}(0.5-z)}{0.5}\bigg)\,dz+o(\al).
\endaligned\eqno(12.12)
$$
By lemma 8.2, 8.3 and 12.3, the sum over $P_1''< dr<P_2$ is equal to
$$\aligned
&\frac{L'(1,\chi)^2\bar{\iota}_3}{(0.504)(0.498)\log^2 P}\sum_{P^{0.496}<n<P^{0.498}}\frac{|\chi(n)|\la_{0j}(n)}{\varphi(n)}
\f_{j6}(P^{0.498}/n)\w_j(n)
\\&+\frac{L'(1,\chi)^2\bar{\iota}_4}{(0.504)(0.5)\log^2 P}\sum_{P^{0.496}<n<P^{0.5}}\frac{|\chi(n)|\la_{0j}(n)}{\varphi(n)}
\f_{j7}(P^{0.5}/n)\w_j(n)+o(\al)\\
=&\frac{\aa\bar{\iota}_3}{(0.504)(0.498)\log P}\int_{0.496}^{0.498}\ff_{j6}(0.498-z)\w_j(P^z)\,dz
\\&+\frac{\aa\bar{\iota}_4}{(0.504)(0.5)\log P}\int_{0.496}^{0.5}\ff_{j7}(0.5-z)\w_j(P^z)\,dz+o(\al).
\endaligned\eqno(12.13)
$$
where
$$
\w_j(n)=-1+(-2\be_6+\be_{j+1}+\be_{j+2})\log(n/P''_1)+\epsilon_{2j}(n).
$$
For $0\le z\le 0.002$, a good approximation to $\ff_{j6}(z)$ is 
$$
1+\pi i(3-j)z\simeq 1+(2\be_6-\be_j)(\log P)z.
$$
Since 
$$
4\be_6-\be_1-\be_2-\be_3\simeq 0,
$$
it follows by simple calculation that
$$
\int_{0.496}^{0.498}\ff_{j6}(0.498-z)\w_j(P^z)\,dz=-0.002+\epsilon/10.
$$
For $0\le z\le 0.004$, a good approximation to $\ff_{j7}(z)$ is
$$
1+\pi i(5-j)z\simeq 1+(2\be_7-\be_j)(\log P)z.
$$
Thus, for $0.496\le z\le 0.5$, the function $\ff_{j7}(0.5-z)\w_j(P^z)$ can be well approximated by
$$
-1-(2\be_7-\be_j)(\log P)(0.5-z)-(2\be_6-\be_{j+1}-\be_{j+2})(\log P)(z-0.496).
$$
Since 
$$
(2\be_6+2\be_7-\be_1-\be_2-\be_3)\log P\simeq 2\pi i,
$$
it follows that
$$
\int_{0.496}^{0.5}\ff_{j6}(0.5-z)\w_j(P^z)\,dz=-0.004-\frac{\pi i}{250^2}+\epsilon/10.
$$
inserting these results into (12.13), we find that the sum over $P_1''< dr<P_2$ is equal to
$$
\frac{\aa}{0.504\log P}\bigg(-\frac{0.002\bar{\iota}_3}{0.498}-0.008\bar{\iota}_4-\frac{2\pi i\bar{\iota}_4}{250^2}+\epsilon/4\bigg).
\eqno(12.14)$$
It follows from (12.12) and (12.14) that
$$
\frac{1}{2\al}S_{1}(\bold{a}_{12},\bold{a}_{25})+\frac{2}{\al}S_{2}(\bold{a}_{12},\bold{a}_{25})+\frac{3}{2\al}S_{3}(\bold{a}_{12},\bold{a}_{25})=\aa(e_1^*+e_2^*+\epsilon/2)\eqno(12.15)
$$
where
$$
e_1^*=-\pi b^*(3e_{11}^*+6e_{12}^*+3e_{13}^*)
$$
with
$$
e_{1j}^*=\int_0^{0.496}\bigg(\frac{\bar{\iota}_3\ff_{j6}(0.498-z)}{0.498}+\frac{\bar{\iota}_4\ff_{j7}(0.5-z)}{0.5}\bigg)\,dz,
$$
and
$$
e_2^*=\frac{4}{0.504\pi}\bigg(-\frac{0.002\bar{\iota}_3}{0.498}-0.008\bar{\iota}_4-\frac{2\pi i\bar{\iota}_4}{250^2}\bigg).
$$

\medskip\noindent

{\it Evaluation of $\Theta_1(\bold{a}_{15},\bold{a}_{22})$.}

\medskip\noindent

We have
$$\aligned
S_{j}(\bold{a}_{15},\bold{a}_{22})=&\sum_{r}\,\sum_d\frac{|\chi(d)||\mu\chi(r)|\la_{0j}(dr)}{dr\varphi(r)}\bigg(\sum_m\frac{\chi(m)\vk_{13}(m)}{m^{1-\be_j}}\bigg)\\
&\times\bigg(\sum_n\frac{\chi(n)\big(\iota_3
\bar{\vk}_3(drn)+\iota_4\bar{\vk}_2(drn)\big)\xi_{0j}(n;d,r)}{n}\bigg).
\endaligned
$$
The sum is split into two sums according to  $ dr\le P''_1$ and $P_1''< dr<P_2$. By lemma 8.2, 8.3 and 12.1, the first sum contributes $o(\al)$; the second sum  is equal to
$$
\aa b^*\int_0^{0.496}\big(\bar{\iota}_3\ff_j6(0.498-z)+\bar{\iota}_4\ff_j6(0.5-z)\big)\,dz+o(\al).
$$
By lemma 8.2, 8.3 and 12.3, the sum over $P_1''< dr<P_2$ is equal to
$$\aligned
&\frac{L'(1,\chi)^2\iota_3}{(0.504)(0.498)\log^2 P}\sum_{P^{0.496}<n<P^{0.498}}\frac{|\chi(n)|\la_{0j}(n)}{\varphi(n)}
\g_{j6}(P^{0.498}/n)\w^*_j(n)
\\&+\frac{L'(1,\chi)^2\iota_4}{(0.504)(0.5)\log^2 P}\sum_{P^{0.496}<n<P^{0.5}}\frac{|\chi(n)|\la_{0j}(n)}{\varphi(n)}
\g_{j7}(P^{0.5}/n)\w^*_j(n)+o(\al)\\
=&\frac{\aa\iota_3}{(0.504)(0.498)\log P}\int_{0.496}^{0.498}\gh_{j6}(0.498-z)\w^*_j(P^z)\,dz
\\&+\frac{\aa\iota_4}{(0.504)(0.5)\log P}\int_{0.496}^{0.5}\gh_{j7}(0.5-z)\w^*_j(P^z)\,dz+o(\al).
\endaligned
$$
where
$$
\w^*(n)=-1+(2\be_6-\be_j)\log(n/P''_1)+\epsilon_{1j}(n).
$$
For $0\le z\le 0.004$ and $\mu=6,7$, the function $\gh_{j\mu}(z)$ is well-approximated by
$$1-(2\be_\mu-\be_{j+1}-\be_{j+2})(\log P)z.$$
Thus, in a way similar to the evaluation of $\Theta_1(\bold{a}_{12},\bold{a}_{25})$, we deduce that
$$
\frac{1}{2\al}S_{1}(\bold{a}_{15},\bold{a}_{23})+\frac{2}{\al}S_{2}(\bold{a}_{15},\bold{a}_{22})+\frac{3}{2\al}S_{3}(\bold{a}_{15},\bold{a}_{22})=\aa(\overline{e_2^*}+\epsilon/2)\eqno(12.16).
$$

Finally, by (12.9), (12.15) and (12.16) we conclude
$$
\Xi_{15}=(e_1^*+2e_2^*+\epsilon)\aa\p.\eqno(12.17)
$$

\medskip\noindent

 \section  *{\centerline{13.  Approximation to $\Xi_{14}$}}

\medskip\noindent

In this section we establish an approximation to $\Xi_{14}$.

Assume that $\psi\in\Psi_1$ and $\rho\in\z(\psi)$. By Lemma 5.2 and (2.2),
$$\aligned
\c^*\ss&=-iZ(\rho+\be_3,\psi)^{-1}\frac{L(\rho+\be_1,\psi)L(\rho+\be_2,\psi)L(\rho+\be_3,\psi)}{L'\ss}(1+O(1/t_0))\\
&=-i\frac{L(\rho+\be_1,\psi)L(\rho+\be_2,\psi)L(1-\rho-\be_3,\bp)}{L'\ss}(1+O(\l^{-123})).
\endaligned\eqno(13.1)
$$
By Lemma 6.1,
$$
L(1-\rho-\be_3,\bp)=\overline{L(\rho+\be_3,\psi)}=K(1-\rho-\be_3,\bp)+Z(\rho+\be_3,\psi)^{-1}N(\rho+\be_3,\psi)
+O(E_1(\rho+\be_3,\psi)),
$$
and, by Lemma 5.1,
$$
Z(\rho+\be_3,\psi)^{-1}=Z(\rho+\be_2,\psi)^{-1}(pt_0)^{\be_1}+O(\l^{-123}).
$$
Hence
$$\aligned
L(\rho&+\be_2,\psi)L(1-\rho-\be_3,\bp)\\=&L(\rho+\be_2,\psi)K(1-\rho-\be_3,\bp)
+(pt_0)^{\be_1}N(\rho+\be_3,\psi)L(1-\rho-\be_2,\bp)\\
&+O(|L(\rho+\be_2,\psi)N(\rho+\be_3,\psi)\l^{-123})+O(|L(\rho+\be_2,\psi)|E_1(\rho+\be_3,\psi)).
\endaligned\eqno(13.2)
$$
By Lemma 6.1 and 5.1,
$$\aligned
(pt_0)^{\be_1}L(1-\rho-\be_2,\bp)=&(pt_0)^{\be_1}K(1-\rho-\be_2,\bp)+(pt_0)^{\be_3}Z(\rho,\psi)^{-1}N(\rho+\be_2,\psi)\\
&+O(N(\rho+\be_2,\psi)\l^{-123})+O(E_1(\rho+\be_2,\psi)).\endaligned
$$
We insert this into (13.2) and then insert the result into (13.1). Thus we obtain
$$\aligned
\c^*\ss=&
-i\frac{L(\rho+\be_1,\psi)L(\rho+\be_2,\psi)}{L'\ss}K(1-\rho-\be_3,\bp)\\
&-i(pt_0)^{\be_1}\frac{L(\rho+\be_1,\psi)}{L'\ss}N(\rho+\be_3,\psi)K(1-\rho-\be_2,\bp)\\
&-i(pt_0)^{\be_3}Z\ss^{-1}\frac{L(\rho+\be_1,\psi)}{L'\ss}N(\rho+\be_2,\psi)N(\rho+\be_3,\psi)\\
&+O(E_1^*\ss)
\endaligned\eqno(13.3)
$$
where
$$\aligned
E_1^*\ss=&\bigg|\frac{L(\rho+\be_1,\psi)L(\rho+\be_2,\psi)}{L'\ss}\bigg|\big(|L(\rho+\be_3,\psi)|t_0^{-1}+|N(\rho+\be_3,\psi)|
\l^{-123}+E_1(\rho+\be_3,\psi)\big)\\
&+\bigg|\frac{L(\rho+\be_1,\psi)N(\rho+\be_3,\psi)}{L'\ss}\bigg|\big(|N(\rho+\be_2,\psi)|\l^{-123}+E_1(\rho+\be_2,\psi)\big).
\endaligned$$
Recall that $B\sp$ is given by (12.2). Multiplying (13.3) by $Z(\rho,\pc)^{-1}B\ss$ and applying Lemma 4.8 we obtain
$$
\c^*\ss Z(\rho,\pc)^{-1}B\ss=-i\big(\k^*_1\ss+(pt_0)^{\be_1}\k^*_2\ss-(pt_0)^{\be_3}\k^*_3\ss\big)+O(E_1^*\ss|B\ss|)
$$
where
$$
\k^*_1\ss=Z(\rho,\pc)^{-1}\frac{L(\rho+\be_1,\psi)L(\rho+\be_2,\psi)}{L'\ss}B\ss K(1-\rho-\be_3,\bp),\eqno(13.4)
$$
$$
\k^*_2\ss=Z(\rho,\pc)^{-1}\frac{L(\rho+\be_1,\psi)}{L'\ss}B\ss N(\rho+\be_3,\psi)K(1-\rho-\be_2,\bp),\eqno(13.5)
$$
$$
\k^*_3\ss=\frac{L(\rho+\be_1,\psi)}{L'\ss}B\ss G\ss N(\rho+\be_2,\psi) N(\rho+\be_3,\psi)F(1-\rho,\bp),\eqno(13.6)
$$
Inserting this into (12.4) we obtain
$$
\Xi_{14}=-i\big(\Phi_1+\Phi_2-\Phi_3\big)+O(\e)\eqno(13.7)
$$
where
$$
\Phi_1=\sum_{\psi\in\Psi_1}\,\sum_{\rho\in\z(\psi)}\k^*_1\ss\omega(\rho),\eqno(13.8)
$$
$$
\Phi_2=\sum_{\psi\in\Psi_1}(p_\psi t_0)^{\be_1}\,\sum_{\rho\in\z(\psi)}\k^*_2\ss\omega(\rho),\eqno(13.9)
$$
$$
\Phi_3=\sum_{\psi\in\Psi_1}(p_\psi t_0)^{\be_3}\,\sum_{\rho\in\z(\psi)}\k^*_3\ss\omega(\rho),\eqno(13.10)
$$
and
$$
\e=\sum_{\psi\in\Psi_1}\,\sum_{\rho\in\z(\psi)}E_1^*\ss |B\ss|\omega(\rho).
$$
Combining  (2.34), Cauchy's inequality, Proposition 7.1,  Lemma 5.9, 6.1 and 3.3, we can verify that
$$
\e=o(\p).\eqno(13.11)
$$
For example, by (2.34),
$$\aligned
\sum_{\psi\in\Psi_1}&\,\sum_{\rho\in\z(\psi)}\bigg|\frac{L(\rho+\be_1,\psi)}{L'\ss}\bigg||L(\rho+\be_2,\psi)|^2\omega(\rho)\\
=&-\frac{1}{2\pi}\sum_{\psi\in\Psi_1}\bigg(\int_{\j(\al)}-\int_{\j(-\al)}\bigg)\frac{M(s+\be_1,\psi)}{M\sp}L(s+\be_2,\psi)L(1-s-\be_2,\bp)
\omega(s)\,ds+O(\ep),
\endaligned
$$
the right side being estimated via  Lemma 5.9, 6.1 and 3.3.

 \section  *{\centerline{14. Mean-value formula II}}
 \medskip\noindent
 
Recall that we always assume $\psi$ is a primitive character $(\bmod\,p)$, $p\sim P$. Sometimes we write $p_\psi$ for the modulus $p$.

Let $\bold{k}^*=\{\kappa^*(m)\}$ and $\bold{a}^*=\{a^*(n)\}$ denote sequences of complex numbers satisfying
$$
\kappa^*(m)\ll\tau_5(m),\eqno (14.1)
$$
$$\quad a^*(n)\ll 1,\quad a^*(n)=0\quad\text{if}\quad n>2P_4.\eqno(14.2)
$$
For $\be\in\mathbf{C}$ write
$$
\Theta_2(\be, \bold{k}^*, \bold{a}^*)=\sum_{\psi\in\Psi_1}\frac{(p_\psi t_0)^\be}{2\pi i}\int_{\j(1)}Z(s,\pc)^{-1}\bigg(\sum_m\frac{\kappa^*(m)\psi(m)}{m^s}\bigg)\bigg(\sum_n\frac{a^*(n)\bp(n)}{n^{1-s}}\bigg)\omega(s)\,ds.
$$ 

The goal of this section is to prove 
\medskip\noindent

{\bf Proposition 14.1.}\,\, {\it Suppose $|\be|<5\al$. Then}
$$
\Theta_2(\be, \bold{k}^*, \bold{a}^*)=\frac{1}{\varphi(D)}\sum_{p\sim P}(p t_0)^\be\,\sum_d\frac{1}{d}\,\sum_k\frac{\mu\chi(k)a^*(dk)}{k\varphi(k)}\sum_{(l,k)=1}\chi(l)\kappa^*(dl)\Delta\bigg(\frac{l}{Dpk}\bigg)+o(\p).
$$
\medskip\noindent

{\it  Proof.}\,\,  We prove this proposition with $\be=0$ only, as the general case is almost identical. Write $\Theta_2$ for $\Theta_2(0, \bold{k}^*, \bold{a}^*)$. Similar to the proof of Proposition 7.1, in the expression for $\Theta_2$, we can extend the sum over $\Psi_1$ to the sum over $\Psi$, with acceptable errors. Namely we have
$$
\Theta_2=\sum_{\psi\in\Psi}\tilde{I}_2(\psi)+o(\p)\eqno(14.3)
$$
where
$$
\tilde{I}_2(\psi)=\frac{1}{2\pi i} \int_{\j(1)}Z(s,\pc)^{-1}\bigg(\sum_m\frac{\kappa^*(m)\psi(m)}{m^s}\bigg)\bigg(\sum_n\frac{a^*(n)\bp(n)}{n^{1-s}}\bigg)\omega(s)\,ds.
$$
Assume $\psi(\bmod\,p)\in\Psi$. We use (2.5) with $\theta=\pc$ and  then replace the segment $\j(1)$ by the vertical line $\sigma=3/2$ with a negligible error. Thus,  by integration term by term,
$$
\tilde{I}_2(\psi)=\frac{\tau(\chi\bp)}{Dp}\sum_m\,\sum_n\frac{\kappa^*(m)a^*(n)\psi(m)\bp(n)}{n}\Delta_1\bigg(\frac{m}{Dpn}\bigg)
+O(\ep).$$
 For $(mn,p)=1$ we have
$$
\sum_{\psi(\bmod\,p)}\tau(\chi\bp)\psi(m)\bp(n)=\tau(\chi)\chi(p)(p-1)e\bigg(\frac{m\overline{Dn}}{p}\bigg).
$$
Note that $(n,p)=1$ if $n<2P_4$ and $|\tau(\chi\psi_p^0)|=\sqrt{D}$. Hence, substituting $d=(m,n)$, $m=dl$, $n=dk$ and inverting the order of summation we obtain
$$\aligned
\sideset{}{^*}\sum_{\psi(\bmod\,p)}\tilde{I}_2(\psi)=& \frac{\tau(\chi)\chi(p)}{D}\sum_m\,\sum_n\frac{\kappa^*(m)a^*(n)}{n}\Delta_1\bigg(\frac{m}{Dpn}\bigg)e\bigg(\frac{m\overline{Dn}}{p}\bigg)+O(PT^{-c})\\
=& \frac{\tau(\chi)\chi(p)}{D}\sum_d\frac{1}{d}\,\sum_k\frac{a^*(dk)}{k}\sum_{(l,k)=1}\kappa^*(dl)\Delta_1\bigg(
\frac{l}{Dpk}\bigg)e\bigg(\frac{l\overline{Dk}}{p}\bigg)+O(PT^{-c}).
\endaligned
$$
(here we have replaced the constraint $(l,pk)=1$ by $(l,k)=1$ with an acceptable error, and used  a bound for $|\Delta_1(x)|=|\Delta(x)|$ given by Lemma  5.3). Since
$$
\frac{\overline{Dk}}{p}\equiv -\frac{\bar{p}}{Dk}+\frac{1}{Dpk}\quad(\bmod\,1),
$$
it follows that
$$
\Delta_1\bigg(\frac{l}{Dpk}\bigg)e\bigg(\frac{l\overline{Dk}}{p}\bigg)=
\Delta\bigg(\frac{l}{Dpk}\bigg)e\bigg(-\frac{l\bar{p}}{Dk}\bigg)
$$
with $\bar{p}p\equiv 1(\bmod\,Dk)$. Hence
$$\aligned
\sideset{}{^*}\sum_{\psi(\bmod\,p)}\tilde{I}_2(\psi)
=&\frac{\tau(\chi)\chi(p)}{D}\sum_d\frac{1}{d}\,\sum_k\frac{a^*(dk)}{k}\sum_{(l,k)=1}\kappa^*(dl)\Delta\bigg(\frac{l}{Dpk}\bigg)e\bigg(-\frac{l\bar{p}}{Dk}\bigg)\\
&+O(PT^{-c}).
\endaligned
\eqno(14.4)
$$
Assume $(k,p)=1$. If $(l,k)=1$ and $(l,D)=D_1$, then $(k, D_1)=1$. Thus, substituting $(l,D)=D_1$ and $l=D_1l_1$, we find that
the innermost sum in (14.4) is equal to
$$
\sum_{\substack{D=D_1D_2\\(D_1,k)=1}}
\sum_{(l_1,D_2k)=1}(\kappa_1*b)(D_1dl_1)\Delta\bigg(\frac{l_1}{D_2pk}\bigg)e\bigg(-\frac{l_1\bar{p}}{D_2k}\bigg).
$$
Inserting this into (14.4) and rearranging the terms we conclude
$$\aligned
\sideset{}{^*}\sum_{\psi(\bmod\,p)}\tilde{I}_2(\psi)
=\frac{\tau(\chi)\chi(p)}{D}\sum_{D=D_1D_2}\s(D_1,D_2;p)+O(PT^{-c}),
\endaligned$$
where
$$
\s(D_1,D_2;p)=
\sum_{d}\frac{1}{d}\sum_{(k, D_1)=1}\frac{a^*(dk)}{k}\sum_{(l,D_2k)=1}\kappa^*(D_1dl)\Delta\bigg(\frac{l}{D_2pk}\bigg)
e\bigg(-\frac{l\bar{p}}{D_2k}\bigg)
$$
(here we have rewritten $l$ for $l_1$). 
Since $\chi(-1)\tau(\chi)^2=D$, the proof of Proposition 14.1 is now reduced to showing that
$$\aligned
\sum_{p\sim P}\chi(p)&\s(1,D;p)\\
=&\frac{\chi(-1)\tau(\chi)}{\varphi(D)}\sum_d\frac{1}{d}\,\sum_k\frac{\mu\chi(k)a^*(dk)}{k\varphi(k)}\sum_{(l,k)=1}\chi(l)\kappa^*(dl)\Delta\bigg(\frac{l}{Dpk}\bigg)\\
&+O(P^2D^{1/2-c})
\endaligned\eqno(14.5)
$$
and, for $D=D_1D_2$, $D_1>1$,
$$
\sum_{p\sim P}\chi(p)\s(D_1,D_2;p)\ll P^2D^{1/2-c}.\eqno(14.6)
$$

We first prove (14.5). Assume $k<2P_4$ and $(l,Dk)=1$. Then
$$
e\bigg(-\frac{l\bar{p}}{Dk}\bigg)=\frac{1}{\varphi(Dk)}\sum_{\theta(\bmod\,Dk)}\tau(\bar{\theta})\theta(-l)\bar{\theta}(p).
$$
Hence
$$
\s(1,D;p)=\sum_{d}\frac{1}{d}\sum_{k}\frac{a^*(dk)}{k\varphi(Dk)}\sum_{\theta(\bmod\,Dk)}\tau(\bar{\theta})\bar{\theta}(p)\sum_{l}\kappa^*(dl)\theta(-l)\Delta\bigg(\frac{l}{Dpk}\bigg)\eqno(14.7)
$$
(here we have remove the constraint $(l,Dk)=1$ as it is superfluous). 

If $\theta$ is the principal character $\psi_{Dk}^0$, then $\tau(\bar{\theta})=\mu(Dk)$. The total contribution from the $\theta=\psi_{Dk}^0$ term to the right side above is, by Lemma 5.3,
$$
\ll \sum_{d}\frac{1}{d}\sum_{k}\sum_{l}\frac{|(\kappa_1*b)(dl)a^*(dk)|}{\varphi(Dk)k}\bigg|\Delta\bigg(\frac{l}{Dpk}\bigg)\bigg|\ll P\l^c
$$
which is admissible for (14.5). Now,  let $\theta_k^1$ be the character $(\bmod\,Dk)$ induced by the primitive character $\chi(\bmod\,D)$ which is real. Then
$$
\tau(\theta_k^1)\theta_k^1(p)\theta_k^1(-l)=\tau(\chi)\mu\chi(k)\chi(-pl).
$$
Hence
$$\aligned
\sum_{d}\frac{1}{d}&\sum_{k}\sum_{(l,Dk)=1}\frac{\kappa^*(dl)a^*(dk)}{\varphi(Dk)k}\Delta\bigg(\frac{l}{Dpk}\bigg)
\tau(\theta_k^1)\theta_k^1(p)\theta_k^1(-l)\\
&=\frac{\chi(-p)\tau(\chi)}{\varphi(D)}\sum_d\frac{1}{d}\,\sum_k\frac{\mu\chi(k)a^*(dk)}{k\varphi(k)}\sum_{(l,k)=1}\chi(l)\kappa^*(dl)\Delta\bigg(\frac{l}{Dpk}\bigg).
\endaligned
$$
By the above discussion, to complete the proof of (14.5), it now suffices to show that the total contribution from the terms with $\theta\ne\psi_{Dk}^0$ and $\theta\ne\theta_k^1$  in (14.7) to the left side of (14.5) is $O(P^2D^{-c})$. Namely we need to prove
$$\aligned
\sum_d\frac{1}{d}\,&\sum_k\frac{|a^*(dk)|}{\varphi(Dk)k}\sideset{}{'}\sum_{\substack{\theta(\bmod\,Dk)\\
\theta\ne\theta_k^1}}|\tau(\bar{\theta})|\bigg|\sum_{l}\kappa^*(dl)\theta(l)\sum_{p\sim P}\chi\bar{\theta}(p)\Delta\bigg(\frac{l}{Dpk}\bigg)\bigg|\\
&\ll P^2D^{-c}.
\endaligned\eqno(14.8)
$$
If $\theta(\bmod\,Dk)$ is induced by a primitive character $\theta^*(\bmod\,r)$, then $r|Dk$ and $|\tau(\bar{\theta})|\le\sqrt{r}$. For $r|Dk$ we have 
$$
\frac{Dk}{r}\equiv 0\bigg(\bmod\,\frac{D}{(D,r)}\bigg).
$$
Thus, substituting $Dk=hr$, we see that the left side of (14.8) is
$$\aligned
\ll\sum_d\frac{1}{d}\,\sum_{1<r<2DP_4}\,\,\sum_{\substack{h<P/r\\h\equiv 0(D/(D,r))}}\frac{D}{\varphi(hr)h\sqrt{r}}
\sideset{}{^*}\sum_{\substack{\theta(\bmod\,r)\\
\theta\ne\chi}}\bigg|\sum_{(l,h)=1}\kappa^*(dl)\theta(l)\sum_{p\sim P}\chi\bar{\theta}(p)\Delta\bigg(\frac{l}{phr}\bigg)\bigg|.
\endaligned
$$
 The range for  $r$ is divided into two parts according to $1<r<D^3$ and $D^3\le r<2DP_4$.
In a way similar to the proof of Proposition 7.1,  for  $1<r<D^3$ we use the Mellin transform, Lemma 5.4 (i) and Lemma 5.6; 
for  $D^3\le r<2DP_4$ we use the Mellin transform, Lemma 5.4 (i)  and the large sieve inequality. Thus we obtain (14.5).

The proof of (14.6) is analogous.  In the case $D=D_1D_2$, $D_1>1$, the major  difference from the proof of (14.5) is that  the constraint $(k, D_1)=1$ 
 implies $D\nmid D_2k$, so that any non-principal character $\theta(\bmod\,D_2k)$ can not be induced by the real character $\chi(\bmod\,D)$. Thus, when the left side of (14.6) is treated in much the same way as (14.5), the main terms involved  in the proof of (14.5) do not appear. This yields (14.6).   $\Box$

 \section  *{\centerline{15. Evaluation of $\Phi_1$}}

 Recall that $\Phi_1$ is given by (13.8).  In view of (12.2), $B\sp$ can be written as
$$
B\sp=\sum_n\frac{b(n)\pc(n)}{n^s}\eqno(15.1)$$
with 
$$
b(n)\ll\tau_2(n), \qquad b(n)=0\quad\text{if}\quad n>PT^{-2}\eta_+.\eqno(15.2)
$$

 Write
 $$
\k_1\sp=Z(s,\pc)^{-1}\frac{L(s+\be_1,\psi)L(s+\be_2,\psi)}{L\sp}B\sp K(1-s-\be_3,\bp).
$$

 Similar to (8.1), for $\psi\in\Psi_1$ we have
 $$
 \sum_{\rho\in\z(\psi)}\k^*_1\ss\omega(\rho)=I_{2}^+(\psi)-I_{2}^-(\psi)+O(\ep)
 $$
 where
 $$
 I_{2}^\pm(\psi)=\tp\int_{\j(\pm\al)}\k_1\sp\omega(s)\,ds.
 $$
 Hence
 $$
 \Phi_1=\sum_{\psi\in\Psi_1}\big(I_{2}^+(\psi)-I_{2}^-(\psi)\big)+O(\ep).\eqno(15.3)
 $$
 
 First we prove that
 $$
 \sum_{\psi\in\Psi_1}I_{2}^-(\psi)=o(\p).\eqno (15.4)
 $$
 
 Assume $\psi\in\Psi_1$ and $s\in\j(-\al)$. By (2.2) and Lemma 5.1,
 $$
 \frac{L(s+\be_1,\psi)L(s+\be_2,\psi)}{L(s,\psi)}=(pt_0)^{-\be_3}Z\sp \frac{L(1-s-\be_1,\bp)L(1-s-\be_2,\bp)}{L(1-s,\bp)}(1+O(\l^{-123})).
 $$
 In a way similar to the proof of Lemma 8.1, it can be shown that the total contribution of the $O(\l^{-123})$ term to the left side of (15.4) is $o(\p)$. On the other hand, by (2.4), (2.5) and the relation
$$
\tau(\pc)=\tau(\chi)\tau(\psi)\psi(D)\chi(p),
$$
we have
 $$
\frac{Z\sp}{Z(s,\pc)}=\tau(\chi)\chi(p)\bp(D)D^{s-1}(1+O(e^{-\pi t})).
$$
 For $\sigma<0$ we can write
 $$
 \frac{L(1-s-\be_1,\bp)L(1-s-\be_2,\bp)K(1-s-\be_3,\bp)}{L(1-s,\bp)}=\sum_m\frac{\tilde{k}(m)\bp(m)}{m^{1-s}}
 $$
 with $\tilde{k}(m)\ll\tau_4(m)$. Hence, moving the segment $\j(-\al)$ to $\j(-1)$  with a negligible error, we find that
 $$
 \sum_{\psi\in\Psi_1} I_{2}^-(\psi)=\tau(\chi) \sum_{\psi\in\Psi_1}\chi(p_\psi)(p_\psi t_0)^{-\be_3}\cdot\tp\int_{\j(-1)}\bigg(\sum_m\frac{\tilde{k}(m)\bp(Dm)}{(Dm)^{1-s}}\bigg)B\sp\omega(s)\,ds+o(\p).
 $$
 In a way similar to the proof of (7.3), the sum over $\psi\in\Psi_1$ can be extended to the sum over $\psi\in\Psi$, and  then the segment $\j(-1)$ can be replaced by the line $\sigma=-1/2$, with  acceptable errors. Thus the right side above is equal to
 $$
 \tau(\chi) \sum_{p\sim P}\chi(p)(pt_0)^{-\be_3}\sideset{}{^*}\sum_{\psi(\bmod\,p)}\tp\int_{(-1/2)}\bigg(\sum_m\frac{\tilde{k}(m)\bp(Dm)}{(Dm)^{1-s}}\bigg)B\sp\omega(s)\,ds+o(\p).
$$
Since
 $$
\tp\int_{(-1/2)}\frac{\omega(s)\,ds}{n^s(Dm)^{1-s}}=\frac{1}{Dm}\bigg(\frac{Dm}{n}\bigg)^{s_0}\exp\bigg\{-\l_2^2\log^2\frac{Dm}{n}\bigg\}
$$
 and, for $n<PT^{-5}$ and $Dm<2n$,
 $$
\sideset{}{^*}\sum_{\psi(\bmod\,p)}\psi(n)\bp(Dm)\ll
\begin{cases} p&\quad\text{if}\quad n=Dm,\\
1&\quad\text{if}\quad n\ne Dm,
\end{cases}
$$
it follows that
$$
\sideset{}{^*}\sum_{\psi(\bmod\,p)}\tp\int_{(-1/2)}\bigg(\sum_m\frac{\tilde{k}(m)\bp(Dm)}{(Dm)^{1-s}}\bigg)B\sp\omega(s)\,ds
\ll PD^{-4/5}.
$$
This yields (15.4).

 Let $\kappa_1(m)$ be given by
$$
\frac{\zeta(s+\be_1)\zeta(s+\be_2)}{\zeta(s)}=\sum_m\frac{\kappa_1(m)}{m^s},\quad \sigma>1.
$$
Regarding $b$ as an arithmetic function, for $\sigma>1$ we have
$$
\frac{L(s+\be_1,\psi)L(s+\be_2,\psi)B\sp}{L\sp}=\sum_m\frac{(\kappa_1*b)(m)\psi(m)}{m^s}.
$$
 On the other hand, we can write
 $$
 K(1-s-\be_3,\bp)=\sum_n\frac{\tg_3(n)\bp(n)}{n^{1-s}}
$$
with
 $$
\tg_3(y)=y^{\be_3}g^*(P_4/y).
$$
 Hence, moving the segment $\j(\al)$ to $\j(1)$ with a negligible error, we obtain
 $$
 \sum_{\psi\in\Psi_1}I_{2}^+(\psi)=\Theta_2(0,\bold{k}^*_1,\bold{a}_1^*)+O(\ep)\eqno(15.5)
 $$
 with
$$
 \kappa^*_1(m)=(\kappa_1*b)(m),\quad a^*_1(n)=\tg_3(n).
 $$
 It follows by (15.3)-(15.5) and Proposition 14.1 that
 $$
\Phi_1=\sum_{p\sim P}\Phi_1(p)+o(\p)
\eqno(15.6)
$$
where
$$
\Phi_1(p)=\frac{1}{\varphi(D)}\sum_k\frac{\mu\chi(k)}{k\varphi(k)}\sum_d\frac{\tg_3(dk)}{d}\,\sum_{(l,k)=1}(\kappa_1*b)(dl)\chi(l)\Delta\bigg(\frac{l}{Dpk}\bigg). \eqno(15.7)
$$

Assume $dk<2P_4$. Similar to (7.17),
$$
(\kappa_1*b)(dl)=\sum_{d=d_1d_2}\,\sum_{\substack{l=l_1l_2\\(l_1,d_2)=1}}b(d_2l_2)\kappa_1(d_1l_1).
$$
This yields
$$\aligned
\sum_{(l,k)=1}(\kappa_1*b)(dl)\chi(l)\Delta\bigg(\frac{l}{Dpk}\bigg)
=\sum_{d=d_1d_2}\,\sum_{(l_2,k)=1}b(d_2l_2)\chi(l_2)\sum_{(l_1,d_2k)=1}\kappa_1(d_1l_1)\chi(l_1)\Delta\bigg(\frac{l_1l_2}{Dpk}\bigg).
\endaligned
$$
 The innermost sum above is, by the Mellin transform, equal to
$$
\tp\int_{(2)}\bigg(\sum_{(l_1,d_2k)=1}\frac{\kappa_1(d_1l_1)\chi(l_1)}{l_1^s}\bigg)\bigg(\frac{Dpk}{l_2}\bigg)^s\delta(s)\,ds.
\eqno(15.8)
$$
Every $l_1$  can be uniquely written as $l_1=hl$ such that $h\in\n(d_1)$ and $(l, d_1)=1$, so that $\kappa_1(d_1l_1)=\kappa_1(d_1h)\kappa_1(l)$.  Hence, for $\sigma>1$,
$$\aligned
\sum_{(l_1,d_2k)=1}\frac{\kappa_1(d_1l_1)\chi(l_1)}{l_1^s}=&\tk_{1}(d_1;d_2k,s)\sum_{(l,d_1d_2k)=1}\frac{\kappa_1(l)\chi(l)}{l^s}
\\=&\tk_{1}(d_1;d_2k,s)\la_1(d_1d_2k,s)\frac{L(s+\be_1,\chi)L(s+\be_2,\chi)}{L(s,\chi)}
\endaligned
$$
where
$$
\tk_{1}(d_1;r,s)=\sum_{\substack{h\in\n(d_1)\\(h, r)=1}}\frac{\kappa_1(d_1h)\chi(h)}{h^s}\eqno(15.9)
$$
and
$$
\la_1(n,s)=\prod_{q|n}\frac{\big(1-\chi(q)q^{-s-\be_1}\big)\big(1-\chi(q)q^{-s-\be_2}\big)}{1-\chi(q)q^{-s}}.\eqno(15.10)
$$
 For notational simplicity we shall write $\tk_1(d;r)$ and $\la_1(n)$  for  $\tk_1(d;r,1)$ and $\la_1(n,1)$ respectively. Note that $Dpk/l_2>T$ if $l_2<PT^{-2}$. In a way similar to the treatment of (7.19), by Lemma 5.5,
we see that  the  expression (15.8) is equal to the residue of the integrand at $s=\tilde{\rho}$ plus an acceptable error. Further, by (5.15), in the expression for this residue, we can replace $\tilde{\rho}$ by $1$ with an acceptable error.  Thus we have
$$
\aligned
\sum_{(l_1,d_2k)=1}\kappa_1(d_1l_1)\chi(l_1)\Delta\bigg(\frac{l_1l_2}{Dpk}\bigg)=\r_1^*\bigg(\frac{Dpk}{l_2}\bigg)\tk_{1}(d_1;d_2k)\la_1(d_1d_2k)+O\bigg(\frac{\al^{100}\tau_3(d_1)Dpk}{l_2}\bigg),
\endaligned
$$
where
$$
\r_1^*=\frac{L(1+\be_1,\chi)L(1+\be_2,\chi)}{L'(1,\chi)}\delta(1).
$$
This yields
$$\aligned
\sum_{(l,k)=1}(\kappa_1*b)(dl)\chi(l)\Delta\bigg(\frac{l}{Dpk}\bigg)=&\r_1^* Dpk
\sum_{d=d_1d_2}\tk_{1}(d_1;d_2k)\la_1(dk)
\sum_{(l_2,k)=1}\frac{b(d_2l_2)\chi(l_2)}{l_2}\\&
+O\big(\al^{50}\tau_3(d_1)Dpk\big).
\endaligned
$$
Hence
$$\aligned
\sum_d\frac{\tg_3(dk)}{d}&\sum_{(l,k)=1}(\kappa_1*b)(dl)\chi(l)\Delta\bigg(\frac{l}{Dpk}\bigg)\\
=&\r_1^* Dpk\sum_{d_1}\,\sum_{d_2}\frac{\tg_3(d_1d_2k)}{d_1d_2}\tk_{1}(d_1;d_2k)\la_1(d_1d_2k)\sum_{(l_2,k)=1}\frac{b(d_2l_2)\chi(l_2)}{l_2}\\
&+O\big(\al^{30}Dpk\big)
\endaligned
$$
Inserting this into (15.7) and rewriting $d$, $l$ and $m$ for $d_2$, $l_2$ and $d_1$ respectively, we obtain
$$
\Phi_{1}(p)=\frac{\r_1^* Dp}{\varphi(D)}\sum_{d}\sum_{l}\frac{b(dl)\chi(l)}{dl}\d_1(d,l)+o(P)\eqno(15.11)
$$
where
$$
\d_1(d,l)=\sum_{(k, l)=1}\frac{\mu\chi(k)}{\varphi(k)}\sum_{m}\frac{\tk_{1}(m; dk)\la_1(mdk)\tg_3(mdk)}{m}.
$$
On substituting $n=mk$ we can write
$$
\d_1(d,l)=\sum_n\frac{\la_1(dn)\tg_3(dn)}{n}\xi_1(n;d,l)\eqno(15.12)
$$
with
$$
\xi_1(n;d,l)=\sum_{\substack{n=mk\\(k,l)=1}}\frac{\mu\chi(k)k}{\varphi(k)}\tk_{1}(m;dk).\eqno(15.13)
$$
It can be verified, for given $d$ and $l$,  that $\xi_1(n;d,l)$ is a multiplicative function of $n$.  On the other hand, we have
$$
\la_1(dn)=\la_1(d)\tilde{\la}_1(n,d)
$$
with
$$
\tilde{\la}_1(n,d)=\prod_{\substack{q|n\\(q,d)=1}}\la_1(q).
$$
Hence
$$
\d_1(d,l)=\la_1(d)\sum_n\frac{\tilde{\la}_1(n,d)\xi_1(n;d,l)\tg_3(dn)}{n}.
$$
On the right side above, we can replace the factor  $\tg_3(dn)$ by $(dn)^{\be_3}g(P_4/(dn))$ with a negligible error. Since
$$
y^{\be_3}g\bigg(\frac{P_4}{y}\bigg)=\tp\int_{(1)}\frac{P_4^{s+\be_3}}{y^s}\frac{\omega_1(s+\be_3)}{s+\be_3}\,ds,
$$
it follows that
$$
\d_1(d,l)=\la_1(d)\cdot\tp\int_{(1)}\bigg(\sum_n\frac{\tilde{\la}_1(n,d)\xi_1(n;d,l)}{n^{1+s}}\bigg)\frac{P_4^{s+\be_3}}{d^s}\,\frac{\omega_1(s+\be_3)}{s+\be_3}\,ds+O(\ep).\eqno(15.14)
$$

Note that $\tilde{\la}_1(q^r,d)=\tilde{\la}_1(q,d)$ for any $r$. If $\sigma>1$, then
$$
\sum_n\frac{\tilde{\la}_1(n,d)\xi_1(n;d,l)}{n^{s}}=\prod_{q}\bigg(1+\tilde{\la}_1(q,d)\sum_r\frac{\xi_1(q^r;d,l)}{q^{rs}}\bigg).
$$
If $(q,dl)=1$, then 
$$
\xi_1(q;d,l)=\sum_{h\in\n(q)}\frac{\kappa_1(qh)\chi(h)}{h}-\frac{\chi(q)q}{q-1}=q^{-\be_1}+q^{-\be_2}-1-\chi(q)+O(1/q),
$$
so that
$$
\frac{(1-q^{-s-\be_1})(1-q^{-s-\be_2})}{(1-q^{-s})(1-\chi(q)q^{-s})}\bigg(1+\tilde{\la}_1(q,d)\sum_r\frac{\xi_1(q^r;d,l)}{q^{rs}}\bigg)
=1+O(q^{-19/10})
$$
for $\sigma>9/10$. In case $(q, dl)>1$ and $\sigma>9/10$, the left side above is trivially $$1+O(q^{-9/10}).$$
It follows that  the function
$$
\m_1(d,l;s):=\frac{\zeta(s)L(s,\chi)}{\zeta(s+\be_1)\zeta(s+\be_2)}\sum_n\frac{\tilde{\la}_1(n,d)\xi_1(n;d,l)}{n^{s}}
$$
is  analytic  and it satisfies
$$
\m_1(d,l;s)\ll \prod_{q|dl}\bigg(1+\frac{c}{q^{9/10}}\bigg)
$$
for $\sigma>9/10$. The right side of  (15.14) can be rewritten as
$$\aligned
\la_1(d)\cdot\tp\int_{(1)}\frac{\zeta(1+s+\be_1)\zeta(1+s+\be_2)\m_1(d,l;1+s)}{\zeta(1+s)L(1+s,\chi)}\frac{P_4^{s+\be_3}}{d^s}\,\frac{\omega_1(s+\be_3)}{s+\be_3}\,ds+O(\ep).
\endaligned
$$
Assume $dl<PT^{-2}$, and  $(dl,D)=1$. Note that $P_4/d>T$. In a way similar to the proof of Lemma 8.4,  we deduce that
$$
\d_1(d,l)=\la_1(d)\sum_{j\le 3}\r_{1j}d^{\be_j}\m_1(d,l;1-\be_j)+O(\ep_1)\eqno(15.15)
$$
where $\r_{1j}$, $1\le j\le 3$, is the residue of the function
$$
 \frac{\zeta(1+s+\be_1)\zeta(1+s+\be_2)}{\zeta(1+s)L(1+s,\chi)}\,\frac{P_4^{s+\be_3}\omega_1(s+\be_3)}{s+\be_3}\eqno(15.16)
$$
at $s=-\be_j$ (the function (15.16) also has a simple pole at $s=\tilde{\rho}-1$, while the residue at this point can be regarded as an acceptable error). Inserting this into (15.11) and substituting $n=dl$ we obtain
$$
\Phi_{1}(p)=\frac{\r_1^* Dp}{\varphi(D)}\sum_{j\le 3}\r_{1j}\s_{1j}+o(p)\eqno(15.17)
$$
where
$$
\s_{1j}=\sum_{n}\frac{b(n)}{n}\sum_{n=dl}\la_1(d)d^{\be_j}\chi(l)\m_1(d,l;1-\be_j).
$$

If $(q,dl)=1$, then $\tilde{\la}_1(q,d)=\tilde{\la}_1(q,1)$ and $\xi_1(q^r;d,l)=\xi_1(q^r;1,1)$ for any $r$. Hence
$$
\frac{\m_1(d,l;s)}{\m_1(1,1;s)}=\prod_{q|dl}\bigg(1+\tilde{\la}_1(q,d)\sum_r\frac{\xi_1(q^r;d,l)}{q^{rs}}\bigg)
\bigg(1+\tilde{\la}_1(q,1)\sum_r\frac{\xi_1(q^r;1,1)}{q^{rs}}\bigg)^{-1}.\eqno(15.18)
$$
We can rewrite
$$
\s_{1j}=\m_1(1,1;1-\be_j)\sum_{n}\frac{b(n)\varpi_{1j}(n)}{n}\eqno(15.19)
$$
with
$$
\varpi_{1j}(n)=\sum_{n=dl}\la_1(d)d^{\be_j}\chi(l)\frac{\m_1(d,l;1-\be_j)}{\m_1(1,1;1-\be_j)}.
$$
In view of (15.18), we see that $\varpi_{1j}(n)$ is a multiplicative function.

Let
$$
\q=\prod_{q< D^4}q.
$$
Every $n$ can be uniquely written as $n=n_1n_2$ with $n_1\in\n(\q)$ and $(n_2,\q)=1$. Hence
$$
\sum_{n}\frac{b(n)\varpi_{1j}(n)}{n}=\sum_{n_1\in\n(\q)}\,\frac{\varpi_{1j}(n_1)}{n_1}\sum_{(n,\q)=1}\frac{b(n_1n)\varpi_{1j}(n)}{n}.\eqno(15.20)
$$
Using the the Rankin trick (see [14, Section 13.2], for example), we can impose the constraint $n_1<T$  to the right side with an acceptable error $O(\ep_1)$. In addition, by (15.18), for $n<PT^{-3}$, $(n,\q)=1$ and $dl=n$, 
$$
 \la_1(d)\frac{\m_1(d,l;s)}{\m_1(1,1;s)}=1+O(D^{-c}),
 $$
so that
$$
\varpi_{1j}(n)=\chi(n)\varrho_j^*(n)+O(\tau_2(n)D^{-c})\eqno(15.21)
$$
where
$$
\varrho_j^*(n)=\sum_{d|n}d^{\be_j}\chi(d).
$$
The following lemma will be proved in Appendix B.
\medskip\noindent

{\bf Lemma 15.1.}\,\,{\it Suppose $n_1\in\n(\q)$ and $n_1<T$. Then
$$
\sum_{(n,\q)=1}\frac{b(n_1n)\chi(n)\varrho_j^*(n)}{n}=\chi(n_1)\tau_2(n_1)\mathfrak{e}_j+O\big(\al_1\tau_2(n_1)\big)
$$
where
$$
\mathfrak{e}_j=(e_{1j}+\iota_2 e_{2j})(\bar{\iota}_3e_{3j}+\bar{\iota}_4e_{2j})
$$
with
$$
e_{2j}=\bigg(1-\frac{2j}{5}+\frac{8j}{25\pi i}\bigg)\exp\{5\pi i/4\}-\frac{8j}{25\pi i},
$$
$$
e_{3j}=\bigg(1-\frac{2j}{3}+\frac{j}{1.1205\pi i}\bigg)\exp\{0.747\pi i\}-\frac{j}{1.1205\pi i}
$$
and}
$$
e_{1j}=e_{1j}'-e_{1j}'',
$$
$$ e_{1j}'=\bigg(1-\frac{2j}{3}+\frac{j}{1.134\pi i}\bigg)\exp\{0.756\pi i\}-\frac{j}{1.134\pi i},
$$
$$
e_{1j}''=\frac{j}{0.756}\int_{0}^{0.004}\big(
\exp\{(3/2)(0.504-z)\pi i\}-\exp\{(3/4)\pi i\}\big)\,dz.
$$

\medskip\noindent

By (15.19)-(15.21) and Lemma 15.1 we obtain
$$
\s_{1j}=\mathfrak{e}_j\m_1(1,1;1-\be_j)\sum_{\substack{n\in\n(\q)\\n<T}}\,\frac{\chi(n)\tau_2(n)\varpi_{1j}(n)}{n}+O(1/\l).\eqno(15.22)
$$
This remains valid if the constraint $n\in\n(\q)$ on the right side is removed, since,  for $D^4<q<T$,
$$
\varrho_j^*(q)=\varrho_j(q)+O(\nu(q))\ll\al_1+\nu(q).
$$
This yields, by (15.21),
$$
\sum_{\substack{n\notin\n(\q)\\n<T}}\,\frac{\chi(n)\tau_2(n)\varpi_{1j}(n)}{n}\ll 1/\l.
$$

To apply (15.22) we need two lemmas which will be proved in Appendix A.

\medskip\noindent

{\bf Lemma 15.2.}\,\, {\it If $|s-1|<5\al$, then}
$$
\m_1(1,1;s)=\prod_{(q,D)=1}\frac{1-\chi(q)q^{-2}}{1-q^{-2}}+O(\al_1).
$$

\medskip\noindent

{\bf Lemma 15.3.}\,\,{\it For $\sigma\ge 9/10$ the function
$$
\u_{1j}(s):=\frac{1}{\zeta(s)^2L(s,\chi)^2}\sum_{n}\frac{\chi(n)\tau_2(n)\varpi_{1j}(n)}{n^s}
$$
is analytic and bounded. Further we have}
$$
\u_{2j}(1)=\frac{\varphi(D)^2}{D^2}\prod_{(q,D)=1}\frac{(1-q^{-2})^2}{1-\chi(q)q^{-2}}+O(\al_1).
$$

By (4.2) and (4.3),
$$\aligned
\sum_{n<T}\frac{\chi(n)\tau_2(n)\varpi_{1j}(n)}{n}=&\sum_{n}\frac{\chi(n)\tau_2(n)\varpi_{1j}(n)}{n}g\bigg(\frac{T}{n}\bigg)+O(\al_1)\\
&=\tp\int_{(1)}\zeta(1+s)^2L(1+s,\chi)^2\u_{2j}(1+s)\frac{T^s\omega_1(s)\,ds}{s}+O(\al_1).
\endaligned
$$
By Lemma 15.3, we can move the contour of integration in the same way as in the proof of Lemma 8.4 to obtain
$$
\sum_{n<T}\frac{\chi(n)\tau_2(n)\varpi_{1j}(n)}{n}=L'(1,\chi)^2\u_{2j}(1)+O(\al_1).
$$
This together with Lemma 15.2 and 15.3 yields
$$
\m_1(1,1;1-\be_j)\sum_{\substack{n\in\n(\q)\\n<T}}\,\frac{\chi(n)\tau_2(n)\varpi_{1j}(n)}{n}=\aa \frac{\varphi(D)}{D}+O(1/\l^3),
$$
since
$$
\frac{\varphi(D)}{D}\prod_{(q,D)=1}(1-q^{-2})=\frac{6}{\pi^2}\prod_{q|D}\frac{q}{q+1}.
$$
It follows by (15.22) that
$$
\s_{1j}=\mathfrak{e}_j\aa\frac{\varphi(D)}{D}+O(1/\l^3)\eqno(15.23)
$$

 By Lemma 5.8,
$$
\r_1^*=\be_1\be_2L'(1,\chi)+O(1/\l^{24})
$$
and
$$
\r_{1j}=\frac{P_4^{\be_3-\be_j}}{(\be_{j+1}-\be_j)(\be_{j+2}-\be_j)L'(1,\chi)}+O(1/\l^3).
$$
Hence, by direct calculation,
$$
\r_1^*\r_{11}=1+O(1/\l),
$$
$$
\r_1^*\r_{12}=2+O(1/\l),
$$
$$
\r_1^*\r_{13}=1+O(1/\l).
$$
Combining these relations with  (15.23) , (15.17) and (15.6) we conclude
$$
\Phi_1=\big(\mathfrak{e}_1+2\mathfrak{e}_2+\mathfrak{e}_3\big)\aa\p+o(\p)\eqno(15.24)
$$

 \section  *{\centerline{16. Evaluation of  $\Phi_2$}}
\medskip\noindent 

Recall that $\Phi_2$ is given by (13.9). Write
$$
\k_2\sp=Z(s,\pc)^{-1}\frac{L(s+\be_1,\psi)}{L\sp}B\sp N(s+\be_3,\psi)K(1-s-\be_2,\bp).
$$
Similar to (15.3),
$$
\Phi_2=\sum_{\psi\in\Psi_1}(p_\psi t_0)^{\be_1}\big(I_3^+(\psi)-I_3^-(\psi)\big)+O(\ep)
$$
where
$$
I_3^\pm(\psi)=\tp\int_{\j(\pm\al)}\k_2\sp\omega(s)\,ds.
$$
Note that the length of the sum $B\sp N(s+\be_3,\psi)$ is $\le PT^{-1}$. For $\psi\in\Psi_1$, moving the segment $\j(-\al)$ to 
$\j(-\l^9)$ we obtain, by simple estimation,
$$
I_3^-(\psi)\ll\ep.
$$
Hence
$$
\Phi_2=\sum_{\psi\in\Psi_1}(p_\psi t_0)^{\be_1}I_3^+(\psi)+O(\ep).
$$
In a way similar to the proof of (15.6), we can move the segment $\j(\al)$ to $\j(1)$ and then extend the sum over $\Psi_1$ to that over $\Psi$ with an acceptable error. Hence
$$
\Phi_2=\sum_{p\sim P}(pt_0)^{\be_1}\,\sideset{}{^*}\sum_{\psi(\bmod\,p)}\frac{1}{2\pi i} \int_{\j(1)}\k_2\sp\omega(s)\,ds+o(\p).
\eqno(16.1)
$$
Let $\kappa_2(n)$ and $b_1(n)$ be given by
$$
\sum_n\frac{\kappa_2(n)}{n^s}=\frac{\zeta(s+\be_1)}{\zeta(s)},\qquad \sigma>1,
$$
$$
\sum_n\frac{b_1(n)\psi(n)}{n^s}=B\sp N(s+\be_3,\psi).
$$
Note that
$$
K(1-s-\be_2,\bp)=\sum_n\frac{\tg_2(n)\bp(n)}{n^{1-s}}
$$
with
 $$
\tg_2(y)=y^{\be_2}g^*(P_4/y).
$$
We can rewrite (16.1) as
$$
\Phi_2=\Theta_2(\be_1, \bold{k}_2^*, \bold{a}_2^*)+o(p)
$$
with $\kappa^*_2=\kappa_2*b_1$ and $a^*_2=\tg_2$. It follows by Proposition 14.1 that
$$
\Phi_2=\sum_{p\sim P}(pt_0)^{\be_1}\Phi_2(p)+o(\p)\eqno(16.2)
$$
where
$$
\Phi_2(p)=
\frac{1}{\varphi(D)}\sum_k\frac{\mu\chi(k)}{k\varphi(k)}\sum_d\frac{\tg_2(dk)}{d}\,\sum_{(l,k)=1}(\kappa_2*b_1)(dl)\chi(l)\Delta\bigg(\frac{l}{Dpk}\bigg).
\eqno(16.3)
$$

Similar to (7.17),
$$\aligned
\sum_{(l,k)=1}(\kappa_2*b_1)&(dl)\chi(l)\Delta\bigg(\frac{l}{Dpk}\bigg)\\
=&\sum_{d=d_1d_2}\,\sum_{(l_2,k)=1}b_1(d_2l_2)\chi(l_2)\sum_{(l_1,d_2k)=1}\kappa_2(d_1l_1)\chi(l_1)\Delta\bigg(\frac{l_1l_2}{Dpk}\bigg).
\endaligned
$$
 The innermost sum is, by the Mellin transform, equal to
$$
\tp\int_{(2)}\bigg(\sum_{(l_1,d_2k)=1}\frac{\kappa_2(d_1l_1)\chi(l_1)}{l_1^s}\bigg)\bigg(\frac{Dpk}{l_2}\bigg)^s\delta(s)\,ds.
\eqno(16.4)$$
For $\sigma>1$,
$$\aligned
\sum_{(l_1,d_2k)=1}\frac{\kappa_2(d_1l_1)\chi(l_1)}{l_1^s}=\tk_{2}(d_1;d_2k,s)\la_2(d_1d_2k,s)\frac{L(s+\be_1,\chi)}{L(s,\chi)}
\endaligned
$$
where
$$
\tk_{2}(d_1;r,s)=\sum_{\substack{h\in\n(d_1)\\(h, r)=1}}\frac{\kappa_2(d_1h)\chi(h)}{h^s}
$$
and
$$
\la_2(n,s)=\prod_{q|n}\frac{1-\chi(q)q^{-s-\be_1}}{1-\chi(q)q^{-s}}.
$$
 For notational simplicity we shall write $\tk_2(m;r)$ and $\la_2(n)$  for  $\tk_2(m;r,1)$ and $\la_2(n,1)$ respectively. Note that $Dpk/l_2>T$ if $l_2<P_2^2$. In a way similar to the treatment of (7.19), 
we see that  the  expression (16.4) is equal to the residue of the integrand at $s=\tilde{\rho}$ plus an acceptable error. Further, by (5.15), we can replace $\tilde{\rho}$ by $1$ in the expression for this residue, with an acceptable error. Thus we have
$$
\aligned
\sum_{(l_1,d_2k)=1}\kappa_2(d_1l_1)\chi(l_1)\Delta\bigg(\frac{l_1l_2}{Dpk}\bigg)=\r_2^*\bigg(\frac{Dpk}{l_2}\bigg)\tk_{2}(d_1;d_2k)\la_2(d_1d_2k)+O\bigg(\frac{\al^{100}\tau_2(d_1)Dpk}{l_2}\bigg)
\endaligned
$$
where
$$
\r_2^*=\frac{L(1+\be_1,\chi)}{L'(1,\chi)}\delta(1).
$$
This yields
$$\aligned
\sum_{(l,k)=1}(\kappa_2*b_1)(dl)\chi(l)\Delta\bigg(\frac{l}{Dpk}\bigg)=&\r_2^* Dpk
\sum_{d=d_1d_2}\tk_{2}(d_1;d_2k)\la_2(dk)
\sum_{(l_2,k)=1}\frac{b_1(d_2l_2)\chi(l_2)}{l_2}\\
&+O\big(\al^{50}\tau_3(d)Dpk\big).
\endaligned
$$
Hence
$$\aligned
\sum_d\frac{\tg_2(dk)}{d}&\sum_{(l,k)=1}(\kappa_2*b_1)(dl)\chi(l)\Delta\bigg(\frac{l}{Dpk}\bigg)\\
=&\r_2^* Dpk\sum_{d_1}\,\sum_{d_2}\frac{\tg_2(d_1d_2k)}{d_1d_2}\tk_{2}(d_1;d_2k)\la_2(d_1d_2k)\sum_{(l_2,k)=1}\frac{b_1(d_2l_2)\chi(l_2)}{l_2}\\
&+O(\al^{20}Dpk).
\endaligned
$$
Inserting this into (16.3) and rewriting $d$, $l$ and $m$ for $d_2$, $l_2$ and $d_1$ respectively, we obtain
$$
\Phi_{2}(p)=\frac{\r_2^* Dp}{\varphi(D)}\sum_{d}\sum_{l}\frac{b_1(dl)\chi(l)}{dl}\d_2(d,l)+o(1)\eqno(16.5)
$$
where
$$
\d_2(d,l)=\sum_{(k, l)=1}\frac{\mu\chi(k)}{\varphi(k)}\sum_{m}\frac{\tk_{2}(m; dk)\la_2(mdk)\tg_2(mdk)}{m}.
$$
On substituting $n=mk$ we can write
$$
\d_2(d,l)=\sum_n\frac{\la_2(dn)\tg_2(dn)}{n}\xi_2(n;d,l)\eqno(16.6)
$$
with
$$
\xi_2(n;d,l)=\sum_{\substack{n=mk\\(k,l)=1}}\frac{\mu\chi(k)k}{\varphi(k)}\tk_{2}(m;dk).\eqno(16.7)
$$
It can be verified, for given $d$ and $l$,  that $\xi_2(n;d,l)$ is a multiplicative function of $n$, and
$$
\d_2(d,l)=\la_2(d)\sum_n\frac{\tilde{\la}_2(n,d)\xi_2(n;d,l)\tg_2(dn)}{n}
$$
with
$$
\tilde{\la}_2(n,d)=\prod_{\substack{q|n\\(q,d)=1}}\la_2(q).\eqno(16.8)
$$
Hence, similar to (15.14),
$$
\d_2(d,l)=\la_2(d)\cdot\tp\int_{(1)}\bigg(\sum_n\frac{\tilde{\la}_2(n,d)\xi_2(n;d,l)}{n^{1+s}}\bigg)\frac{P_4^{s+\be_2}}{d^s}\,\frac{\omega_1(s+\be_2)}{s+\be_2}\,ds+O(\ep).\eqno(16.9)
$$
It can be verified, for $\sigma>9/10$,  that the function
$$
\m_2(d,l;s):=\frac{\zeta(s)L(s,\chi)}{\zeta(s+\be_1)}\sum_n\frac{\tilde{\la}_2(n,d)\xi_2(n;d,l)}{n^{s}}
$$
is  analytic  and it satisfies
$$
\m_2(d,l;s)\ll \prod_{q|dl}\bigg(1+\frac{c}{q^{9/10}}\bigg).
$$
 Thus we can rewrite $\d_2(d,l)$ as
$$\aligned
\la_2(d)\cdot\tp\int_{(1)}&\frac{\zeta(1+s+\be_1)\m_2(d,l;1+s)}{\zeta(1+s)L(1+s,\chi)}\frac{P_4^{s+\be_2}}{d^s}\,\frac{\omega_1(s+\be_2)}{s+\be_2}\,ds+O(\ep).
\endaligned
$$
Assume $dl<P_2^2$, and  $(dl,D)=1$. Note that $P_4/d>T$. In a way similar to the proof of (15.15),  we deduce that
$$
\d_2(d,l)=\la_2(d)\sum_{j=1,2}\r_{2j}d^{\be_j}\m_2(d,l;1-\be_j)+O(\ep_1)\eqno(16.10)
$$
where $\r_{2j}$, $ j=1,2$, is the residue of the function
$$
 \frac{\zeta(1+s+\be_1)}{\zeta(1+s)L(1+s,\chi)}\,\frac{P_4^{s+\be_2}\omega_1(s+\be_2)}{s+\be_2}\eqno(16.11)
$$
at $s=-\be_j$. Inserting this into (16,5) and substituting $n=dl$ we obtain
$$
\Phi_{2}(p)=\frac{\r_2^* Dp}{\varphi(D)}\sum_{j=1,2}\r_{2j}\s_{2j}+o(p)\eqno(16.12)
$$
where
$$
\s_{2j}=\sum_{n}\frac{b_1(n)}{n}\sum_{n=dl}\la_2(d)d^{\be_j}\chi(l)\m_2(d,l;1-\be_j).
$$

Let $\m_2^*(s)$ be given by
$$
\m_2^*(s)=\m_2(1,1;s)\quad\text{if}\quad \chi(2)\ne 1,
$$
$$
\m_2^*(s)=2\prod_{q>2}\frac{1-q^{-s-\be_1}}{(1-q^{-s})(1-\chi(q)q^{-s})}\bigg(1+\tilde{\la}_2(q,1)\sum_r\frac{\xi_2(q^r;1,1)}{q^{rs}}\bigg)\quad\text{if}\quad \chi(2)=1.
$$

The following result will be proved in Appendix A.

\medskip\noindent

{\bf Lemma 16.1.}\,\, {\it Suppose $dl<PT^{-2}$ and $|s-1|<5\al$. Then
$$
\m_2^*(s)=\mathfrak{p}+O(1/\l^8)
$$
where
$$
\mathfrak{p}=\prod_{q}\frac{1}{1-\chi(q)q^{-1}}\bigg(1-\frac{\chi(q)}{q-1}\bigg)\quad\text{if}\quad \chi(2)\ne 1,
$$
and}
$$
\mathfrak{p}=2\prod_{q>2}\frac{1}{1-\chi(q)q^{-1}}\bigg(1-\frac{\chi(q)}{q-1}\bigg)\quad\text{if}\quad \chi(2)=1.
$$

As a direct consequence of Lemma 16.1 we have $|\m_2^*(s)|\gg 1$ if $|s-1|<5\al$. Thus we can write
$$
\s_{2j}=\m_2^*(1-\be_j)\sum_{n}\frac{b_1(n)\varpi_{2j}(n)}{n}\eqno(16.13)
$$
with 
$$
\varpi_{2j}(n)=\sum_{n=dl}\la_2(d)d^{\be_j}\chi(l)\frac{\m_2(d,l;1-\be_j)}{\m_2^*(1-\be_j)}.
$$

Note that $\varpi_{2j}(n)$ is multiplicative. In a way similar to the proof of (15.20) we have
$$
\sum_{n}\frac{b_1(n)\varpi_{2j}(n)}{n}=\sum_{\substack{n_1\in\n(\q)\\n_1<T}}\frac{\varpi_{2j}(n_1)}{n_1}\sum_{(n,\q)=1}\frac
{b_1(n_1n)\varpi_{2j}(n)}{n}+O(D^{-c}).
$$
Assume $n_1\in\n(\q)$ and  $n_1<T$.  For  $(n,\q)=1$ we have
$$
b_1(n_1n)=\sum_{n_1=l_1m_1}\,\sum_{n=lm}(l_1l)^{-\be_3}g^*\bigg(\frac{T^2}{l_1l}\bigg)b(m_1m).
$$
Hence
$$
\sum_{(n,\q)=1}\frac{b_1(n_1n)\varpi_{2j}(n)}{n}=\sum_{n_1=l_1m_1}\sum_{(l,\q)=1}\frac{\varpi_{2j}(l)}{l(l_1l)^{\be_3}}g^*\bigg(\frac{T^2}{l_1l}\bigg)\sum_{(m,\q)=1}\frac{b(m_1m)\varpi_{2j}(m)}{m}+O(D^{-c}),
$$
since the terms with $(l,m)>1$ above contribute $\ll D^{-c}$. Further, for $(m,\q)=1$ and $m<P$ we have trivially
$$
\varpi_{2j}(m)=\chi(m)\varrho^*_j(m)+O\big(\tau_2(m)D^{-c}\big).
$$
Hence, for $m_1<T$,
$$
\sum_{(m,\q)=1}\frac{b(m_1m)\varpi_{2j}(m)}{m}=\mathfrak{e}_j\chi(m_1)\tau_2(m_1)+O(D^{-c})
$$
by Lemma 15.1. It follows that
$$
\sum_{(n,\q)=1}\frac{b_1(n_1n)\varpi_{2j}(n)}{n}=\mathfrak{e}_j\sum_{n_1=l_1m_1}\chi(m_1)\tau_2(m_1)\sum_{(l,\q)=1}\frac{\varpi_{2j}(l)}{l(l_1l)^{\be_3}}g^*\bigg(\frac{T^2}{l_1l}\bigg)+O(D^{-c}).\eqno(16.14)
$$
On the other hand, if $1<l<T^5$ and $(l,\q)=1$, then for any $q|l$,
$$
\varrho_{2j}(l)\ll\al_1+O(\nu(q)).
$$
Thus the innermost sum in (16.14) is equal to $1+O(1/\l^{7})$. It follows that
$$
\sum_{(n,\q)=1}\frac{b_1(n_1n)\varpi_{2j}(n)}{n}=\mathfrak{e}_j\sum_{m_1|n_1}\chi(m_1)\tau_2(m_1)+O(1/\l^4).
$$
Since $1*\chi\tau_2=\nu*\chi$ (here $1$ denotes the arithmetic function identically equal to $1$), it follows that
$$
\sum_{n}\frac{b_1(n)\varpi_{2j}(n)}{n}=\mathfrak{e}_j\sum_{\substack{n_1\in\n(\q)\\n_1<T}}\frac{\varpi_{2j}(n_1)(\nu*\chi)(n_1)}{n_1}+O(1/\l).\eqno(16.15)
$$
On the right side above, the constraint $n_1\in\n(\q)$ can be removed with an acceptable error.

The following lemma will be proved in Appendix A.
\medskip\noindent

{\bf Lemma 16.2.}\,\, {\it The function
$$
\u_{2j}(s)=\frac{1}{\zeta(s)^3L(s,\chi)^3}\sum_n\frac{\varpi_{2j}(n)(\nu*\chi)(n)}{n^s}
$$
is analytic and bounded for $\sigma>9/10$. Further we have}
$$
\u_{2j}(1)=\frac{6}{\pi^2}\frac{\varphi(D)}{D\mathfrak{p}}\prod_{q|D}\frac{q}{q+1}+O(1/\l^4).
$$
\medskip\noindent

We have
$$\aligned
\sum_{n<T}\frac{\varpi_{2j}(n)(\nu*\chi)(n)}{n}=&\sum_{n}\frac{\varpi_{2j}(n)(\nu*\chi)(n)}{n}g\bigg(\frac{T}{n}\bigg)+O(1/\l^{10})\\
=&\tp\int_{(1)}\zeta(1+s)^3L(1+s,\chi)^3\u_{2j}(1+s)T^s\omega_1(s)\,\frac{ds}{s}+O(1/\l^{10}).
\endaligned
$$
The contour of integration is moved in the same way as in the proof of Lemma 8.4. Thus the right side above is, by Lemma 16.2,  equal to
$$
L'(1,\chi)^3\u_{2j}(1)+O(1/\l^4)=\frac{\aa}{\mathfrak{p}}\frac{\varphi(D)}{D}L'(1,\chi)+O(1/\l^4).
$$
Hence, by (16.15),
$$
\sum_{n}\frac{b_1(n)\varpi_{2j}(n)}{n}=\frac{\aa\mathfrak{e}_j}{\mathfrak{p}}\frac{\varphi(D)}{D}L'(1,\chi)+O(1/\l^4).
$$
Inserting this into (16.13) and applying Lemma 16.1 we obtain
$$
\s_{2j}=\aa\mathfrak{e}_j\frac{\varphi(D)}{D}L'(1,\chi)+O(1/\l^4).\eqno(16.16)
$$
On the other hand, by Lemma 5.8 and direct calculation,
$$
\r_2^*=\be_1+O(1/\l^{10}),\quad \r_{2j}=-\frac{1}{\be_1L'(1,\chi)}+O(\l^6),\quad j=1,2,
$$
so that
$$
\r_2^*\r_{2j}=-\frac{1}{L'(1,\chi)}+O(\l^{-1}L'(1,\chi)^{-1}).
$$
This together with (16.16) and (16.12) yields
$$
\Phi_2(p)=-(\mathfrak{e}_1+\mathfrak{e}_2)\aa p+o(p).
$$
Since $(pt_0)^{\be_1}=-1+O(\al_1)$, it follows by (16.2) that
$$
\Phi_2=(\mathfrak{e}_1+\mathfrak{e}_2)\aa \p+o(\p).\eqno(16.17)
$$

 \section  *{\centerline{17. Evaluation of $\Phi_3$}}
\medskip\noindent

Recall that $\Phi_3$ is given by (13.10). Write
$$
\k_3\sp=\frac{L(s+\be_1,\psi)}{L\sp}B\sp G\sp N(s+\be_2,\psi) N(s+\be_3,\psi)F(1-s,\bp)
$$
and
$$
I_4^\pm(\psi)=\tp\int_{\j(\pm\al)}\k_3\sp\omega(s)\,ds.
$$
Similar to the treatment of $\Phi_2$, we have
$$
\Phi_3=\sum_{\psi\in\Psi_1}(p_\psi t_0)^{\be_3}\big(I_4^+(\psi)-I_4^-(\psi)\big)+O(\ep).\eqno(17.1)
$$

To treat the sum of  $I_4^+(\psi)$ we move the segment $\j(\al)$ to $\j(1)$, and then extend the sum over $\Psi_1$ to the sum over $\Psi$ with an acceptable error. Hence
$$
\sum_{\psi\in\Psi_1}(p_\psi t_0)^{\be_3}I_4^+(\psi)=\sum_{p\sim P}(p t_0)^{\be_3}\Phi_3^+(p)+o(\p)\eqno(17.2)
$$
where
$$
\Phi_3^+(p)=\sideset{}{^*}\sum_{\psi\bmod\,p}\tp\int_{\j(1)}\k_3\sp\,\omega(s)\,ds.
$$
For $\sigma>1$ we can write
$$
\frac{L(s+\be_1,\psi)}{L\sp}B\sp G\sp N(s+\be_2,\psi) N(s+\be_3,\psi)=\sum_m\frac{\nu^*(m)\psi(m)}{m^s}.
$$
Thus, replacing the segment $\j(1)$ by the line $\sigma=3/2$ and integrating term by term give
$$
\Phi_3^+(p)=\sum_m\,\sum_{n<D^4}\frac{\nu^*(m)\nu(n)}{n}\bigg(\frac{n}{m}\bigg)^{s_0}
\bigg(\sideset{}{^*}\sum_{\psi\bmod\,p}\psi(m)\bp(n)\bigg)\exp\bigg\{-\l^2_2\log^2\frac{n}{m}\bigg\}+O(\ep).
$$
By trivial estimation, the contribution from the terms with $m\ne n$ above is $o(p)$. Hence
$$
\Phi_3^+(p)=p\sum_{n<D^4}\frac{\nu^*(n)\nu(n)}{n}+o(p).\eqno(17.3)
$$
Note that for $n<D^4$
$$
\nu^*(n)=\sum_{n=n_1...n_6}\kappa_2(n_1)\chi(n_2n_3)\big(\vk_1(n_2)+\iota_2\vk_2(n_2)\big)
\big(\bar{\iota}_3\vk_3(n_3)+\bar{\iota_4}\vk_4(n_3)\big)\upsilon(n_4)\tilde{g}_2(n_5)\tilde{g}_3(n_6).
$$
A detailed analysis shows that (17.2) remains valid if, in the expression for $\nu^*(n)$, the factor $\kappa_2(n_1)$is replaced by $(1*\mu)(n_1)$, the factor 
$$\chi(n_2n_3)\big(\vk_1(n_2)+\iota_2\vk_2(n_2)\big)
\big(\bar{\iota}_3\vk_3(n_3)+\bar{\iota_4}\vk_4(n_3)\big)
$$
is replaced by $\mathfrak{e}_0\chi(n_2n_3)$ with
$$
\mathfrak{e}_0=\big(\vk_1(1)+\iota_2\vk_2(1)\big)
\big(\bar{\iota}_3\vk_3(1)+\bar{\iota_4}\vk_4(1)\big),\eqno(17.4)
$$
and the factors $\tilde{g}_2(n_5)$ and $\tilde{g}_3(n_6)$ are replaced by $1$ respectively. Since
$$
1*\mu*\chi*\chi*\upsilon*1*1=\nu,
$$
it follows that
$$
\Phi_3^+(p)=\mathfrak{e}_0p\sum_{n<D^4}\frac{\nu(n)^2}{n}+o(p).\eqno(17.5)
$$
The following lemma will be proved in Appendix B.

\medskip\noindent

{\bf Lemma 17.1.}\,\,{\it We have}
$$
\sum_{n<D^4}\frac{\nu(n)^2}{n}=\mathfrak{a}+o(1)
$$

Since $(pt_0)^{\be_3}=-1+O(\al_1)$, it follows by (17.5), (17.2) and Lemma 17.1 that
$$
\sum_{\psi\in\Psi_1}(p_\psi t_0)^{\be_3}I_3^+(\psi)=-\mathfrak{e}_0\aa\p+o(\p).\eqno(17.6)
$$

 To treat the sum involving $I_4^-(\psi)$ on $\j(-\al)$ we first apply Lemma 5.1 to obtain
$$
\frac{L(s+\be_1,\psi)}{L\sp}=(pt_0)^{-\be_1}\frac{L(1-s-\be_1,\bp)}{L(1-s,\bp)}(1+O(\al^6))
$$
for $\psi\in\Psi_1$ and $s\in\j(-\al)$.
Thus, in a way similar to the proof of (15.4),
$$
\sum_{\psi\in\Psi_1}(p_\psi t_0)^{\be_3}I_3^-(\psi)=\sum_{\psi\in\Psi_1}(p_\psi t_0)^{\be_2}\cdot\tp\int_{\j(-\al)}\k_3^*(s,\psi)
\omega(s)\,ds+o(\p)
$$
where
$$
\k_3^*(\psi)=\frac{L(1-s-\be_1,\bp)}{L1-s,\bp)}F(1-s,\bp)
B\sp G\sp N(s+\be_2,\psi) N(s+\be_3,\psi).
$$
Moving the segment $\j(-\al)$ to $\j(-1)$ and then extend the sum over $\Psi_1$ to then sum over $\Psi$ we obtain
$$
\sum_{\psi\in\Psi_1}(p_\psi t_0)^{\be_3}I_3^-(\psi)=\sum_{p\sim P}(pt_0)^{\be_2}\Phi_3^-(p)+o(\p)\eqno(17.7)
$$
where
$$
\Phi_3^-(p)=\sideset{}{^*}\sum_{\psi\bmod\,p}\tp\int_{\j(-1)}\k_3^*\sp\,\omega(s)\,ds.
$$
Write
$$
G\sp N(s+\be_2,\psi) N(s+\be_3,\psi)=\sum_l\frac{\nu_1^*(l)\psi(l)}{l^s},
$$
so that
$$
B\sp G\sp N(s+\be_2,\psi) N(s+\be_3,\psi)=\sum_n\frac{(b*\nu_1^*)(n)\psi(n)}{n^s}.
$$
On the other hand, for $\sigma<0$,
$$
\frac{L(1-s-\be_1,\bp)}{L(1-s,\bp)}=\sum_m\frac{\bar{\kappa}_2(m)\bp(m)}{m^{1-s}}
$$
with $\bar{\kappa}_2(m)=\overline{\kappa_2(m)}$, so that
$$
\frac{L(1-s-\be_1,\bp)}{L(1-s,\bp)}F(1-s,\bp)=\sum_n\frac{\varrho^*(n)\bp(n)}{n^{1-s}}
 $$
 with
 $$
 \varrho^*(n)=\sum_{\substack{n=lm\\l<D^4}}\nu(l)\bar{\kappa}_2(m).
 $$
Note that the arithmetic function $(b*\nu_1^*)(n)$ is supported on $n<PT^{-2}$. In a way similar to the proof of (17.3) we deduce that
$$
\Phi_3^-(p)=p\sum_{n}\frac{(b*\nu_1^*)(n)\varrho^*(n)}{n}+o(p).\eqno(17.8)
$$
We have
$$
\sum_{n}\frac{(b*\nu_1^*)(n)\varrho^*(n)}{n}=\sum_{l<D^4}\frac{\nu(l)}{l}\sum_m\frac{(b*\nu_1^*)(lm)\bar{\kappa}_2(m)}{m}.
$$
The inner sum is equal to
$$
\sum_{l=l_1l_2}\sum_{m_1}\,\sum_{(m_2,l_1)=1}\frac{b(l_1m_1)\nu_1^*(l_2m_2)\bar{\kappa}_2(m_1m_2)}{m_1m_2}.
$$
Note that $\nu_1^*(n)$ is supported on $n<T^5$. On the right side above, we can drop the terms with $m_2>1$ or $(m_1,\q)>1$ with an acceptable error. Hence
$$
\sum_{n}\frac{(b*\nu_1^*)(n)\varrho^*(n)}{n}=\sum_{l<D^4}\frac{\nu(l)}{l}\sum_{l=l_1l_2}\sum_{(m_1,\q)=1}\frac{b(l_1m_1)\nu_1^*(l_2)\bar{\kappa}_2(m_1)}{m_1}+o(1).
$$
If $m_1$ is square-free, then
$$
\bar{\kappa}_2(m_1)=\mu(m_1)\varrho_1(m_1).
$$
Hence, by Lemma 15.1 (see (15.)),
$$
\sum_{(m_1,\q)=1}\frac{b(l_1m_1)\bar{\kappa}_2(m_1)}{m_1}=\mathfrak{e}_1\chi(l_1)\tau_2(l_1)+O
$$
It follows that
$$
\sum_{n}\frac{(b*\nu_1^*)(n)\varrho^*(n)}{n}=\mathfrak{e}_1\sum_{l<D^4}\frac{\nu(l)}{l}\sum_{l=l_1l_2}\chi(l_1)\tau_2(l_1)\nu_1^*(l_2)+o(1).
$$
In a way similar to the proof of (17.5), by Lemma 17.1, we find that the right side is equal to
$$
\mathfrak{e}_1\sum_{l<D^4}\frac{\nu(l)^2}{l}+o(1)=\mathfrak{e}_1\aa+o(1).
$$
Inserting this into (17.8) gives
$$
\Phi_3^-(p)=\mathfrak{e}_1\aa p+o(p).
$$
Since $(pt_0)^{\be_2}=1+O(\al_1)$, it follows by (17.7) that
$$
\sum_{\psi\in\Psi_1}(p_\psi t_0)^{\be_3}I_3^-(\psi)=\mathfrak{e}_1\aa\p+o(\p). \eqno(17.9)
$$

Finally, from (17.1), (17.6) and (17.9) we conclude
$$
\Phi_3=-(\mathfrak{e}_0+\mathfrak{e}_1)\aa\p+o(\p). \eqno(17.10)
$$

\medskip\noindent
\medskip\noindent

 \section  *{\centerline{18. Proof of Proposition 2.5}}
\medskip\noindent

By the discussion at the end of Section 2, it suffices to prove (2.32) and (2.33).
\medskip\noindent

{\it Proof of (2.32).}
\medskip\noindent

By (12.3), (12,17), (13.7), (15.24), (16.17) and (17.10),
$$
\Xi_{13}=\mathfrak{c}_3\aa\p\eqno(18.1)
$$
with
$$
\mathfrak{c}_3=-i(3\mathfrak{e}_1+3\mathfrak{e}_2+\mathfrak{e}_3+\mathfrak{e}_0)+e^*_1+2e^*_2
+\epsilon.
$$
In view of (15.), we can write
$$
\mathfrak{e}_j=\mathfrak{e}'_j-\mathfrak{e}''_j,\quad 1\le j\le 3
$$
with
$$
\mathfrak{e}_j'=(e_{1j}'+\iota_2 e_{2j})(\bar{\iota}_3e_{3j}+\bar{\iota}_4e_{2j})
$$
$$
\mathfrak{e}_j''=e_{1j}''(\bar{\iota}_3e_{3j}+\bar{\iota}_4e_{2j}).
$$
By calculation (there is a theoretical interpretation),
$$
i(3\mathfrak{e}''_1+3\mathfrak{e}''_2+\mathfrak{e}''_3)+e^*_1=\epsilon.
$$
Hence
$$
\mathfrak{c}_3=-i(3\mathfrak{e}'_1+3\mathfrak{e}'_2+\mathfrak{e}'_3+\mathfrak{e}_0)+2e^*_2
+2\epsilon.
$$
Direct calculation shows that
$$
\Re\{\mathfrak{c}_3\}<-6.9951 \eqno(18.2)
$$

It follows from (8.24), (9.8) and (18.2) that
$$
\mathfrak{c}_1+\mathfrak{c}_2+2\Re\{\mathfrak{c}_3\}<0.001.
$$
This with together (8.23), (9.7) and (18.1)  yields (2.32).
\medskip\noindent

{\it Proof of (2.33).}
\medskip\noindent

By Lemma 8.1,
$$
\sum_{\psi\in\Psi_1}\,\sum_{\rho\in\z(\psi)}|J_1\ss|^2\omega(\rho)=2\Re\{\Theta_1(\bold{a}_{23}, \bold{a}_{23}\}+o(\p).\eqno(18.3)
$$
We have
$$\aligned
S_{j}(\bold{a}_{23},\bold{a}_{23})=&\sum_{r}\,\sum_d\frac{|\chi(d)||\mu\chi(r)|\la_{0j}(dr)}{dr\varphi(r)}\bigg(\sum_m\frac{\chi(m)
\tilde{f}\big(\log(drm)/\log P\big)}{m^{1-\be_j}}\bigg)\\&\times\bigg(\sum_n\frac{\chi(n)\tilde{f}\big(\log(drn)/\log P\big)\xi_{0j}(n;d,r)}{n}\bigg).
\endaligned
$$

The right side is split into three sums according to
 $$dr<P^{0.5},\quad P^{0.5}\le dr<P^{0.502},\quad
P^{0.502}\le dr<P^{0.504}.
$$
By the discussion in Section 8 and 10, the sum over $dr<P^{0.5}$ contributes $o(\al)$; the sum over $P^{0.5}\le dr<P^{0.502}$
is equal to
$$
\bigg(\frac{500L'(1,\chi)}{\log P}\bigg)^2\sum_{P^{0.5}\le n<P^{0.502}}\frac{|\chi(n)|\la_{0j}(n)}{\varphi(n)}
\big(-1-\be_j\log(n/P^{0.5})\big)\big(-1+\y_{1j}(n)\big)+o(\al);
$$
the sum over $P^{0.502}\le dr<P^{0.504}$ is equal to
$$
\bigg(\frac{500L'(1,\chi)}{\log P}\bigg)^2\sum_{P^{0.502}\le n<P^{0.504}}\frac{|\chi(n)|\la_{0j}(n)}{\varphi(n)}
\big(1-\be_j\log(P^{0.504}/n)\big)\big(1+\y_{1j}(n)\big)+o(\al).
$$
The factors $\be_j\log(n/P^{0.5})$, $\be_j\log(P^{0.504}/n)$ and $\y_{1j}(n)$ in the above sums make minor contribution. If we disregard these factors, the major contribution to $S_{j}(\bold{a}_{23},\bold{a}_{23})$ will be
$$
\bigg(\frac{500L'(1,\chi)}{\log P}\bigg)^2\sum_{P^{0.5}\le n<P^{0.504}}\frac{|\chi(n)|\la_{0j}(n)}{\varphi(n)}=\frac{1000\aa}{\log P}
+o(\al).$$
Thus we have the crude bound
$$
\Re\{S_{j}(\bold{a}_{23},\bold{a}_{23})\}<\frac{1100\aa}{\log P},
$$
so that
$$\Re\bigg\{\frac{1}{2\al}S_{1}(\bold{a}_{23},\bold{a}_{23})+\frac{2}{\al}S_{2}(\bold{a}_{23},\bold{a}_{23})
+\frac{3}{2\al}S_{3}(\bold{a}_{23},\bold{a}_{23})\bigg\}<\frac{4400\aa}{\pi}.
$$
This yields (2.33) by Lemma 8.1 and Proposition 7.1. $\Box$

\medskip\noindent

 \section  *{\centerline{Appendix A. Some Euler Products}}
\medskip\noindent

This appendix is devoted to proving Lemma 8.3, 15.2, 15.3, 16.1 and 16.2. For notational simplicity we shall write 
$$u=q^{-1},\qquad v=\chi(q).
$$
\medskip\noindent

{\it Proof of Lemma 8.3.}\,\, Note that
$$
\la(km,s)=\la(k,s)\tilde{\la}(m,k;s)
$$
with
$$
\tilde{\la}(m,k;s)=\prod_{\substack{q|m\\(q,k)=1}}\frac{(1-q^{-s-\be_1})(1-q^{-s-\be_2})(1-q^{-s-\be_3})}{1-q^{-s}}.
$$
For given $d$ and $h$, the function $\tilde{\la}(m,h,d;1-\be_j)\xi_j(m;d,h)$ is multiplicative in $m$, and  $\tilde{\la}(q^r,k;s)=\tilde{\la}(q,k;s)$ for any $r$. Thus, for $\sigma>1$,
$$
 \u_{j}(d,h;s)=\la(dh,1-\be_j)\prod_q\t_j(d,1,s;q)
 $$
 with
 $$
 \t_j(d,h,s;q)=\frac{(1-\chi(q)q^{-s-\be_{j+1}})(1-\chi(q)q^{-s-\be_{j+2}})}{1-\chi(q)q^{-s}}\bigg(1+
\tilde{\la}(q,dh;1-\be_j)\sum_r\frac{\chi(q^r)\xi_j(q^r;d,h)}{q^{rs}}\bigg).
$$
Since
$$\aligned
\tilde{\la}(q;d,1-\be_j)\chi(q)\xi_j(q^r;d,h)&=\chi(q)(\kappa(q)-q^{-\be_j})+O(q^{-1})\\
&=\chi(q)(q^{-\be_{j+1}}+q^{-\be_{j+2}}-1)+
O(q^{-1}),
\endaligned
$$
we have $ \t_j(d,h,s;q)=1+O(q^{-9/5})$ for $\sigma>9/10$, so $\u_{j\mu}(d,s)$ is analytic in this region, and
$$
 \u_{j}(d,h;s)=\la(dh,1-\be_j)\prod_{q<D}\t_j(d,h_1,s;q)+O(D^{-c}).
 $$
Since
 $$
 \la(dh,1-\be_j)=\frac{\varphi(dh)^2}{(dh_1)^2}+O(\al_1)
 $$
 and $\t_j(d,h,s;q)=1$ if $q|D$, the proof is reduced to showing that
$$
\t_j(d,h,s;q)=1+O(\al\log q/q)\quad\text{if}\quad (q,dh_1)=1,\eqno(A.1)
$$
$$
\t_j(d,h,s;q)=\frac{1}{1-vu}+O(\al\log q/q)\quad\text{if}\quad q|h_1,\eqno(A.2)
$$
and
$$
\t_j(d,h,s;q)=\frac{1-u-vu}{(1-vu)(1-u)}+O(\al\log q/q)\quad\text{if}\quad q|d\quad (q,h_1)=1,\eqno(A.3)
$$
provided
$$
 |s-1|<5\al,\quad q<D\quad\text{and}\quad (q,D)=1,
$$
which are henceforth assumed. We discuss in three cases.
\medskip\noindent

{\it Case 1.}\,\, $(q,dh)=1$. 
\medskip\noindent

 We have
$$
\xi_j(q^r,dh)=\tk(q^r;dh,1-\be_j)-\frac{\kappa(q^{r-1})q^{1-\be_j}}{q-1}.
$$
because $\tk(q^{r-1};dq,1-\be_j)=\kappa(q^{r-1})$.
Note that
$$
\kappa(q^r)=\sum_{\mu_1+\mu_2+\mu_3=r}q^{\mu_1\be_1+\mu_2\be_2+\mu_3\be_3}-\sum_{\mu_1+\mu_2+\mu_3=r-1}q^{\mu_1\be_1+\mu_2\be_2+\mu_3\be_3}=\tau_2(q^r)(1+O(\al r\log q))
$$
where $\mu_1$, $\mu_2$ and $\mu_3$ run through non-negative integers. Hence
$$
\xi_j(q^r,dh)=\sum_{\eta=0}^{\infty}(r+\eta+1)u^{\eta}-\frac{r}{1-u}+O(\al r\log q/q)=\frac{1}{(1-u)^2}+O(\al r\log q/q).
$$
It follows that
$$\aligned
\sum_r\frac{\chi(q^r)\xi_j(q^r;d,h)}{q^{rs}}=\frac{1}{(1-u)^2}\bigg(\frac{1}{1-vu}-1\bigg)+O(\al \log q/q).
\endaligned
$$
This together with the relations
$$
\tilde{\la}(q,dh;1-\be_j)=(1-u)^2+O(\al \log q/q)
$$
and
$$
 \frac{(1-\chi(q)q^{-s-\be_{j+1}})(1-\chi(q)q^{-s-\be_{j+2}})}{1-\chi(q)q^{-s}}=1-vu+O(\al \log q/q),
$$
yields (A.1). 

\medskip\noindent

{\it Case 2.}\,\,$q|h$.
\medskip\noindent

We have 
$$
\xi_j(q^r,dh)=\kappa(q^r)=r+1+O(\al r\log q/q),
$$
so that
$$
\sum_r\frac{\chi(q^r)\xi_j(q^r;d,h)}{q^{rs}}=\frac{1}{(1-vu)^2}-1+O(\al \log q/q).
$$
This yields (A.2) since $\tilde{\la}(q,dh;1-\be_j)=1$.

\medskip\noindent
{\it Case 3.}\,\,$q|d$, $(q,h_1)=1$.

\medskip\noindent

We also have $\tilde{\la}(q,dh;1-\be_j)=1$ and
$$
\xi_j(q^r,dh)=\kappa(q^r)-\frac{\kappa(q^{r-1})q}{q-1}=r+1-\frac{r}{1-u}+O(\al r\log q/q),
$$
so that
$$
\sum_r\frac{\chi(q^r)\xi_j(q^r;d,h)}{q^{rs}}=\frac{1}{(1-vu)^2}-1-\frac{vu}{(1-u)(1-vu)^2}+O(\al \log q/q).
$$
This yields (A.3). $\Box$

\medskip\noindent

To prove Lemma 15.2 and 15.3 we need several prerequisite results.
 In what follows assume $dl<P_2^2$,  $(dl,D)=1$ and $|s-1|<5\al$. We have
$$\aligned
\m_1(d,l;s)=&\prod_{q<D}
\frac{(1-q^{-s-\be_1})(1-q^{-s-\be_2})}{(1-q^{-s})(1-\chi(q)q^{-s})}\bigg(1+\tilde{\la}_1(q,d)\sum_r\frac{\xi_1(q^r;d,l)}{q^{rs}}\bigg)
\\&\times\big(1+O(D^{-c})\big).
\endaligned\eqno(A.4)$$
Note that 
$$
\tilde{\kappa}_1(q^{r-1},dq)=\kappa_1(q^{r-1}),
$$
so that
$$
\xi_1(q^r;d,l)=
\begin{cases}\tilde{\kappa}_1(q^r,d)-\frac{\chi(q)q}{q-1}\kappa_1(q^{r-1})&\quad\text{if}\quad (q,l)=1\\
\tilde{\kappa}_1(q^r,d)&\quad\text{if}\quad q|l
\end{cases}
\eqno(A.5)
$$
for any $q$, $d$, $l$ and $r$.

{\it Proof of Lemma 15.2.}\,\, By (A.4), (A.5) and the relations
$$ \tilde{\la}_1(q^r,1)=\la_1(q),\quad \frac{(1-q^{-s-\be_1})(1-q^{-s-\be_2})}{(1-q^{-s})(1-\chi(q)q^{-s})}=\frac{1-q^{-1}}{1-\chi(q)q^{-1}}+O\bigg(\frac{\al\log q}{q}\bigg),
$$
it suffices to show that
$$\aligned
1+\la_1(q)\sum_r\frac{1}{q^{rs}}\bigg(\tk_1(q^r,1)-\frac{\chi(q)q}{q-1}\kappa_1(q^{r-1})\bigg)
=&\frac{1-\chi(q)q^{-1}}{1-q^{-1}}\frac{1-\chi(q)q^{-2}}{1-q^{-2}}\\&+O\bigg(\frac{\al\log q}{q}\bigg)
\endaligned\eqno(A.6)
$$
if $q<D$, $(q,D)=1$, and
$$\aligned
1+\la_1(q)\sum_r\frac{1}{q^{rs}}\bigg(\tk_1(q^r,1)-\frac{\chi(q)q}{q-1}\kappa_1(q^{r-1})\bigg)
=\frac{1}{1-q^{-1}}+O\bigg(\frac{\al\log q}{q}\bigg)
\endaligned\eqno(A.7)
$$
if $q|D$.

 Assume  $q<D$. Since
$$
\kappa_1(q^r)=\sum_{\mu_1+\mu_2=r}q^{-\mu_1\be_1-\mu_2\be_2}-\sum_{\mu_1+\mu_2=r-1}q^{-\mu_1\be_1-\mu_2\be_2},
$$
where $\mu_1$ and $\mu_2$ run through non-negative integers, it follows that
$$
\kappa_1(q^r)=1+O(\al r\log q)
$$
 if $q^r<P$. Hence
$$
\tk_1(q^r,1)=\sum_{h\in\n(q)}\frac{\kappa_1(q^rh)\chi(h)}{h}=\sum_{h\in\n(q)}\frac{\chi(h)}{h}+O(\al r\log q)=\frac{1}{1-vu}+O(\al r\log q).
$$
It follows that
$$
\sum_r\frac{1}{q^{rs}}\bigg(\tk_1(q^r,1)--\frac{\chi(q)q}{q-1}\kappa_1(q^{r-1})\bigg)=\frac{u}{1-u}\bigg(\frac{1}{1-vu}-\frac{v}{1-u}\bigg)+O(\al\log q/q).
$$
This together with the relation
$$
\la_1(q)=1-vu+O(\al\log q/q)
$$
yields
$$
1+\la_1(q)\sum_r\frac{1}{q^{rs}}\bigg(\tk_1(q^r,1)-\frac{\chi(q)q}{q-1}\kappa_1(q^{r-1})\bigg)=\frac{1}{1-u}-\frac{uv(1-vu)}{(1-u)^2}+O(\al\log q/q).
$$
It is direct to verify that
$$
\frac{1}{1-u}-\frac{uv(1-vu)}{(1-u)^2}=\frac{1-\chi(q)q^{-1}}{1-q^{-1}}\frac{1-\chi(q)q^{-2}}{1-q^{-2}}
$$
if $\chi(q)=\pm 1$, and
$$
\frac{1}{1-u}-\frac{uv(1-vu)}{(1-u)^2}=\frac{1}{1-q^{-1}}
$$
if $\chi(q)=0$, whence (A.6) and (A.7) follow. $\Box$

\medskip\noindent

{\it Proof of Lemma 15.3.}\,\, We give a sketch only, as the situation is similar to Lemma 15.2. it can be verified that the factor
$\u_{1j}(q,s)$ in the Euler product representation
$$
\u_{1j}(s)=\prod_q \u_{1j}(q,s)
$$
satisfies $\u_{1j}(q,s)=1+O(q^{-2\sigma})$ if $(q,D)=1$. Further, in the case $q\le D$ we have
$$
\u_{1j}(q,1)=(1-q^{-1})^2+O(\al_1/q)\quad\text{if}\quad q|D,
$$
and
$$
\u_{1j}(q,1)=\frac{(1-q^{-2})^2}{1-\chi(q)q^{-2}}+O(\al_1/q)\quad\text{if}\quad (q,D)=1.
$$
This completes the proof. $\Box$

\medskip\noindent

{\it Proof of Lemma 16.1.}\,\, For any $q$, $r$, $d$ and $l$ we have
$$
\xi_2(q^r;d,l)=
\begin{cases}
\tk_2(q^r;d)-\frac{\chi(q)q}{q-1}\kappa_2(q^{r-1})&\quad\text{if}\quad (q,l)=1\\
\tk_2(q^r;d)&\quad\text{if}\quad q|l,
\end{cases}
$$
$$
|\kappa_2(q^r)|=|q^{-\be_1}-1|\ll\al\log q,\qquad \tilde{\la}_2(q,d)=1+O(\al\log q/q).
$$
Hence
$$
1+\tilde{\la}_2(q,d)\sum_r\frac{\xi_2(q^r;d,l)}{q^{rs}}=1-\frac{\chi(q)}{q-1}+O(\al\log q/q)\quad\text{if}\quad (q,l)=1
$$
and
$$
1+\tilde{\la}_2(q,d)\sum_r\frac{\xi_2(q^r;d,l)}{q^{rs}}=1+O(\al\log q/q)\quad\text{if}\quad q|l.
$$
On the other hand we have
$$
\frac{1-q^{-s-\be_1}}{(1-q^{-s})(1-\chi(q)q^{-s})}=\frac{1}{1-\chi(q)q^{-1}}+O(\al\log q/q).
$$
It follows that
$$
\m_2(d,l;s)=\prod_{\substack{q<D\\(q,l)=1}}\frac{1}{1-\chi(q)q^{-1}}\bigg(1-\frac{\chi(q)}{q-1}\bigg)\prod_{\substack{q<D\\q|l}}\frac{1}{1-\chi(q)q^{-1}}+O(\al_1).
$$
It is direct to verify that in either case the assertion holds. $\Box$

\medskip\noindent

{\it Proof of Lemma 16.2.}\,\,We give a sketch only. If $dl<D$, $(dl,D)=1$ and $|s-1|\le 5\al$, then
$$
\frac{\m_2(d,l;s)}{\m_2^*(s)}=\Pi_2(l)+O(\al_1\tau_2(dl))
$$
with
$$
\Pi_2(l)=\prod_{q|l}\bigg(1-\frac{\chi(q)}{q-1}\bigg)^{-1}\quad\text{if}\quad \chi(2)\ne 1,
$$
$$
\Pi_2(l)=\prod_{\substack{q|l\\q>2}}\bigg(1-\frac{\chi(q)}{q-1}\bigg)^{-1}\quad\text{if}\quad \chi(2)=1,\,\,2|l,
$$
$$
\Pi_2(l)=0\quad\text{if}\quad \chi(2)=(l,2)=1.
$$
The assertion follows by discussing the cases $\chi(2)\ne 1$ and $\chi(2)=1$ respectively. $\Box$

\medskip\noindent

 \section  *{\centerline{Appendix B. Some arithmetic sums}}
\medskip\noindent

{\it Proof of Lemma 15.1.}\,\, Put
$$
\varrho_j(n)=\sum_{d|n}\mu(d)d^{\be_j}.
$$
First we claim that
$$
\sum_{\substack{n<P\\(n,\q)=1}}\frac{|\varrho_j(n)-\varrho^*_j(n)|}{n}\ll\l^{-8}. \eqno(B.1)
$$
Since $\chi=\mu*\nu$, it follows that
$$
\chi(d)=\mu(d)+O\bigg(\sum_{\substack{h|d\\h>1}}\nu(h)\bigg).
$$
Hence
$$
\varrho_j(n)-\varrho^*_j(n)\ll\sum_{\substack{h|n\\h>1}}\nu(h)\tau_2(n/h).
$$
If $h>1$ and $(h,\q)=1$, then $h>D^4$. Hence, by substituting $n=hm$,
$$
\sum_{\substack{n<P\\(n,\q)=1}}\frac{|\varrho_j(n)-\varrho^*_j(n)|}{n}\ll\sum_{h>D^4}\frac{\nu(h)}{h}\sum_{m<P/h}\frac{\tau_2(m)}{m}\ll(\log P)^2\sum_{D^4<h<P}\frac{\nu(h)}{h}.
$$
This together with Lemma 3.2 yields (B.1).

Next we claim that
$$
\sum_{\substack{n<P\\(n,\q)>1}}\frac{|\varrho_j(n)|}{n}\ll\l^{-8}.\eqno(B.2)
$$
Note that $\varrho_j(q^r)=\varrho_j(q)$ for any $r$, and $\varrho_j(q)\ll\al\log q$ if $q<P$. Thus the left side above is
$$
\le\bigg(\sum_{q<D^4}\frac{|\varrho_j(q)|}{q}\bigg)\prod_{q<P}\big(1+|\varrho_j(q)|/q+O(1/q^2)\big)
\ll\al\sum_{q<D^4}\frac{\log q}{q}.
$$
This yields (B.2).

By (B.1) and (B.2), for $\mu=2, 3$,
$$
\sum_{(l,\q)=1}\frac{\vk_{\mu}(l_1l)\varrho^*_j(l)}{l}=\sum_{l}\frac{\vk_{\mu}(l_1l)\varrho_j(l)}{l}+O(\l^{-8}).
$$

We proceed  to prove theassertion with $\mu=2$. Since
$$
\sum_l\frac{\varrho_j(l)}{l^s}=\frac{\zeta(s)}{\zeta(s-\be_j)}
$$
for $\sigma>1$ and 
$$
\vk_2(m)=\tp\int_{(1)}\frac{(P_2/m)^s}{(\log P_2)(s-\be_7)^2}\,ds,
$$
it follows that
$$
\sum_{l}\frac{\vk_{2}(l_1l)\varrho_j(l)}{l}=\tp\int_{(1)}\frac{\zeta(1+s)}{\zeta(1+s-\be_j)}\frac{(P_2/l_1)^s}{(\log P_2)(s-\be_7)^2}\,ds.
$$
In a way similar to the proof of Lemma 8.1, we see that the right side is equal to the sum of the residues of the integrand at $s=0$ and $s=\be_7$ plus an acceptable error. By direct calculation, the residue at $s=0$ is
$$
-\frac{\be_j}{\be_7^2\log P_2}=-\frac{8j}{25\pi i}+O(\al_1);
$$
the residue at $s=\be_7$ is, by the Cauchy integral formula,
$$\aligned
\frac{\zeta(1+\be_7)}{\zeta(1+\be_7-\be_j)}&\frac{\log(P_2/l_1)}{\log P_2}\bigg(\frac{P_2}{l_1}\bigg)^{\be_7}+\frac{1}{\log P_2}\bigg(\frac{P_2}{l_1}\bigg)^{\be_7}\frac{d}{ds}\frac{\zeta(1+s)}{\zeta(1+s-\be_j)}|_{s=\be_7}\\
=&\bigg(1-\frac{2j}{5}+\frac{8j}{25\pi i}\bigg)\exp\{5\pi i/4\}+O(\al_1).
\endaligned
$$
These together  complete the proof in case $\mu=2$.
In case  $\mu=3$ the proof  can be obtained with $\be_6$ and $P_3$ in place of $\be_7$ and $P_2$ respectively. The same argument also gives
$$
\sum_{l}\frac{\vk_{1}(l_1l)\varrho_j(l)}{l}=\bigg(1-\frac{2j}{3}+\frac{j}{1.134\pi i}\bigg)\exp\{0.756\pi i\}
-\frac{j}{1.134\pi i}+O(\al_1).
$$
For $\mu=1$ the proof  is therefore reduced to showing that
$$
\sum_{l>P^{12}/l_1}\frac{\vk_{1}(l_1l)\varrho_j(l)}{l}=\frac{j}{0.756}\int_{0}^{0.004}\big(
\exp\{(3/2)(0.504-z)\pi i\}-\exp\{(3/4)\pi i\}\big)\,dz+O(\al_1).\eqno(B.3)
$$

By (4.2) and (4.3), the left side of (B.3) is equal to
$$
\sum_{l}\frac{\varrho_j(l)}{l}\bigg(\frac{P_1}{l_1l}\bigg)^{\be_6}\frac{1}{0.504}\int_{0.5}^{0.504}\bigg\{g\bigg(\frac{P^z}{l_1l}\bigg)-g\bigg(\frac{P^{0.5}}{l_1l}\bigg)\bigg\}\,dz+O(\al_1).
$$
By a change of variable, for $0.5\le z\le 0.504$,
$$
\bigg(\frac{P_1}{y}\bigg)^{\be_6}\bigg\{g\bigg(\frac{P^z}{y}\bigg)-g\bigg(\frac{P^{0.5}}{y}\bigg)\bigg\}=\tp\int_{(1)}\big(P^{\be_6(0.504-z)}P^{zs}-P^{0.004\be_6}P^{0.5s}\big)\frac{\omega_1(s-\be_6)\,ds}{y^{s}(s-\be_6)}.
$$
Hence, in a way similar to the proof of, we find that the left side of (B.3) is
$$\aligned
\frac{1}{0.504}&\int_{0.5}^{0.504}\,\bigg\{\tp\int_{(1)}\big(P^{\be_6(0.504-z)}P^{zs}-P^{0.004\be_6}P^{0.5s}\big)
\frac{\zeta(1+s)}{\zeta(1+s-\be_j)}\frac{\omega_1(s-\be_6)\,ds}{l_1^{s}(s-\be_6)}\bigg\}dz+O\\
=&\frac{1}{0.504}\frac{\be_j}{\be_6}\int_{0.5}^{0.504}\big(P^{\be_6(0.504-z)}-P^{0.004\be_6}\big)\,dz+O(\al_1).\quad\Box
\endaligned
$$

{\it Proof of Lemma 17.1.}\,\,By Lemma 3.1, 
$$
\sum_{n<D^4}\frac{\nu(n)^2}{n}=\sum_{n}\frac{\nu(n)^2}{n}g\bigg(\frac{T}{n}\bigg)+o(1).
$$
The sum on the right side is equal to
$$
\tp\int_{(1)}\bigg(\sum_{n}\frac{\nu(n)^2}{n^{1+s}}\bigg)\frac{T^s\omega_1(s)}{s}\,ds.
$$
Assume $\sigma>1$. We have
$$
\sum_{n}\frac{\nu(n)^2}{n^{s}}=\prod_p\bigg(1+\sum_r\frac{\nu(p^r)^2}{p^{rs}}\bigg).
$$
If $\chi(p)=1$, then (see [19, (1.2.10)])
$$
1+\sum_r\frac{\nu(p^r)^2}{p^{rs}}=(1-p^{-s})^{-4}(1-p^{-2s})=(1-p^{-s})^{-2}(1-\chi(p)p^{-s})^{-2}(1-p^{-2s});
$$
if $\chi(p)=-1$, then
$$
1+\sum_r\frac{\nu(p^r)^2}{p^{rs}}=1+\sum_r\frac{1}{p^{2rs}}=(1-p^{-2s})^{-1}=(1-p^{-s})^{-2}(1-\chi(p)p^{-s})^{-2}(1-p^{-2s});
$$
if $\chi(p)=0$, then
$$
1+\sum_r\frac{\nu(p^r)^2}{p^{rs}}=1+\sum_r\frac{1}{p^{rs}}=(1-p^{-s})^{-1}=(1-p^{-s})^{-2}(1-\chi(p)p^{-s})^{-2}(1-p^{-s}).
$$
Hence
$$
\sum_{n}\frac{\nu(n)^2}{n^{s}}=\zeta(s)^2L(s,\chi)^2\prod_{(p,D)=1}(1-p^{-2s})\prod_{p|D}(1-p^{-s}).
$$
In a way similar to the proof of, by (A) and simple estimate, we find that the integral (14) is equal to the residue of the function
$$
\zeta(1+s)^2L(1+s,\chi)^2\prod_{(p,D)=1}(1-p^{-2(1+s)})\prod_{p|D}(1-p^{-(1+s)})\frac{T^s\omega_1(s)}{s}
$$
at $s=0$, plus an acceptable error $O$, which is equal to
$$
L'(1,\chi)^2\prod_p(1-p^{-2})\prod_{p|D}(1-p^{-1})(1-p^{-2})^{-1}+o(1)=\mathfrak{a}+o(1).\quad\Box
$$

\begin{flushleft}
\medskip\noindent
\begin{tabbing}
XXXXXXXXXXXXXXXXXXXXXXXXXX*\=\kill
Yitang Zhang\\
Department of Mathematics\\
University of California at Santa Barbara\\
Santa Barbara, CA 93106\\

E-mail: yitang.zhang@math.ucsb.edu

\end{tabbing}

\end{flushleft}
\end {document}